\documentstyle{amsppt}
\pageheight{19cm}
\magnification=1200
\baselineskip=18pt
\nologo
\TagsOnRight
%\NoBlackBoxes
\document
\define \ts{\thinspace}
\define \oa{\frak o}
\define \aaa{\frak a}
\define \g{\frak g}
\define \spa{\frak {sp}}
\define \Ra{\Rightarrow}
\define \gl{\frak {gl}}
\define \sll{\frak {sl}}
\define \qdet{\operatorname {qdet}}
\define \sdet{\operatorname {sdet}}
\define \U{\operatorname {U}}
\define \Norm{\operatorname {Norm}}
\define \Y{\operatorname {Y}}
\define \SY{\operatorname {SY}}
\define \Sym{{\frak S}}
\define \ZZ{{\Bbb Z}}
\define \sgn{\text{{\rm sgn}}}

\define \dime{\text{\rm dim}\ts}
\define \ot{\otimes}
\define \hra{\hookrightarrow}
\define \wh{\widehat}
\define \wt{\widetilde}
\define \sss{S}
\define \sst{S^*}
\define \Z{\operatorname {Z}}

\define \End{\operatorname {End}}

\define \C{\Bbb C}

\define \Proof{\noindent {\it Proof. }}

\heading{\bf FINITE-DIMENSIONAL IRREDUCIBLE REPRESENTATIONS}
\endheading
\heading{\bf OF TWISTED YANGIANS}
\endheading 
\bigskip
\bigskip
\heading{Alexander Molev}\endheading
\bigskip
\bigskip
\bigskip
\bigskip
\noindent
Centre for Mathematics and its Applications\newline
Australian National University\newline
Canberra, ACT 0200, Australia\newline
(e-mail: molev\@pell.anu.edu.au)
\bigskip
\bigskip
\noindent
{\bf November 1997}
\bigskip
\bigskip
\noindent
{\bf Abstract}\newline
We study quantized enveloping algebras called twisted Yangians.
They are analogues of the Yangian $\Y(\gl(N))$ for the classical
Lie algebras of $B$, $C$, and $D$ series. The twisted Yangians are
subalgebras in $\Y(\gl(N))$ and coideals with respect to the coproduct
in $\Y(\gl(N))$. We give a complete description of their 
finite-dimensional irreducible representations. Every such
representation is highest weight and we give necessary and sufficient
conditions for an irreducible highest weight
representation to be finite-dimensional. The result is analogous to
Drinfeld's theorem for the ordinary Yangians. Its detailed proof
for the $A$ series is also reproduced.
For the simplest twisted Yangians we construct an explicit realization 
for each finite-dimensional irreducible representation
in tensor products of representations of the corresponding Lie
algebras. 


\newpage
\noindent
{\bf 0. Introduction}
\bigskip

The Yangian $\Y(\aaa)$ associated to a simple complex Lie algebra $\aaa$
is a deformation of the universal enveloping algebra $\U(\aaa[t])$ 
in the class of Hopf algebras [D1], [D2]. In the case of the $A$ series
one can also define the Yangian $\Y(N)=\Y(\gl(N))$ for the reductive
complex Lie algebra $\aaa=\gl(N)$ (see, e.g., [TF]). The Yangian
$\Y(\sll(N))$ can be identified with a (Hopf) subalgebra in $\Y(\gl(N))$.
The Yangians are
closely related to rational solutions of the Yang-Baxter equation
[D1] and play the role of quantum symmetry groups of certain
integrable systems; see, e.g., [BL], [KS1], [TU], [UK].
Finite-dimensional irreducible
representations of the Yangians were classified by V.~G.~Drinfeld [D2].
In the case of $\aaa=\gl(2)$ this had been done by V.~O.~Tarasov
[T1], [T2]; see also [K].\footnote{In these papers the language of
the quantum inverse scattering theory and integrable models was used.
Representations of
the Yangian $\Y(2)$ appear there in the form of the {\it monodromy
matrix\/}
for the XXX-model.}
This result plays the main role in the proof
in the general case. Representations of the Yangians were further studied
in [CP1--CP3], [C1], [C2], [KR], [M2], [NT2].

For any $\aaa$ the Yangian $\Y(\aaa)$ contains the universal
enveloping algebra $\U(\aaa)$. However,
the $A$ series case seems to be exceptional since only in this case does
there exists an algebra homomorphism $\Y(\aaa)\to\U(\aaa)$ identical on
$\U(\aaa)$. This homomorphism is important for applications of
the `quantum theory' to representations of the
Lie algebra $\gl(N)$ such as
constructions of Laplace operators and characteristic identities 
[C2], [M4], [NT1]; Capelli identities [N1], [N2], [Ok]. 

Furthermore, as was shown by G. Olshanski [O1],
the Yangian $\Y(N)$
can be constructed as a projective limit of 
$\gl(M)$-centralizers 
in $\U(\gl(N+M))$ as $M\to\infty$.
It turned out that Olshanski's construction applied to other series of
classical Lie algebras leads not to the corresponding Yangian 
$\Y(\aaa)$, but to
different algebras which were called twisted Yangians [O2]; see also [S].
Defining relations for the twisted Yangian can be written in the form of
a reflection-type equation (see (1.16) below) which allows one to use
a special $R$-matrix technique for their study. Some other algebras
associated to the reflection equation were studied in [KK], [KS2], [KJC].

Let $\g(N)$ denote the orthogonal Lie algebra $\oa(N)$ or symplectic Lie
algebra $\spa(N)$. We regard $\g(N)$ as a fixed point subalgebra in 
$\gl(N)$ with respect to an involutive automorphism 
$\sigma$. The twisted Yangian is a deformation of the universal enveloping
algebra $\U(\gl(N)[t]^{\sigma})$, where
$$
\gl(N)[t]^{\sigma}=\{A(t)\in\gl(N)[t]\ |\ \sigma(A(t))=A(-t)\}.
$$
We denote by $\Y^+(N)$ and $\Y^-(N)$ the twisted Yangians corresponding
to the orthogonal and symplectic Lie algebras, respectively.
Similarly to the $A$ series case the algebra $\Y^{\pm}(N)$
contains the corresponding universal enveloping algebra $\U(\g(N))$
and one has a natural algebra homomorphism $\Y^{\pm}(N)\to \U(\g(N))$
identical on $\U(\g(N))$. 

A detailed description of the algebraic structure of the Yangian $\Y(N)$ 
and the twisted Yangians $\Y^{\pm}(N)$ is contained in [MNO].
An explicit formula for the formal series (the Sklyanin determinant)
whose coefficients generate the center
of $\Y^{\pm}(N)$ was given in [M3]. Bethe-type
commutative subalgebras in $\Y^{\pm}(N)$ were constructed in [NO]. Some
applications of these results
to the Laplace operators for the orthogonal and symplectic
Lie algebras, characteristic identities and
Capelli identities can be found in [M3], [M4], [MN].

In this paper we give a complete description of the
finite-dimensional irreducible representations for all twisted Yangians.
Our approach is very close to the $\gl(N)$-case and we use some 
ideas and results in that case in our constructions.

Section 1 contains definitions and preliminary results on the Yangian
$\Y(N)$ and the twisted Yangians $\Y^{\pm}(N)$.

In Section 2 we give detailed proofs of Drinfeld's and Tarasov's theorems
describing the finite-dimensional irreducible representations
of $\Y(N)$. These representations can be constructed as follows. First,
using
the homomorphism $\Y(N)\to\U(\gl(N))$ one extends any
finite-dimensional irreducible representation $L(\lambda)$
of the Lie algebra $\gl(N)$ to the Yangian $\Y(N)$. Using the Hopf algebra
structure on $\Y(N)$ one can consider representations of the form
$$
L(\lambda^{(1)})\ot\cdots\ot L(\lambda^{(k)}). \tag 0.1
$$
Then any
finite-dimensional irreducible representation of $\Y(\sll(N))$
(or $\Y(N)$, up to tensoring with one-dimensional 
representations of the center of $\Y(N)$)
can be realized
as the irreducible quotient of the cyclic span in (0.1) 
of the tensor product of the highest weight vectors of the 
$L(\lambda^{(i)})$.

In the case $N=2$ irreducible representations (0.1)
constitute a complete list of  
finite-dimensional irreducible representation
of $\Y(\sll(2))$. A criterion of irreducibility of representation (0.1)
of $\Y(\sll(2))$
in terms of `strings' (subsets in $\C$ of the form $\{a,a+1,\dots,b\}$)
was given by Chari and Pressley [CP1], [CP3]. We give a proof
of this result below (Corollary 2.11).

Alternatively, every finite-dimensional irreducible representation
of $\Y(N)$ is highest weight (see Definition 2.1). 
It is isomorphic to
a representation of the form $L(\lambda(u))$, where 
$\lambda(u)=(\lambda_1(u),\dots,\lambda_N(u))$ is
a collection of formal series in $u^{-1}$. Moreover, given a set
$\lambda(u)$
there exists a unique (up to an isomorphism) irreducible
highest weight module $L(\lambda(u))$. Drinfeld's theorem [D2] 
(see Theorem 2.12 below) provides
necessary and sufficient conditions for $L(\lambda(u))$ to be
finite-dimensional.

In Sections 3--6 we give analogues of these descriptions for the
twisted Yangians. Namely, using the homomorphism 
$\Y^{\pm}(N)\to \U(\g(N))$ we extend any 
finite-dimensional irreducible representation $V(\mu)$
of the Lie algebra $\g(N)$ to the twisted Yangian $\Y^{\pm}(N)$.
The subalgebra $\Y^{\pm}(N)\subset \Y(N)$ is a left coideal,
$$
\Delta(\Y^{\pm}(N))\subset \Y(N)\ot \Y^{\pm}(N), \tag 0.2
$$
and so,
we may construct representations of $\Y^{\pm}(N)$ of the form
$$
L(\lambda^{(1)})\ot\cdots\ot L(\lambda^{(k)})\ot V(\mu) . \tag 0.3
$$
As in the $\gl(N)$-case, any
finite-dimensional irreducible representation of $\Y^{\pm}(N)$
(up to tensoring with one-dimensional representations of the center) 
can be realized
as the irreducible quotient of the cyclic span in (0.3) 
of the tensor product of the highest weight vectors of the 
$L(\lambda^{(i)})$ and $V(\mu)$.

On the other hand, we prove that any
finite-dimensional irreducible representations of $\Y^{\pm}(N)$
is highest weight, that is, isomorphic to a representation
of the form $V(\mu(u))$ (see Definition 3.1), and
we find necessary and sufficient 
conditions for $V(\mu(u))$ to be finite-dimensional; see
Theorems 4.8, 5.9, 6.7 for the cases of $C$, $D$, $B$ series,
respectively.

The key part of the proof is a detailed study of the simplest cases
$\Y^{\pm}(2)$ and $\Y^{+}(3)$.
We use a modified Tarasov's argument
[T1], [T2] and also use some ideas from [CP1], [CP3], as well as
[Z1], [Z2] where representations of
the super-Yangians were studied.

We give explicit realizations
for all
finite-dimensional irreducible representations of the simplest
twisted Yangians $\Y^{-}(2)$ and $\Y^{+}(2)$ and prove
analogues of Chari--Pressley's theorem for these algebras;
see Corollaries 4.5, 4.7 and Corollaries 5.5, 5.6, respectively.

The $B$ series case requires an extra care. We construct raising
and lowering operators in a $\Y^{+}(3)$-module $V(\mu(u))$ which
preserve the subspace of
$\Y^{+}(2)$-singular vectors. 
To get the necessary
conditions for $V(\mu(u))$ to be finite-dimensional
we examine the weights of the
$\Y^{+}(2)$-singular vectors and use the results for the case of
$\Y^{+}(2)$.

Some of the results of this paper (mainly for the case of $\Y^-(2n)$)
were announced without complete proofs in [M1].
\medskip

During this work I benefited from collaboration with G.~Olshanski
to whom I would like to express my sincere gratitude.
I would also like to thank M.~Nazarov, V.~Tarasov and R.~Zhang
for valuable discussions.

\bigskip
\bigskip
\noindent
{\bf 1. Definitions and preliminary results}
\bigskip

Our main reference for the algebraic structure of the Yangian $\Y(N)$
and the twisted Yangians $\Y^{\pm}(N)$ is the expository paper [MNO].
In this section we give necessary definitions and formulate some
results to be used in the subsequent sections.
\medskip

We shall assume that given positive integer $N$, indices $i,j$ run
through the set $\{-n,\dots,-1,0,1,\dots n\}$ if $N=2n+1$ and
run through the set $\{-n,\dots,-1,1,\dots n\}$ if $N=2n$.

The {\it Yangian\/} $\Y(N)=\Y(\gl(N))$ is the
complex associative algebra with the
generators $t_{ij}^{(1)},t_{ij}^{(2)},\dots$ where $-n\leq i,j\leq n$,
and the defining relations
$$
[t_{ij}^{(r+1)},t_{kl}^{(s)}]-
[t_{ij}^{(r)},t_{kl}^{(s+1)}]=
t_{kj}^{(r)}t_{il}^{(s)}-t_{kj}^{(s)}t_{il}^{(r)},
$$
where $r,s=0,1,2,\ldots\;$ and $t_{ij}^{(0)}:=\delta_{ij}$. 
For any $i,j$ define the formal power series
$$
t_{ij}(u) = \delta_{ij} + t^{(1)}_{ij} u^{-1} + t^{(2)}_{ij}u^{-2} +
\cdots \in \Y(N)[[u^{-1}]].
$$
Then the defining relations can be also written as follows:
$$
[t_{ij}(u),t_{kl}(v)]={1\over u-v}(t_{kj}(u)t_{il}(v)-t_{kj}(v)t_{il}(u)).
\tag 1.1
$$
They also can be written as a single 
{\it ternary relation\/} (see, e.g., [MNO, Section 1])
for the $T$-{\it matrix\/}
$$
T(u):=\sum_{i,j} t_{ij}(u)\ot E_{ij}
\in \Y(N)[[u^{-1}]]\ot \End\C^N,
$$
where the $E_{ij}$ are the standard matrix units. It has the form
$$
R(u-v)T_1(u)T_2 (v) = T_2(v)T_1(u)R(u-v). \tag 1.2
$$
Here both sides are regarded as elements of
$\Y(N)[[u^{-1}]]\ot \End\C^N\ot \End\C^N$, the subindex of $T(u)$
indicates to which copy of $\End\C^N$ this matrix corresponds;
$R(u)=1-u^{-1}P$, where $P$ is the permutation operator in
$\C^N\ot \C^N$.

The following analogue of the Poincar\'e--Birkhoff--Witt theorem
holds for the algebra $\Y(N)$ [MNO, Corollary 1.23]:
\medskip

given an arbitrary linear order on the set of
the generators
$t^{(r)}_{ij}$, any element of the algebra $\Y(N)$ is uniquely written
as a linear combination of ordered monomials in the generators.
\medskip
\noindent
In what follows we shall  always assume that
a linear order $\prec$ is chosen in such a way that
$$
t^{(r_1)}_{i_1j_1}\prec t^{(r_2)}_{i_2j_2}\prec t^{(r_3)}_{i_3j_3}
\tag 1.3
$$
as soon as $i_1>j_1$,  $i_2=j_2$, and $i_3<j_3$.

We shall regard the set of the matrix units $\{E_{ij}\}$ as a basis of the
Lie
algebra $\gl(N)$. The mapping
$$
E_{ij}\to t_{ij}^{(1)} \tag 1.4
$$
defines an inclusion $\U(\gl(N))\hra\Y(N)$.

The mapping
$$
t_{ij}(u)\mapsto \delta_{ij}+E_{ij}u^{-1} \tag 1.5
$$
defines an algebra homomorphism
$
\Y(N)\to \U(\gl(N)).
$

The {\it quantum determinant\/} $\qdet T(u)$
of the matrix $T(u)$ is a formal series in $u^{-1}$
with coefficients from $\Y(N)$ defined as follows: 
$$
\qdet T(u)=\sum_{p\in \Sym_N} \sgn(p)\ts t_{p(-n),-n}(u)\cdots
t_{p(n),n}(u-N+1),
$$
where $\Sym_N$ is the group of permutations of the indices
$\{-n,\dots,n\}$. 

The coefficients of the quantum
determinant $\qdet T(u)$ are algebraically independent generators of
the center of the algebra $\Y(N)$.

Let us define the {\it quantum comatrix\/}
$\widehat T(u)=(\widehat t_{ij}(u))$ by the
following formula:
$$
\qdet T(u)=\widehat T(u)T(u-N+1).\tag 1.6
$$
Then $\widehat t_{ij}(u)$ equals $(-1)^{i+j}$ times the quantum determinant
of the submatrix of $T(u)$ obtained by removing the 
$i$th column and $j$th row.

Each of the following mappings defines an antiautomorphism of the
algebra $\Y(N)$: change of sign in $u$
$$
\sigma:\ts T(u)\to T(-u);\tag 1.7
$$
a matrix transposition
$$
T(u)\to T^t(u);\tag 1.8
$$
comatrix transformation
$$
T(u)\to \widehat T(u). \tag 1.9
$$

The mapping
$$
T(u)\to T(u+a),\quad a\in \C; \tag 1.10
$$
and
multiplication by a formal series $\varphi(u)=1+\varphi^{(1)}u^{-1}+
\varphi^{(2)}u^{-2}+\cdots$:
$$
T(u)\to \varphi(u)T(u); \tag 1.11
$$
define automorphisms of $\Y(N)$.

The {\it Yangian\/} $\Y(\sll(N))$ for the Lie algebra $\sll(N)$ can
be defined as the subalgebra of $\Y(N)$ consisting of the elements fixed by
all automorphisms of the form (1.11). One has a 
tensor product decomposition
$$
\Y(N)\simeq \Z(N)\ot  \Y(\sll(N)), \tag 1.12
$$
where $\Z(N)$ denotes the center of $\Y(N)$.

Using the automorphism of $\Y(N)$ of the
form $T(u)\to T^t(-u)$, one can show that
$$
\qdet T(u)=\widehat T^t(u-1)T^t(u).\tag 1.13
$$

Coproduct 
$$
\Delta : \Y(N)\to
\Y(N)\ot\Y(N) 
$$ 
is defined by the formula
$$
\Delta (t_{ij}(u)):=\sum_{a} t_{ia}(u)\otimes
t_{aj}(u).\tag 1.14
$$

\medskip

Let us now define
the {\it twisted Yangians\/} $\Y^{+}(N)$ and $\Y^{-}(N)$
corresponding to the orthogonal Lie algebra $\oa(N)$ and
symplectic Lie algebra $\spa(N)$, respectively.
To consider both cases $\Y^{+}(N)$ and $\Y^{-}(N)$
simultaneously,
it will be convenient to use
the symbol $\theta_{ij}$ which is defined as follows:
\medskip
$$
\theta_{ij}:=\cases 1,\quad&\text{in the orthogonal case};\\
\sgn(i)\sgn(j),\quad&\text{in the symplectic case}.\endcases
$$
\bigskip
\noindent
Whenever the double sign $\pm{}$ or $\mp{}$ occurs,
the upper sign corresponds to the orthogonal case and the lower sign to
the symplectic one.
By $X\mapsto X^t$ we will denote the matrix 
transposition such that
$
(E_{ij})^t=\theta_{ij}E_{-j,-i}.
$
Let us introduce the $S$-{\it matrix\/} $S(u)=(s_{ij}(u))$ by setting
$
S(u):=T(u)T^t(-u),
$
or, in terms of matrix elements,
$$
s_{ij}(u)=\sum_a \theta_{aj}t_{ia}(u)t_{-j,-a}(-u). \tag 1.15
$$
Write
$$
s_{ij}(u)=\delta_{ij}+s_{ij}^{(1)}u^{-1}+s_{ij}^{(2)}u^{-2}+\cdots.
$$
The twisted Yangian $\Y^{\pm}(N)$ is the subalgebra of $\Y(N)$ generated by
the elements $s_{ij}^{(1)},s_{ij}^{(2)},\dots$, where $-n\leq i,j\leq n$.

One can show that the $S$-matrix satisfies the following
{\it quaternary relation\/} and {\it symmetry relation\/} 
(see [MNO, Section 3]):
$$
\aligned
R(u-v)S_1(u)R^t(-u-v)S_2(v)={}&S_2(v)R^t(-u-v)S_1(u)R(u-v), \\
S^t(-u)={}&S(u)\pm {S(u)-S(-u)\over 2u}.
\endaligned
\tag 1.16
$$
Here we use the same notation as for the ternary relation (1.2),
where $R^t(u)$ is obtained from $R(u)$ by applying the transposition $t$
in either of the two copies of $\End\C^N$.
 
Relations (1.16) are defining relations for the algebra $\Y^{\pm}(N)$ and
they can be rewritten in terms of 
the generating series $s_{ij}(u)$ as follows:
$$
\aligned
[s_{ij}(u),s_{kl}(v)]={}&{1\over
u-v}(s_{kj}(u)s_{il}(v)-s_{kj}(v)s_{il}(u))\\
-{}&{1\over u+v}(\theta_{k,-j}s_{i,-k}(u)s_{-j,l}(v)-
\theta_{i,-l}s_{k,-i}(v)s_{-l,j}(u))\\
+{}&{1\over u^2-v^2}(\theta_{i,-j}s_{k,-i}(u)s_{-j,l}(v)-
\theta_{i,-j}s_{k,-i}(v)s_{-j,l}(u)) 
\endaligned
\tag 1.17
$$
and
$$
\theta_{ij}s_{-j,-i}(-u)=s_{ij}(u)\pm {s_{ij}(u)-s_{ij}(-u)\over 2u}. 
\tag 1.18
$$

The
elements
$$
s_{ij}^{(2k)},\quad i+j\geq 0;\qquad  s_{ij}^{(2k-1)},\quad
i+j>0;\qquad k=1,2,\dots, \tag 1.19
$$
in the case of $\Y^{+}(N)$, and the elements
$$
s_{ij}^{(2k)},\quad i+j> 0;\qquad  s_{ij}^{(2k-1)},\quad
i+j\geq 0;\qquad k=1,2,\dots, \tag 1.20
$$
in the case of $\Y^{-}(N)$, 
constitute a system of linearly independent generators.
An analogue of
the Poincar\'e--Birkhoff--Witt theorem for the algebra
$\Y^{\pm}(N)$ can be formulated as follows:
\medskip

given an arbitrary linear order on the set of
the above generators,
any element of the algebra $\Y^{\pm}(N)$ is uniquely written
as a linear combination of ordered monomials in the generators;
\medskip
\noindent
see [MNO, Remark 3.14]. Similarly to the case of $\Y(N)$ we
shall always assume that a linear order ${}\prec$ is chosen
in such a way that
$$
s^{(r_1)}_{i_1j_1}\prec s^{(r_2)}_{i_2j_2}\prec s^{(r_3)}_{i_3j_3}
\tag 1.21
$$
as soon as $i_1>j_1$,  $i_2=j_2$, and $i_3<j_3$.

Introduce the following elements of the Lie algebra $\gl(N)$:
$$
F_{ij}=E_{ij}-\theta_{ij}E_{-j,-i},\qquad -n\leq i,j\leq n.
$$
Denote by $\g(N)$ the Lie subalgebra of $\gl(N)$ spanned 
by the elements $F_{ij}$.
Then $\g(N)$ is isomorphic to $\oa(2n)$ or
$\spa(2n)$ if $N=2n$, and to $\oa(2n+1)$ if $N=2n+1$. 

The mapping
$$
F_{ij}\to s_{ij}^{(1)}\tag 1.22
$$
defines an inclusion $\U(\g(N))\hra\Y^{\pm}(N)$.

The mapping
$$
s_{ij}(u)\mapsto \delta_{ij}+F_{ij}
\left(u\pm {1\over 2}\right)^{-1} \tag 1.23
$$
defines an algebra homomorphism
$
\Y^{\pm}(N)\to \U(\g(N)).
$

There is an analogue of the quantum determinant for the
twisted Yangians. It is denoted by $\sdet S(u)$ and is called the
{\it Sklyanin determinant\/}. This is a formal series in $u^{-1}$
and its coefficients generate the center of the algebra
$\Y^{\pm}(N)$. The Sklyanin determinant is related
to the quantum determinant by the formula:
$$
\sdet S(u)=\gamma_N(u)\ts\qdet T(u)\ts\qdet T(-u+N-1),\tag 1.24
$$
where $\gamma_N(u)\equiv1$ for $\Y^+(N)$ and $\gamma_N(u)=
{(2u+1)}/{(2u-N+1)}$ for $\Y^-(N)$. 
An explicit determinant-type expression for
$\sdet S(u)$ in terms of the generators $s_{ij}(u)$ was given in [M3].

Define the {\it Sklyanin comatrix\/}
$\widehat S(u)=(\widehat s_{ij}(u))$ by the
following formula:
$$
\sdet S(u)=\widehat S(u)S(u-N+1).\tag 1.25
$$
An explicit expression for $\widehat s_{ij}(u)$ can also be found from 
[M3, Section 6].

The restriction of the antiautomorphism (1.7)
to the subalgebra $\Y^{\pm}(N)$
gives an antiautomorphism
$$
\sigma:\ S(u)\to S^t(u). \tag 1.26
$$

We shall need the following result.

\proclaim
{\bf Proposition 1.1} The mapping
$$
S(u)\mapsto \gamma_N(u)\ts \wh S(-u+\frac N2-1),\tag 1.27
$$
defines an automorphism of the algebra $\Y^{\pm}(N)$.
\endproclaim

\Proof Note that the mapping $T(u)\to \wt T(u)$ with
$$
\wt T(u)=\wh T^t(u+\frac N2-1)=\qdet T(u+\frac N2)
(T^t(u+\frac N2))^{-1}
$$
defines an automorphism of the algebra $\Y(N)$. To see this,
it suffices to present this map as a composition
of antiautomorphisms of form (1.8), (1.9) and
an automorphism of form (1.10).
Denote by $\wt S(u)$ the image of the matrix $S(u)=T(u)T^t(-u)$
under this automorphism. That is,
$$
\wt S(u)=\wt T(u)\wt T^t(-u)=
\wh T^t(u+\frac N2-1)\wh T(-u+\frac N2-1).
$$
By (1.6) and (1.13) we have
$$
\wt S(-u+\frac N2-1)S(u-N+1)=\wh T^t(-u+N-2)\wh T(u)
T(u-N+1)T^t(-u+N-1)
$$
$$
=\qdet T(u)\ts\wh T^t(-u+N-2)T^t(-u+N-1)
=\qdet T(u)\ts\qdet T(-u+N-1),
$$
which equals 
$\gamma_N^{-1}(u)\ts\sdet S(u)$, by (1.24).
Note that
$$
\gamma_N(u)\gamma_N(-u+\frac N2-1)=1,
$$
and so, using (1.25) we obtain that
$$
\wt S(u)=\gamma_N(u)\ts \wh S(-u+\frac N2-1).
$$
Thus, the subalgebra $\Y^{\pm}(N)\subset\Y(N)$ is stable under
the automorphism $T(u)\to \wt T(u)$ and its restriction to
$\Y^{\pm}(N)$ yields the automorphism (1.27) which completes the proof.
$\square$
\medskip

We shall use the following commutation relations between
the matrix elements of the matrices $\wh S(u)$ and $S(u)$:
$$
\aligned
&[\wh s_{ij}(u),s_{kl}(v)]=-
{1\over u-v-N+1}
\sum_a\left(\delta_{kj}
\wh s_{ia}(u)s_{al}(v)-\delta_{il} s_{ka}(v)\wh s_{aj}(u)\right)\\ 
+&{1\over u+v+1}
\sum_a \left(\delta_{i,-k}\theta_{ak}
\wh s_{-a,j}(u)s_{al}(v)-\delta_{j,-l} \theta_{al}
s_{ka}(v)\wh s_{i,-a}(u)\right)\\ 
+&{1\over (u-v-N+1)(u+v+1)}
\sum_a \left(\delta_{i,-k}\theta_{jk}
\wh s_{-j,a}(u)s_{al}(v)-\delta_{j,-l}\theta_{il}
s_{ka}(v)\wh s_{a,-i}(u)\right). 
\endaligned
\tag 1.28
$$
They easily follow from the quaternary relation; see (1.16).
Indeed, inverting its both sides and multiplying 
from the left and right by $S_2(v)$ we obtain
$$
R^t(u+v-N)S^{-1}_1(u)R(v-u)S_2(v)=S_2(v)R(v-u)S^{-1}_1(u)R^t(u+v-N). 
$$
To get (1.28)
rewrite this first
in terms of the matrix elements of $S^{-1}(u)$ and $S(u)$, then use
(1.25) and the centrality of the Sklyanin determinant in $\Y^{\pm}(N)$.
(In a slightly different form this calculation was performed in
[MNO, Section 7]).

Multiplication by an even formal series $\psi(u)=1+\psi^{(2)}u^{-2}+
\psi^{(4)}u^{-4}+\cdots$ defines an automorphism of $\Y^{\pm}(N)$:
$$
S(u)\to \psi(u)S(u). \tag 1.29
$$

The {\it special twisted Yangian\/} $\SY^{\pm}(N)$ 
is defined by the equality $\SY^{\pm}(N)=\Y(\sll(N))\cap\Y^{\pm}(N)$.
In other words, $\SY^{\pm}(N)$ is the subalgebra of $\Y^{\pm}(N)$
consisting of elements which are stable under all automorphisms
of the form (1.29). One has the following tensor product decomposition
analogous to (1.12):
$$
\Y^{\pm}(N)\simeq \Z^{\pm}(N)\ot \SY^{\pm}(N), \tag 1.30
$$
where $\Z^{\pm}(N)$ denotes 
the center of $\Y^{\pm}(N)$.

The images of the generators $s_{ij}(u)$ under coproduct $\Delta$
are given by the formula
$$
\Delta(s_{ij}(u))=\sum_{a,b}\theta_{bj}t_{ia}(u)t_{-j,-b}(-u)\otimes
s_{ab}(u).
\tag 1.31
$$
This implies that the subalgebra $\Y^{\pm}(N)\subset\Y(N)$ is a left
coideal; 
see (0.2).

\newpage

\noindent
{\bf 2. Finite-dimensional irreducible representations of $\Y(N)$ and 
$\Y(\sll(N))$}
\bigskip

For this section only we adopt a more standard enumeration of the
generators $t_{ij}(u)$ of the Yangian $\Y(N)$: we let the indices $i,j$
run through the set $\{1,\dots,N\}$. All formulas from the previous section
should be understood accordingly. In the following definitions we follow
[D2] and [CP1--CP3].
\medskip

\noindent
{\bf Definition 2.1.} A representation $L$ of the Yangian $\Y(N)$
is called {\it highest weight\/} if there exists a nonzero vector 
$\xi\in L$ such that $L$ is generated by $\xi$,
$$
t_{ij}(u)\xi=0 \qquad \text{for} \quad 1\leq i<j\leq N \tag 2.1
$$
and
$$
t_{ii}(u)\xi=\lambda_i(u)\xi \qquad \text{for} \quad 1\leq i\leq N \tag 2.2
$$
for some formal series $\lambda_i(u)\in 1+u^{-1}\C[[u^{-1}]]$.
In this case the vector $\xi$ is called the {\it highest weight vector\/}
of $L$ and the set $\lambda(u):=(\lambda_1(u),\dots,\lambda_N(u))$
is the {\it highest weight\/} of $L$. $\square$
\medskip

Sometimes we shall also refer to
$L$ as a representation {\it with the highest weight\/} $\lambda(u)$.

The definition of a {\it lowest weight representation\/} is obtained by
replacing (2.1) with the relations
$$
t_{ji}(u)\xi=0 \qquad \text{for} \quad 1\leq i<j\leq N.
$$

\bigskip

\noindent
{\bf Definition 2.2.} Let $\lambda(u)=(\lambda_1(u),\dots,\lambda_N(u))$
be any set of formal series where $\lambda_i(u)\in 1+u^{-1}\C[[u^{-1}]]$.
The {\it Verma module\/} $M(\lambda(u))$ is the quotient of $\Y(N)$
by the left ideal generated by all the coefficients of the series
$t_{ij}(u)$ for $1\leq i<j\leq N$ and $t_{ii}(u)-\lambda_i(u)$ for 
$1\leq i\leq N$. $\square$
\medskip

The following properties of $M(\lambda(u))$ are
immediate from
the Poincar\'e--Birkhoff--Witt theorem
for $\Y(N)$. First,
ordered monomials in the coefficients of the series
$t_{ij}(u)$ with $i>j$ (assuming that (1.3) is satisfied) 
form a basis of 
$M(\lambda(u))$. The Verma module $M(\lambda(u))$ is obviously
a representation of $\Y(N)$ with 
the highest weight $\lambda(u)$ and any representation with
the highest weight $\lambda(u)$ is isomorphic to a quotient
of $M(\lambda(u))$. Further, $M(\lambda(u))$ has a unique irreducible
quotient which we denote by $L(\lambda(u))$. 
Indeed, the sum of all submodules of $M(\lambda(u))$ which do not contain
the highest weight vector $1$ is a unique maximal proper submodule $K$,
and $L(\lambda(u))$ is the quotient $M(\lambda(u))/K$.
Moreover, if an irreducible $\Y(N)$-module $L$ contains a nonzero
vector $\xi$ satisfying (2.1) and (2.2) then $L$ is isomorphic to
$L(\lambda(u))$. 

Next, using the inclusion (1.4) we may regard any $\Y(N)$-module
as a $\gl(N)$-module. Relations (1.1)
imply that
$$
[E_{ij},t_{kl}(u)]=\delta_{kj}t_{il}(u)-\delta_{il}t_{kj}(u), \tag 2.3
$$
where we identify the elements $E_{ij}\in\gl(N)$ with their
images in $\Y(N)$.
Let $\frak h$ be the diagonal subalgebra in $\gl(N)$ with the basis
$E_{11},\dots, E_{NN}$ and let $\varepsilon_1,\dots, \varepsilon_N$
be the dual basis in $\frak h^*$. 
Choose the standard system of positive roots $R_+$ for 
$\gl(N)$ corresponding to
the upper triangular Borel subalgebra. Then the elements
$\alpha_1=\varepsilon_1-\varepsilon_2$,$\dots$,$\alpha_{N-1}=
\varepsilon_{N-1}-\varepsilon_N$ form
a basis in $R_+$.
If $L$ is a representation of $\Y(N)$ with the highest weight
$\lambda(u)$, and $\xi$ is its highest weight vector, then $\xi$ is also
a $\gl(N)$-highest weight vector of the $\gl(N)$-weight
$\lambda=(\lambda_1,\dots,\lambda_N)$ where $\lambda_i$ is the coefficient
at $u^{-1}$ in the series $\lambda_i(u)$. By (2.3) 
all $\gl(N)$-weights of $L$ have the form 
$\lambda-\omega$, where $\omega$ is a linear combination of the $\alpha_i$
with coefficients from $\ZZ_+$. Moreover, the subspace in $L$ of vectors
of $\gl(N)$-weight $\lambda$ if one-dimensional and spanned by $\xi$.

The following result is contained in [D2].
\bigskip

\proclaim
{\bf Theorem 2.3}
Every finite-dimensional irreducible representation $L$
of $\Y(N)$ is highest weight. Moreover, $L$ contains a unique
(up to scalar multiples) highest weight vector.
\endproclaim

\Proof The $\gl(N)$-module $L$ admits a weight space decomposition
$$
L=\underset{\lambda\in\frak{h}^*}\to{\bigoplus}L_{\lambda},\qquad
L_{\lambda}=\{\eta\in L\ |\ E_{ii}\eta=\lambda_i\eta,\quad
i=1,\dots,N\}.
\tag 2.4
$$
To see this, we note 
that by (2.3) the direct sum $L'$ of all weight subspaces
in $L$ is a $\Y(N)$-submodule. This submodule is nontrivial
because $L$ is finite-dimensional and so it contains an eigenvector
for the commutative subalgebra $\frak h$. Since $L$ is irreducible
as a $\Y(N)$-module we conclude that $L'=L$.
 
Introduce the subspace in $L$ 
$$
L^0:=\{\eta\in L\ |\ t_{ij}(u)\eta=0\qquad\text{for}\quad 
1\leq i<j\leq N\}. \tag 2.5
$$
By (2.3) the weight space decomposition (2.4)
induces that of $L^0$:
$$
L^0=\underset{\lambda\in\frak{h}^*}\to{\bigoplus}L^0_{\lambda},
\qquad L^0_{\lambda}=L^0\cap L_{\lambda}.
\tag 2.6
$$
Since $\gl(N)$-module $L$ is finite-dimensional there
exists a nonzero vector $\eta_0\in L$ of $\gl(N)$-weight
$\lambda'$ such that $\lambda'+\alpha$ is not a weight of $L$
for any $\alpha\in R_+$.
We see from (2.3) that the vectors $t_{ij}^{(r)}\eta_0$ with $i<j$
and $r\geq 1$ have the $\gl(N)$-weight $\lambda'+\varepsilon_i-
\varepsilon_j$ and therefore are equal to
zero. This proves that the subspace $L^0$ is nontrivial.

Relations (1.1) and (2.3) imply that $\Y(N)$ can be equipped with
a $\ZZ^N$-grading
$$
\Y(N)=\underset{\nu\in\ZZ^N}\to{\bigoplus} \Y(N)_{\nu},
$$
where
$$
\Y(N)_{\nu}=\{y\in\Y(N)\ |\ [E_{ii},y]=\nu_i\ts y,\quad i=1,\dots,N\}.
$$
In particular, $\Y(N)_0$ is a subalgebra in $\Y(N)$.

Denote by $I$ and $I'$ the left and right
ideals in $\Y(N)$ respectively generated by the coefficients
of the series $t_{ij}(u)$ and $t_{ji}(u)$,
where $1\leq i<j\leq N$. Set $I_0=\Y(N)_0\cap I$.
Exactly as in
the case of semisimple Lie algebras (see, e.g., [Di, Section 7.4])
one derives from the Poincar\'e--Birkhoff--Witt theorem for $\Y(N)$ 
that $I_0=I'\cap \Y(N)_0$ and that $I_0$ is a two-sided ideal in
$\Y(N)_0$. Moreover, $\Y(N)_0=\Y\oplus \ts I_0$, where $\Y$ is the subspace
in $\Y(N)_0$ with the basis formed by ordered monomials in the
generators $t_{ii}^{(r)}$. 

Consider the algebra $\Cal Y=\Y(N)_0/I_0$.
We see from (1.1) that $[t_{ii}(u),t_{ii}(v)]=0$ and
$$
[t_{ii}(u),t_{jj}(v)]={1\over
{u-v}}(t_{ji}(u)t_{ij}(v)-t_{ji}(v)t_{ij}(u)).
$$
So, if $i<j$ each commutator $[t_{ii}^{(r)},t_{jj}^{(s)}]$
belongs to $I_0$. This implies that the algebra $\Cal Y$
is commutative and freely generated by the elements 
$t_{ii}^{(r)}\mod I_0$.

Further, we note that the subspace $L^0$ has a natural structure of
$\Y(N)_0$-module and hence also of a
$\Cal Y$-module because $L^0$ is annihilated by $I_0$.
Therefore, $L^0$ contains a nonzero vector $\xi$ which is
a common eigenvector for all operators $t_{ii}^{(r)}$.
So, $\xi$ generates a highest weight $\Y(N)$-submodule in $L$.
Since $L$ is irreducible this submodule coincides with 
$L$ which proves that $L$ is highest weight.

Finally, if $\xi$ has $\gl(N)$-weight $\lambda_0$ then nonzero
weight subspaces in (2.6) can only correspond to the 
$\gl(N)$-weights $\lambda$ of the form $\lambda_0-\omega$, 
where $\omega$ is a $\ZZ_+$-linear combination of the $\alpha_i$.
Moreover, the subspace 
$L^0_{\lambda_0}$ is one-dimensional. However, each 
subspace $L^0_{\lambda}$
is $\Cal Y$-invariant. If it is nonzero for a certain 
$\lambda\ne\lambda_0$ then
a common eigenvector 
in $L^0_{\lambda}$
for all $t_{ii}^{(r)}$ generates a proper submodule is $L$,
which makes a contradiction. $\square$

\bigskip
\noindent
{\it Remark 2.4.} It follows from the proof that if $L$ is an irreducible
highest weight representation of $\Y(N)$ (not necessary finite-dimensional)
then the subspace (2.5) is one-dimensional and spanned by the
highest weight vector $\xi$ of $L$. $\square$
\medskip

Due to Theorem 2.3, to describe the finite-dimensional irreducible
representations of $\Y(N)$ it suffices to find necessary and sufficient
conditions for an irreducible highest weight
representation $L(\lambda(u))$
to be finite-dimensional.
The case of 
the simplest Yangian $\Y(2)$ plays
a key role in the proof. We shall be assuming now that $N=2$. 

The following result is due to V.~O.~Tarasov [T1], [T2].

\bigskip
\proclaim
{\bf Proposition 2.5} If $\dim L(\lambda(u))<\infty$ then
there exists a formal series 
$\varphi(u)\in 1+u^{-1}\C[[u^{-1}]]$ such that
$\varphi(u)\lambda_1(u)$ and 
$\varphi(u)\lambda_2(u)$ are polynomials in $u^{-1}$.
\endproclaim

\Proof Let $\xi$ denote the highest weight vector of $L(\lambda(u))$.
Since $L(\lambda(u))$ is finite-dimensional there exists
a nonnegative integer $k$ such that the vector $t_{21}^{(k+1)}\xi$
is a linear combination of
$t_{21}^{(1)}\xi,\dots,t_{21}^{(k)}\xi$. Take the minimum $k$ with this
property. Let us show that then for any vector $t_{21}^{(r)}\xi$ 
with $r\geq k+1$ we have
$$
t_{21}^{(r)}\xi=a_1^{(r)}\xi_1+\cdots+a_k^{(r)}\xi_k
\tag 2.7
$$
for some complex coefficients $a_i^{(r)}$, where $\xi_i:=t_{21}^{(i)}\xi$,
$i=1,\dots,k$.
Indeed, this is true for $r=k+1$
by our choice of $k$.
Relations (1.1) imply that
$$
t_{11}^{(2)}t_{21}(u)=t_{21}(u)(t_{11}^{(2)}-t_{11}^{(1)}-u)+
t_{21}^{(1)}t_{11}(u), \tag 2.8
$$
and taking the coefficient at $u^{-p}$ we 
obtain
$$
t_{11}^{(2)}t_{21}^{(p)}=-t_{21}^{(p+1)}+t_{21}^{(1)}t_{11}^{(p)}
+t_{21}^{(p)}(t_{11}^{(2)}-t_{11}^{(1)}).
\tag 2.9
$$
Write $\lambda_1(u)=1+\lambda_1^{(1)}u^{-1}+\lambda_1^{(2)}u^{-2}+\cdots$.
Then by (2.9), for any $p\geq 1$,
$$
t_{21}^{(p+1)}\xi=-t_{11}^{(2)}t_{21}^{(p)}\xi+
(\lambda_1^{(2)}-\lambda_1^{(1)})t_{21}^{(p)}\xi+
\lambda_1^{(p)}t_{21}^{(1)}\xi.
\tag 2.10
$$
Hence, for $p=1,\dots,k-1$ we have
$$
t_{11}^{(2)}\xi_p=\lambda_1^{(p)}\xi_1+
(\lambda_1^{(2)}-\lambda_1^{(1)})\xi_p-\xi_{p+1}
\tag 2.11
$$
and applying (2.7) with $r=k+1$, we get
$$
t_{11}^{(2)}\xi_k=\lambda_1^{(k)}\xi_1+
(\lambda_1^{(2)}-\lambda_1^{(1)})\xi_k-
(a_1^{(k+1)}\xi_1+\cdots+a_k^{(k+1)}\xi_k).
\tag 2.12
$$
Now (2.7) follows from (2.10)--(2.12) by an obvious induction on $r$.

So, rewriting (2.7) in terms of generating series we obtain
$$
t_{21}(u)\xi=a_1(u)\xi_1+\cdots+a_k(u)\xi_k,
\tag 2.13
$$
where $a_i(u)=u^{-i}+a_i^{(k+1)}u^{-k-1}+\cdots$.
Applying both sides of (2.8) to $\xi$ and using (2.13) we obtain
$$
(t_{11}^{(2)}-\lambda_1^{(2)}+\lambda_1^{(1)}+u)
(a_1(u)\xi_1+\cdots+a_k(u)\xi_k)=\lambda_1(u)\xi_1. \tag 2.14
$$
Using (2.11) and (2.12) and taking the coefficient at $\xi_p$,
$p=2,\dots,k$ in this relation we get
$$
-a_{p-1}(u)+ua_p(u)-a_p^{(k+1)}a_k(u)=0.
$$
This implies that for any $p=1,\dots,k$ one has $a_p(u)=A_p(u)a_k(u)$
for a monic polynomial $A_p(u)$ of degree $k-p$ in $u$. Further,
taking the coefficient at $\xi_1$ in (2.14) we obtain 
$\lambda_1(u)=B(u)a_k(u)$ for 
a monic polynomial $B(u)$ of degree $k$ in $u$.

Similarly, repeating this argument with the use of the relation
$$
t_{22}^{(2)}t_{21}(u)=t_{21}(u)(t_{22}^{(2)}+t_{22}^{(1)}+u)-
t_{21}^{(1)}t_{22}(u)
$$
instead of (2.8) we obtain that $\lambda_2(u)=C(u)a_k(u)$ for 
a monic polynomial $C(u)$ of degree $k$ in $u$. Thus, the desired
series is $\varphi(u)=(a_k(u)u^k)^{-1}$. $\square$
\medskip

Given $\lambda=(\lambda_1,\dots,\lambda_N)\in \frak h^*\simeq\C^N$ 
we shall denote by
$L(\lambda)$ the irreducible 
$\gl(N)$-module with the highest weight $\lambda$ with respect to
the basis $\{E_{11},\dots, E_{NN}\}$ in
$\frak h$.
The representation $L(\lambda)$ is finite-dimensional if and only if
$\lambda_i-\lambda_{i+1}\in\ZZ_+$ for all $i=1,\dots,N-1$.
The homomorphism (1.5) allows us to consider any
representation $L(\lambda)$
as a $\Y(N)$-module. It is obviously highest weight with 
$\lambda_i(u)=1+\lambda_iu^{-1}$, $i=1,\dots, N$.
The Hopf algebra structure on $\Y(N)$ (see (1.14)) allows one to 
naturally equip any tensor product $L_1\ot L_2$ of $\Y(N)$-modules
with a $\Y(N)$-action by the rule
$$
y\cdot (\xi_1\ot\xi_2):=\Delta(y)(\xi_1\ot\xi_2), 
\qquad y\in\Y(N),\quad \xi_i\in L_i. \tag 2.15
$$
\medskip

Let us now return to the case $N=2$. 
A basis in the irreducible $\gl(2)$-module $L(\alpha,\beta)$
with the highest weight $(\alpha,\beta)$ is formed by the elements
$(E_{21})^r\xi$, where $\xi$ is the highest weight vector and
$r=0,\dots,\alpha-\beta$ if $\alpha-\beta\in\ZZ_+$; or
$r\in\ZZ_+$ otherwise.
Using (2.15) we can construct representations of $\Y(2)$ of the form
$$
L=L(\alpha_1,\beta_1)\ot\cdots\ot L(\alpha_k,\beta_k).
\tag 2.16
$$
In accordance with (1.14) the generators of $\Y(2)$ act by the rule
$$
t_{ij}(u)(\eta_1\ot\cdots\ot \eta_k)
=\sum_{a_1,\dots,a_{k-1}}t_{ia_1}(u)\eta_1\ot
t_{a_1a_2}(u)\eta_2\ot\cdots\ot t_{a_{k-1}j}(u)\eta_k, \tag 2.17
$$
where $\eta_i\in L(\alpha_i,\beta_i)$. By (1.5) all generators 
$t_{ab}^{(r)}$ with $r\geq 2$ act as zero operators in each
$L(\alpha_i,\beta_i)$ and so, by (2.17) we have for any $\eta\in L$
$$
t_{ij}^{(r)}\eta=0,\qquad \text{for}\quad r\geq k+1.\tag 2.18
$$
For a set of complex numbers
$\{\alpha_1,\dots,\alpha_k,\beta_1,\dots,\beta_k\}$ we can find
re-enumerations of the
$\alpha_i$ and of the $\beta_i$ such that the following
condition is satisfied: 
$$
\aligned
&\text{for every}\ \  i=1,\dots, k\ \  \text{we have: } \\
&\text{if the set}\ \ \{\alpha_p-\beta_q\ |\ i\leq p,q \leq k\}\cap \ZZ_+ 
\ \ \text{is not empty}\\
&\text{then}\ \  \alpha_i-\beta_i\ \  
\text{is its minimal element.}
\endaligned
\tag 2.19
$$
Indeed, choose first a minimal nonnegative integer difference among
all differences $\alpha_p-\beta_q$ if it exists. Re-enumerating the 
indices if necessary, we may assume that this difference is 
$\alpha_1-\beta_1$. Then we proceed by induction, considering
the differences $\alpha_p-\beta_q$
with $p,q\geq 2$ and so on.

The following proposition was proved in [T2].

\bigskip
\proclaim
{\bf Proposition 2.6} Let $L(\lambda_1(u),\lambda_2(u))$ 
be an irreducible highest weight representation
of $\Y(2)$ and $\lambda_1(u)$ and $\lambda_2(u)$ be polynomials
of degree ${}\leq k$ in $u^{-1}$. Fix decompositions
$$
\align
\lambda_1(u)&=(1+\alpha_1u^{-1})\cdots (1+\alpha_ku^{-1}),
\tag 2.20\\
\lambda_2(u)&=(1+\beta_1u^{-1})\cdots (1+\beta_ku^{-1})
\tag 2.21
\endalign
$$
and assume that condition (2.19) is satisfied.
Then $L(\lambda_1(u),\lambda_2(u))$ is isomorphic to the tensor product 
(2.16).
\endproclaim

\Proof Let $\xi_i$ be the highest weight vector of 
$L(\alpha_i,\beta_i)$. We have 
$$
t_{11}(u)\xi_i=(1+\alpha_i u^{-1})\xi_i,\quad
t_{22}(u)\xi_i=(1+\beta_i u^{-1})\xi_i,
\quad\text{and}\quad t_{12}(u)\xi_i=0. \tag 2.22
$$
Therefore, the vector $\xi:=\xi_1\ot\cdots\ot \xi_k$ satisfies
the relations
$$
t_{11}(u)\xi=\lambda_1(u)\xi,\quad
t_{22}(u)\xi=\lambda_2(u)\xi\quad\text{and}\quad
t_{12}(u)\xi=0,\tag 2.23
$$
which follows immediately from (2.17). 
Thus, the cyclic  $\Y(2)$-span of the vector $\xi$ in $L$ is
a highest weight module with the highest weight 
$(\lambda_1(u),\lambda_2(u))$.
Therefore, to prove Proposition 2.6 we only need to show that
the representation $L$ is irreducible. We use induction on $k$.

For $k=1$ this is obvious, because
$L(\alpha,\beta)$ is an irreducible
$\gl(2)$-module and hence also irreducible as a $\Y(2)$-module.

Suppose now that $k> 1$. Let us prove first that 
any nonzero submodule $\wt L$ of $L$ contains the vector $\xi$.
Note that $\wt L$ contains a nonzero vector $\eta$ such that 
$t_{12}(u)\eta=0$. Indeed, $\gl(2)$-weights of $L$ have the form
$(\sum_i\alpha_i-r, \sum_i\beta_i+r)$ with $r\in\ZZ_+$. Therefore,
there exists a $\gl(2)$-weight $(\mu_1,\mu_2)$ of $\wt L$ such that
$(\mu_1+1,\mu_2-1)$ is not a weight. Then a vector $\eta$ of 
$\gl(2)$-weight
$(\mu_1,\mu_2)$ is annihilated by $t_{12}(u)$.

Write $\eta$ in the form
$$
\eta=\sum_{r=0}^p(E_{21})^r \xi_1\ot \eta_r,
$$
with $p\geq 0$, 
where $\eta_r\in L(\alpha_2,\beta_2)\ot\cdots\ot L(\alpha_k,\beta_k)$; and
if $\alpha_1-\beta_1\in\ZZ_+$ then $p\leq  \alpha_1-\beta_1$. We may
assume that $\eta_p\ne 0$.

Applying $t_{12}(u)$ to $\eta$ with the use of (1.14) and (2.15) we get
$$
\sum_{r=0}^p\left(t_{11}(u)(E_{21})^r \xi_1\ot t_{12}(u) \eta_r
+t_{12}(u)(E_{21})^r \xi_1\ot t_{22}(u) \eta_r\right)=0.\tag 2.24
$$
Using (1.5) and the commutation relations in $\gl(2)$ we can write
$$
t_{11}(u)(E_{21})^r \xi_1=(1+(\alpha_1-r)u^{-1})(E_{21})^r \xi_1
\tag 2.25
$$
and
$$
t_{12}(u)(E_{21})^r \xi_1=u^{-1}\ts r(\alpha_1-\beta_1-r+1)
(E_{21})^{r-1} \xi_1.
\tag 2.26
$$
Taking the coefficient at $(E_{21})^p \xi_1$ in (2.24) we get
$$
(1+(\alpha_1-p)u^{-1})t_{12}(u)\eta_p=0.
$$
This implies that $t_{12}(u)\eta_p=0$. On the other hand, by the induction
hypotheses, the representation 
$L(\alpha_2,\beta_2)\ot\cdots\ot L(\alpha_k,\beta_k)$ is irreducible, and
so
$$
\eta_p=\ \text{const}\cdot \xi_2\ot\cdots\ot\xi_k,\tag 2.27
$$
where `const' is a nonzero constant; see Remark 2.4.

Suppose now that $p\geq 1$ and take the coefficient at 
$(E_{21})^{p-1} \xi_1$ in (2.24). By (2.25) and (2.26) we have
$$
(1+(\alpha_1-p+1)u^{-1})t_{12}(u)\eta_{p-1}+
u^{-1}\ts p(\alpha_1-\beta_1-p+1) t_{22}(u)\eta_p=0.
\tag 2.28
$$
Relations (2.27) and (2.17) imply that
$$
t_{22}(u)\eta_p=(1+\beta_2u^{-1})\cdots (1+\beta_ku^{-1}) \eta_p.
\tag 2.29
$$
Multiply both sides of (2.28) by $u^k$:
$$
(u+\alpha_1-p+1)u^{k-1}t_{12}(u)\eta_{p-1}+
\ts p(\alpha_1-\beta_1-p+1)u^{k-1} t_{22}(u)\eta_p=0.
\tag 2.30
$$
The operators $u^{k-1}t_{ij}(u)$ in 
$L(\alpha_2,\beta_2)\ot\cdots\ot L(\alpha_k,\beta_k)$
are polynomials in $u$ by (2.18) and we may put $u=-\alpha_1+p-1$
in (2.30). Taking into account (2.29) we get
$$
p(\alpha_1-\beta_1-p+1)(\alpha_1-\beta_2-p+1)
\cdots (\alpha_1-\beta_k-p+1)=0.
$$
However, this is impossible because of condition (2.19).
Therefore, $p$ has to be equal to $0$ and so, $\eta$ is
a multiple of $\xi$.

To complete the proof we have to show that the submodule of $L$
generated by $\xi$ coincides with $L$.

Using the antiautomorphism (1.7) we equip the vector space $V^*$
dual to a $\Y(2)$-module $V$ with a structure of
$\Y(2)$-module. Namely, set
$$
(y\cdot f)(v)=f(\sigma(y)v),\qquad y\in\Y(2),\quad f\in V^*,\quad v\in V.
$$
Then one can easily see that 
$L(\alpha_i,\beta_i)^*\simeq L'(-\alpha_i,-\beta_i)$,
where $L'(-\alpha_i,-\beta_i)$ is the lowest weight representation
of $\gl(2)$ with the lowest weight vector $\xi'_i$ dual to $\xi_i$,
of the weight
$(-\alpha_i,-\beta_i)$.
This implies that
$$
L^*\simeq  L'(-\alpha_1,-\beta_1)\ot\cdots\ot L'(-\alpha_k,-\beta_k).
\tag 2.31
$$
Similarly to the case of representation $L$, any nonzero submodule of $L^*$
contains a nonzero vector $\eta'$ such that $t_{21}(u)\eta'=0$.
Modifying the previous argument for lowest weight representations
we prove that $\eta'$ coincides (up to a nonzero factor) with
the vector $\xi':=\xi'_1\ot\cdots\ot \xi'_k$.
Now, if the submodule in $L$ generated by $\xi$ is proper, then
its annihilator in $L^*$ is a nonzero submodule which does not
contain $\xi'$. Contradiction. $\square$
\medskip

Proposition 2.6 allows one to get necessary and sufficient conditions
for $L(\lambda(u))$ to be finite-dimensional.

\bigskip
\proclaim
{\bf Proposition 2.7} Let $L(\lambda(u))$ be an irreducible representation
of $\Y(2)$ with the highest weight 
$\lambda(u)=(\lambda_1(u),\lambda_2(u))$. Then $L(\lambda(u))$ is
finite-dimensional if and only if there exists a formal series
$\varphi(u)\in 1+u^{-1}\C[[u^{-1}]]$ such that 
$\varphi(u)\lambda_1(u)$ and $\varphi(u)\lambda_2(u)$ are polynomials
in $u^{-1}$ with decompositions
$$
\align
\varphi(u)\lambda_1(u)&=(1+\alpha_1u^{-1})\cdots (1+\alpha_ku^{-1}),
\tag 2.32\\
\varphi(u)\lambda_2(u)&=(1+\beta_1u^{-1})\cdots (1+\beta_ku^{-1})
\tag 2.33
\endalign
$$
and $\alpha_i-\beta_i\in\ZZ_+$, $i=1,\dots,k$ for certain re-enumerations
of the $\alpha_i$ and $\beta_i$.
\endproclaim

\Proof Let $\dim L(\lambda(u))<\infty$.
By Proposition 2.6,  
taking the composition of our representation with
an automorphism of the form (1.11) we may assume that (2.32) and (2.33)
are satisfied for a certain series $\varphi(u)$. Re-enumerate the
$\alpha_i$
and $\beta_i$ in such a way that condition (2.19)
is satisfied. Now use Proposition 2.6. Since representation
(2.16) is finite-dimensional this gives the relations
$\alpha_i-\beta_i\in\ZZ_+$ for all $i=1,\dots,k$.

Conversely, let (2.32) and (2.33) hold with 
$\alpha_i-\beta_i\in\ZZ_+$, $i=1,\dots,k$. Then representation
(2.16) is finite-dimensional and the composition of $L(\lambda(u))$
with the automorphism (1.11) corresponding to the series $\varphi(u)$
is isomorphic to the irreducible
quotient of the $\Y(2)$-cyclic span in $L$ of the vector $\xi$;
see (2.23).
Therefore, $\dim L(\lambda(u))<\infty$. $\square$
\medskip

The following reformulation of this result is due to V.~G.~Drinfeld [D2].

\bigskip
\proclaim
{\bf Theorem 2.8} The irreducible highest weight representation
$L(\lambda_1(u),\lambda_2(u))$ 
of $\Y(2)$ is
finite-dimensional if and only if there exists a monic
polynomial $P(u)\in\C[u]$ such that
$$
\frac{\lambda_1(u)}{\lambda_2(u)}=\frac{P(u+1)}{P(u)}.\tag 2.34
$$
In this case $P(u)$ is unique.
\endproclaim

\Proof If $\dim L(\lambda_1(u),\lambda_2(u))<\infty$ then by 
Proposition 2.7 relations (2.32) and (2.33) hold for a certain
series $\varphi(u)$.
If $\alpha_i-\beta_i\in\ZZ_+$ for all $i$ we put
$$
P(u)=\prod_{i=1}^k (u+\beta_i)\cdots (u+\alpha_i-1).
$$
Conversely, suppose (2.34) is satisfied for $P(u)=(u+\gamma_1)\cdots
(u+\gamma_s)$. Define the polynomials
$$
\lambda'_1(u) =(1+(\gamma_1+1)u^{-1})\cdots (1+(\gamma_s+1)u^{-1})
$$
and
$$
\lambda'_2(u) =(1+\gamma_1 u^{-1})\cdots (1+\gamma_s u^{-1}).
$$
By Proposition 2.7, the irreducible highest weight representation
$L(\lambda'_1(u),\lambda'_2(u))$ is finite-dimensional.
However, $\lambda'_1(u)$ and $\lambda'_2(u)$ satisfy (2.34) and so,
there exists an automorphism of $\Y(2)$ of the form (1.11)
such that its composition with the representation
$L(\lambda'_1(u),\lambda'_2(u))$ is isomorphic to
$L(\lambda_1(u),\lambda_2(u))$. Thus, the latter
is also finite-dimensional.

To complete the proof of the theorem, suppose that $Q(u)$ is another
monic polynomial in $u$ and
$$
\frac{P(u+1)}{P(u)}=\frac{Q(u+1)}{Q(u)}.
$$
Then $P(u)/Q(u)$  is periodic in 
$u$ which is only possible when $P(u)=Q(u)$.
$\square$
\medskip

The polynomial $P(u)$ defined by (2.34)
is called the {\it Drinfeld polynomial\/}; see also Corollary 2.13 below.

Relation (2.34) motivates the following notation which will be
used in the sequel. For two formal series $\lambda_1(u)$ and 
$\lambda_1(u)$
in $u^{-1}$ we shall write
$$
\lambda_1(u)\to\lambda_2(u), \tag 2.35
$$
if there exists a monic polynomial $P(u)$ in $u$ such that (2.34)
holds.

Due to the decomposition (1.12) we can establish a correspondence
between representations of $\Y(N)$ and $\Y(\sll(N))$
analogous to the correspondence between 
the representations of the Lie algebras $\gl(N)$ and $\sll(N)$.
Consider the following {\it similarity classes\/} of 
finite-dimensional irreducible representations of $\Y(N)$:
two representations belong to the same class if
one can be obtained from the other by the composition
with an automorphism (1.11). Since the subalgebra $\Y(\sll(N))$
is stable with respect to all automorphisms of the form (1.11),
any representation
from such a class, restricted to the subalgebra $\Y(\sll(N))$
gives the same representation of $\Y(\sll(N))$. 
Moreover, all finite-dimensional irreducible representations
of $\Y(\sll(N))$ can be obtained in this way.

In particular, using Theorem 2.8
we obtain a complete description of finite-dimensional irreducible
representations of the Yangian $\Y(\sll(2))$.

\bigskip
\proclaim
{\bf Corollary 2.9} There is a one-to-one correspondence
between
finite-dimensional irreducible representations
of the Yangian $\Y(\sll(2))$ and monic polynomials in $u$. Every such
representation is isomorphic to a representation of the form (2.16).
\endproclaim

\Proof The first claim follows from the fact that the ratios 
$\lambda_1(u)/\lambda_2(u)$ take the same value for all
representations $L(\lambda_1(u),\lambda_2(u))$ of $\Y(2)$ of a given
similarity class and hence by Theorem 2.8 the similarity
classes are parameterized by monic polynomials in $u$.

Further, given a monic polynomial $P(u)$ we may find
polynomials $\lambda_1(u)$ and $\lambda_2(u)$ in $u^{-1}$ such that
relation (2.34) holds; see the proof of Theorem 2.8. 
Then by Proposition 2.6, 
the representation $L(\lambda_1(u),\lambda_2(u))$ of $\Y(2)$
is isomorphic to a tensor product (2.16) and its restriction to
$\Y(\sll(2))$ gives the finite-dimensional irreducible representation
corresponding to $P(u)$. $\square$

\bigskip
\noindent
{\it Remark 2.10.} Tensoring $L(\lambda_1(u),\lambda_2(u))$ 
by factors of the form 
$L(\gamma,\gamma)$ is equivalent to multiplying both
$\lambda_1(u)$, and $\lambda_2(u)$ simultaneously by
$\gamma(u)=1+\gamma u^{-1}$.
However, the representations 
$L(\lambda_1(u),\lambda_2(u))$ and 
$L(\gamma(u)\lambda_1(u),\gamma(u)\lambda_2(u))$
of $\Y(2)$ belong to the same similarity class.
Therefore, we might add to Corollary 2.9 the condition that in (2.16)
all differences $\alpha_i-\beta_{i}$ are positive integers. $\square$
\medskip

Since all finite-dimensional irreducible representation
of $\Y(\sll(2))$ have the form of (2.16) it is useful
to have a criterion of irreducibility of such representations.

Following [CP1] we define the {\it string\/} corresponding
to a pair of complex numbers $(\alpha,\beta)$ with $\alpha-\beta\in\ZZ_+$
as the set
$$
S(\alpha,\beta)=\{\beta,\beta+1,\dots,\alpha-1\}\subset\C.
$$
If $\alpha=\beta$ the set $S(\alpha,\beta)$ is considered to be 
empty.
We say that two strings $S_1$ and $S_2$ are 
{\it in general position\/}
if either

(i) $S_1\cup S_2$ is not a string, or

(ii) $S_1\subseteq S_2$, or $S_2\subseteq S_1$.

The following result is contained in [CP1] and [CP3].

\bigskip
\proclaim
{\bf Corollary 2.11} Consider the representation $L$ 
of $\ts\Y(2)$ or $\Y(\sll(2))$ given by (2.16) where
all the differences $\alpha_i-\beta_i$ are nonnegative integers.
Then $L$ is irreducible if and only if the strings
$S(\alpha_1,\beta_1),\dots,S(\alpha_k,\beta_k)$ are pairwise in
general position.
\endproclaim

\Proof Suppose that the strings are pairwise in general position and
assume first that
$\alpha_1-\beta_1\leq\cdots\leq \alpha_k-\beta_k$.
Then one easily checks that (2.19) holds and so,
by Proposition 2.6 representation (2.16) is irreducible.
To complete the proof of the `if' part note that
if representation (2.16)
is irreducible, then any permutation of the tensor factors gives
an isomorphic representation $L'$. Indeed, the submodule
in $L'$ generated by the tensor product of the highest weight
vectors of the $L(\alpha_i,\beta_i)$ 
is highest weight with the same $\lambda_1(u)$ and
$\lambda_2(u)$ given by (2.20) and (2.21). Therefore, $L$ is
isomorphic to a subquotient of $L'$. As $L$ and $L'$ have the same
dimension, they have to be isomorphic.

Conversely, let $k=2$ and $L(\alpha_1,\beta_1)\ot L(\alpha_2,\beta_2)$
be irreducible. Suppose that the strings 
$S(\alpha_1,\beta_1)$ and $S(\alpha_2,\beta_2)$
are not in general position. Then $\alpha_2-\beta_1$, $\alpha_1-\beta_2
\in \ZZ_+$ and
either both of
$\beta_1-\beta_2$ and $\alpha_1-\alpha_2$,
or both of $\beta_2-\beta_1$ and $\alpha_2-\alpha_1$
are positive integers. The strings
$S(\alpha_1,\beta_2)$ and $S(\alpha_2,\beta_1)$ are in general position
and so, the representation $L(\alpha_1,\beta_2)\ot L(\alpha_2,\beta_1)$
is irreducible. Since it has the same highest weight
as the irreducible representation 
$L(\alpha_1,\beta_1)\ot L(\alpha_2,\beta_2)$ these two
representations have to be
isomorphic and, in particular, have the same dimension:
$$
(\alpha_1-\beta_1+1)(\alpha_2-\beta_2+1)=
(\alpha_1-\beta_2+1)(\alpha_2-\beta_1+1).
$$
This implies that $(\alpha_1-\alpha_2)(\beta_1-\beta_2)=0$,
which is impossible.

Now consider the general case. Suppose that $L$ is irreducible but
two strings $S(\alpha_i,\beta_i)$ and $S(\alpha_j,\beta_j)$
are not in general position. Permuting the factors in (2.16) 
if necessary, we may assume that $i$ and $j$ are adjacent. However,
the representation $L(\alpha_i,\beta_i)\ot L(\alpha_j,\beta_j)$
is reducible. Contradiction.
$\square$
\medskip

We are now able to prove Drinfeld's classification theorem for the
Yangian $\Y(N)$, $N$ is an arbitrary positive integer; see [D2].
Theorem 2.3 together with the following theorem give a complete description
of finite-dimensional irreducible representations of $\Y(N)$.
We use notation (2.35). 

\bigskip
\proclaim
{\bf Theorem 2.12} The irreducible highest weight representation
$L(\lambda(u))$, $\lambda(u)=(\lambda_1(u),\dots,\lambda_N(u))$ of 
$\Y(N)$ is finite-dimensional if and only if 
the following relation holds:
$$
\lambda_1(u)\to  \lambda_2(u)\to\cdots\to\lambda_N(u). \tag 2.36
$$
\endproclaim

\Proof Suppose that $\dim L(\lambda(u))<\infty$.
By the Poincar\'e--Birkhoff--Witt theorem for
$\Y(N)$, 
for any $i\in\{1,\dots,N-1\}$ the subalgebra in $\Y(N)$ 
generated by the coefficients of the series 
$t_{ii}(u)$, $t_{i,i+1}(u)$, $t_{i+1,i}(u)$, $t_{i+1,i+1}(u)$ 
is isomorphic to
$\Y(2)$. The cyclic span of the highest weight vector of $L(\lambda(u))$
with respect to this subalgebra is a representation with
the highest weight $(\lambda_i(u),\lambda_{i+1}(u))$. Its irreducible
quotient
is finite-dimensional and so, applying Theorem 2.8 we get (2.36).

Conversely, let
$$
\frac{\lambda_i(u)}{\lambda_{i+1}(u)}=\frac{P_i(u+1)}{P_i(u)}\tag 2.37
$$
for certain 
polynomials $P_i(u)$, $i=1,\dots,N-1$ given by
$$
P_i(u)=(u+\gamma_1^{(i)})\cdots (u+\gamma_{s_i}^{(i)}).
$$
Consider first the irreducible highest weight representation 
$L(\lambda'(u))$ of $\Y(N)$
with the highest weight $\lambda'(u)=(\lambda'_1(u),\dots,\lambda'_N(u))$,
where
$$
\align
\lambda'_i(u)={}&\prod_{a=1}^{i-1}(1+\gamma_1^{(a)}u^{-1})
\cdots(1+\gamma_{s_a}^{(a)}u^{-1})\\
\times{}&\prod_{a=i}^{N-1}(1+(\gamma_1^{(a)}+1)u^{-1})
\cdots(1+(\gamma_{s_a}^{(a)}+1)u^{-1}).
\endalign
$$
These polynomials have the form
$$
\lambda'_i(u)=(1+\alpha_1^{(i)}u^{-1})
\cdots(1+\alpha_{k}^{(i)}u^{-1})
$$
with $k=\sum s_i$
and the condition $\alpha_r^{(i)}-\alpha_r^{(i+1)}\in\ZZ_+$
is satisfied for each $r=1,\dots,k$ and $i=1,\dots,N-1$.
Consider the tensor product
$$
L(\alpha_1)\ot\cdots\ot L(\alpha_k),\tag 2.38
$$
where $L(\alpha_r)$ is the irreducible representation of
$\gl(N)$ with the highest weight 
$\alpha_r=(\alpha_r^{(1)},\dots,\alpha_r^{(N)})$.
The generators $t_{ij}(u)$ act in representation (2.38) by the rule
(2.17). So, the tensor product of the highest weight vectors of
the $L(\alpha_r)$ generates a representation of $\Y(N)$ with the highest
weight $\lambda'(u)$. Hence, its irreducible quotient is isomorphic to
$L(\lambda'(u))$.
Representation (2.38) is finite-dimensional, 
and therefore $L(\lambda'(u))$ is.

Finally, note that the $\lambda'_i(u)$ satisfy (2.37) and
so, there exists an automorphism of $\Y(N)$ of the form (1.11)
such that its composition with the representation
$L(\lambda'(u))$ is isomorphic to
$L(\lambda(u))$ which completes the proof. $\square$
\medskip

We conclude this section with a description of 
finite-dimensional irreducible representations of the Yangian 
$\Y(\sll(N))$ implied by Theorem 2.12; see [D2].

\bigskip
\proclaim
{\bf Corollary 2.13} There is a one-to-one correspondence
between
finite-dimensional irreducible representations
of the Yangian $\Y(\sll(N))$ 
and the families of
monic polynomials $\{P_1(u),\dots,P_{N-1}(u)\}$ in $u$. 
Every such representation is
isomorphic to a subquotient of a representation of the form (2.38).
$\square$
\endproclaim

The $P_i(u)$ are called the 
{\it Drinfeld polynomials\/} associated to the representation
$L(\lambda(u))$ of $\Y(\sll(N))$. 
Note that our definition of $P_i(u)$ differs from
that of [D2] by a shift in $u$.

\bigskip
\noindent
{\it Remark 2.14.} Contrary to the case $N=2$, it is not true for $N\geq 3$
that every finite-dimensional irreducible representation of 
$\Y(\sll(N))$ is isomorphic to a tensor product (2.38).
However, this is true for `generic' representations; see, e.g., [KR], [M2],
[NT2].

\bigskip
\bigskip
\noindent
{\bf 3. Highest weight representations of the twisted Yangians}
\bigskip

Here we define highest weight representations of the
algebras $\Y^{\pm}(N)$ and
prove an analogue of Theorem 2.3 for $\Y^{\pm}(N)$. We also find
some necessary conditions for an irreducible highest weight
representation to be finite-dimensional and construct
a family of representations of $\Y^{\pm}(N)$.

From now on till the
end of the paper we shall be assuming that the generators
$t_{ij}(u)$ of $\Y(N)$ as well as the basis elements $E_{ij}$
of the Lie algebra $\gl(N)$
are parameterized by the indices
$i,j\in\{-n,-n+1,\dots,n\}$, where $n=[N/2]$ and the index $0$ is skipped
in $N=2n$ (see Section 1).
\medskip

\noindent
{\bf Definition 3.1.} A representation $V$ of the twisted Yangian
$\Y^{\pm}(N)$
is called {\it highest weight\/} if there exists a nonzero vector 
$\xi\in V$ such that $V$ is generated by $\xi$,
$$
s_{ij}(u)\xi=0 \qquad \text{for} \quad -n\leq i<j\leq n, \tag 3.1
$$
and
$$
s_{ii}(u)\xi=\mu_i(u)\xi \tag 3.2
$$
for some formal series $\mu_i(u)\in 1+u^{-1}\C[[u^{-1}]]$, 
where $i=1,\dots,n$ if $N=2n$ and $i=0,\dots,n$ if $N=2n+1$.
In this case the vector $\xi$ is called the {\it highest weight vector\/}
of $V$ and the set $\mu(u):=(\mu_1(u),\dots,\mu_n(u))$ or
$\mu(u):=(\mu_0(u),\dots,\mu_n(u))$
is the {\it highest weight\/} of $V$. $\square$
\medskip

Due to the symmetry relation (1.18), the vector $\xi$ is also an
eigenvector for the generators
$s_{-i,-i}(u)$, $i=1,\dots,n$. Moreover,
if $N=2n+1$ the component $\mu_0(u)$ of the highest weight
has to be even: $\mu_0(u)=\mu_0(-u)$.

If we replace (3.1) by the relations
$$
s_{ji}(u)\xi=0 \qquad \text{for} \quad -n\leq i<j\leq n.
$$
we get the definition of a {\it lowest weight representation\/}.

\bigskip

\noindent
{\bf Definition 3.2.} Let $\mu(u)=(\mu_1(u),\dots,\mu_n(u))$
or $\mu(u)=(\mu_0(u),\dots,\mu_n(u))$
be any set of formal series where $\mu_i(u)\in 1+u^{-1}\C[[u^{-1}]]$
for $i=1,\dots,n$ and $\mu_0(u)\in 1+u^{-2}\C[[u^{-2}]]$.
The {\it Verma module\/} $M(\mu(u))$ is the quotient of $\Y^{\pm}(N)$
by the left ideal generated by all the coefficients of the series
$s_{ij}(u)$ for $-n\leq i<j\leq n$ and $s_{ii}(u)-\mu_i(u)$ for 
$i=1,\dots,n$ or $i=0,\dots,n$, respectively. $\square$
\medskip

The following properties of $M(\mu(u))$ are
immediate from
the Poincar\'e--Birkhoff--Witt theorem
for $\Y^{\pm}(N)$ (see Section 1).
First,
ordered monomials in the generators (1.19) or (1.20)
with the additional condition $i>j$
(assuming that (1.21) is satisfied)
form a basis of 
$M(\mu(u))$. Further, $M(\mu(u))$ is a representation of $\Y^{\pm}(N)$
with the highest weight $\mu(u)$ and any representation of
the highest weight $\mu(u)$ is isomorphic to a quotient
of $M(\mu(u))$. The module $M(\mu(u))$ has a unique irreducible
quotient which we denote by $V(\mu(u))$. 
Moreover, if an irreducible $\Y^{\pm}(N)$-module $V$ contains a nonzero
vector $\xi$ satisfying (3.1) and (3.2) then $V$ is isomorphic to
$V(\mu(u))$. 

Using the inclusion (1.22) we may regard any $\Y^{\pm}(N)$-module
as a $\g(N)$-module. Relations (1.17)
imply that
$$
[F_{ij},s_{kl}(u)]=\delta_{kj}s_{il}(u)-\delta_{il}s_{kj}(u)
-\theta_{k,-j}\delta_{i,-k}s_{-j,l}(u)+
\theta_{i,-l}\delta_{-l,j}s_{k,-i}(u),
\tag 3.3
$$
where we identify the elements $F_{ij}\in\g(N)$ with their
images in $\Y^{\pm}(N)$. 
Set $\frak{h}=\ \text{span of}\ \ \{F_{11},\dots,F_{nn}\}$.
Then $\frak{h}$ is a Cartan subalgebra in $\g(N)$.
Choose the standard system of positive roots $R_+$ for 
$\g(N)$ with respect to $\frak h$ corresponding to
the upper triangular Borel subalgebra and let
$\alpha_1,\dots,\alpha_{n}$ be a basis in $R_+$.
If $V$ is a representation of $\Y^{\pm}(N)$ with the highest weight
$\mu(u)$ then its highest weight vector $\xi$ is also a $\g(N)$-highest
weight vector of the $\g(N)$-weight
$\mu=(\mu_1,\dots,\mu_n)$ where $\mu_i$ is the coefficient
at $u^{-1}$ in the series $\mu_i(u)$. By (3.3) 
all $\g(N)$-weights of $V$ have the form 
$\mu-\omega$, where $\omega$ is a $\ZZ_+$-linear combination 
of the $\alpha_i$. Moreover, the subspace in $V$ of vectors
of $\g(N)$-weight $\mu$ if one-dimensional and spanned by $\xi$.

\bigskip

\proclaim
{\bf Theorem 3.3}
Every finite-dimensional irreducible representation $V$
of $\ts\Y^{\pm}(N)$ is highest weight. Moreover, $V$ contains a unique
(up to scalar multiples) highest weight vector.
\endproclaim

\Proof The proof follows the same steps as that of Theorem 2.3.
We set
$$
V^0:=\{\eta\in V\ |\ s_{ij}(u)\eta=0\qquad\text{for}\quad 
-n\leq i<j\leq n\}. \tag 3.4
$$
By (3.3) we have the weight space decomposition
$$
V^0=\underset{\mu\in\frak{h}^*}\to{\bigoplus}V^0_{\mu},
\tag 3.5
$$
where
$$
V^0_{\mu}=\{\eta\in V^0\ |\ F_{ii}\eta=\mu_i\eta\},
\qquad \mu_i=\mu(F_{ii}).
$$
Introduce the $\frak h$-centralizer in $\Y^{\pm}(N)$:
$$
\Y^{\pm}(N)_{0}=\{y\in\Y^{\pm}(N)\ |\ 
[F_{ii},y]=0\quad\text{for all}\ \ i\}.
$$
Denote by $I$ the left
ideal in $\Y^{\pm}(N)$ generated by the coefficients
of the series $s_{ij}(u)$,
where $-n\leq i<j\leq n$ and set $I_0=\Y^{\pm}(N)_0\cap I$.
We derive from the Poincar\'e--Birkhoff--Witt theorem for $\Y^{\pm}(N)$ 
that $I_0$ is a two-sided ideal in
$\Y^{\pm}(N)_0$. Next, using (1.17) we verify that the algebra
$\Cal Y:=\Y^{\pm}(N)_0/I_0$ is commutative
and generated by the elements $s_{ii}^{(r)}\mod I_0$.

Further, $V^0$ is a nonzero
$\Cal Y$-module (cf. the proof of Theorem 2.3)
and the weight spaces $V^0_{\mu}$ are
$\Cal Y$-invariant. If $V^0_{\mu_0}\ne 0$ for a certain $\mu_0$
then there exists a nonzero element $\xi\in V^0_{\mu_0}$ which is
a common eigenvector for the generators $s_{ii}^{(r)}$. So,
$\xi$ generates a highest weight $\Y^{\pm}(N)$-submodule in $V$.
Since $V$ is irreducible this submodule coincides with 
$V$ which proves that $V$ is highest weight.

Finally, the subspace 
$V^0_{\mu_0}$ is one-dimensional
and nonzero
weight subspaces in (3.5) can only correspond to the 
$\g(N)$-weights $\mu$ of the form $\mu_0-\omega$, 
where $\omega$ is a $\ZZ_+$-linear combination of the $\alpha_i$.
If $V^0_{\mu}$ is nonzero for a certain $\mu\ne\mu_0$ then
a common eigenvector in $V^0_{\mu}$
for all $s_{ii}^{(r)}$ generates a proper submodule is $V$,
which makes a contradiction. $\square$

\bigskip
\noindent
{\it Remark 3.4.} The above proof implies that if $V$ is an irreducible
highest weight representation of $\Y^{\pm}(N)$ 
then the subspace (3.4) is one-dimensional and spanned by the
highest weight vector $\xi$ of $V$. $\square$
\medskip

Due to Theorem 3.3, to describe the finite-dimensional irreducible
representations of $\Y^{\pm}(N)$ 
it suffices to find necessary and sufficient
conditions for an irreducible highest weight
representation $V(\mu(u))$
to be finite-dimensional. Some
necessary conditions are provided by the following proposition.
We use notation (2.35).

\bigskip
\proclaim
{\bf Proposition 3.5} If the irreducible highest weight
representation $V(\mu(u))$ of the algebra $\Y^{\pm}(N)$
is finite-dimensional then the following relation holds:
$$
\mu_1(u)\to  \mu_2(u)\to\cdots\to\mu_n(u). 
$$
\endproclaim

\Proof Let $J$ be the left ideal of $\Y^{\pm}(N)$ generated by
the coefficients of the series $s_{-i,j}(u)$, where $i,j=1,\dots,n$.
Denote by $\Norm J$ the normalizer of $J$:
$$
\Norm J=\{y\in \Y^{\pm}(N)\ |\ Jy\subseteq J\}.
$$
The annihilator $V(\mu(u))^{J}$ of $J$ in $V(\mu(u))$ 
has a natural structure
of a $\Norm J/J$-module. Note that the highest weight
vector of $V(\mu(u))$ is contained in $V(\mu(u))^{J}$.
Defining relations (1.17) show that the mapping
$$
t_{ij}^{(r)}\mapsto s_{ij}^{(r)} \mod\ J,\qquad 
1\leq i,j\leq n
$$
defines an algebra homomorphism
$$
\Y(n)\to  \Norm J/J.
$$
By making use of this homomorphism
we may regard $V(\mu(u))^{J}$ as a
$\Y(n)$-module. The cyclic span of the vector $\xi$ is a finite-dimensional
representation of $\Y(n)$ with the highest weight
$(\mu_1(u),\dots,\mu_n(u))$. Applying Theorem 2.12 we complete the proof.
$\square$
\medskip

Given $\mu=(\mu_1,\dots,\mu_n)\in \frak h^*$ 
consider the irreducible 
$\g(N)$-module with the highest weight $\mu$ with respect to
the basis $\{F_{11},\dots,F_{nn}\}$ in $\frak h^*$.
It is finite-dimensional if and only if the following
conditions are satisfied (see, e.g., [Di, Section 7.2]):
$$
\mu_{i+1}-\mu_i\in\ZZ_+\qquad\text{for}\quad i=1,\dots,n-1
\tag 3.6
$$
and
$$
\align
-\mu_1-\mu_2\in\ZZ_+ \qquad&\text{for}\quad \oa(2n), \tag 3.7\\
-\mu_1\in\ZZ_+ \qquad&\text{for}\quad \spa(2n), \tag 3.8\\
-2\mu_1\in\ZZ_+ \qquad&\text{for}\quad \oa(2n+1). \tag 3.9
\endalign
$$
The homomorphism (1.23) allows us to extend it to
a $\Y^{\pm}(N)$-module. This module is obviously highest weight with
$$
\mu_i(u)=\frac{1+(\mu_i\pm 1/2)u^{-1}}{1\pm1/2u^{-1}},
\qquad i=1,\dots, n,
\tag 3.10
$$
and $\mu_0(u)=1$ if $N=2n+1$.

Using (1.31) we can
equip any tensor product $L\ot V$ of a $\Y(N)$-module $L$
and a $\Y^{\pm}(N)$-module $V$
with a $\Y^{\pm}(N)$-action by the rule
$$
y\cdot (\xi\ot\eta):=\Delta(y)(\xi\ot\eta), 
\qquad y\in\Y^{\pm}(N),\quad \xi\in L,\quad \eta\in V. \tag 3.11
$$
In particular, we can construct representations of $\Y^{\pm}(N)$
of the form (0.3).
\medskip

Decomposition (1.30) enables one to establish a correspondence
between representations of the twisted Yangian
$\Y^{\pm}(N)$ and the special twisted Yangian $\SY^{\pm}(N)$
analogous to the $A$ series case (see Section 2).
Consider the {\it similarity classes\/} of 
finite-dimensional irreducible representations of $\Y^{\pm}(N)$:
two representations belong to the same class if
one can be obtained from the other by the composition
with an automorphism (1.29). Since the subalgebra $\SY^{\pm}(N)$
is stable with respect to all automorphisms of the form (1.29),
any representation
from such a class, restricted to the subalgebra $\SY^{\pm}(N)$
gives the same representation of the latter. 
Moreover, all finite-dimensional irreducible representations
of $\SY^{\pm}(N)$ can be obtained in this way.
\medskip

In the next three sections we consider separately the twisted Yangians
corresponding to the series C, D, and B, respectively. We obtain
a description of finite-dimensional irreducible representations
for the algebras $\Y^{\pm}(N)$ and $\SY^{\pm}(N)$.

\bigskip
\bigskip
\noindent
{\bf 4. Representations of $\Y^-(2n)$}
\bigskip

As for the Yangian $\Y(N)$ the representations of the simplest
twisted Yangian $\Y^-(2)$ play an important role in the general case.
\medskip

Irreducible highest weight representations $V(\mu(u))$ of $\Y^-(2)$
are parameterized by formal series $\mu(u)\in 1+u^{-1}\C[[u^{-1}]]$.
The following is an analogue of Proposition 2.5 for the algebra $\Y^-(2)$.

\bigskip
\proclaim
{\bf Proposition 4.1} If $\ts\dim V(\mu(u)) <\infty$ then there exists
an even formal series $\psi(u)\in 1+u^{-2}\C[[u^{-2}]]$ such that
$\psi(u)\mu(u)$ is a polynomial in $u^{-1}$.
\endproclaim

\Proof The symmetry relation (1.18) implies that
$$
(2u-1)s_{1,-1}(u)=(-2u-1)s_{1,-1}(-u),
$$
and so, for any $r\geq 1$
$$
s_{1,-1}^{(2r-1)}=2s_{1,-1}^{(2r)}.
\tag 4.1
$$
Denote by $\xi$ the highest weight vector of $V(\mu(u))$.
Take the minimum nonnegative integer $k$ such that
$s_{1,-1}^{(2k+2)}\xi$ is a linear combination of the vectors
$s_{1,-1}^{(2)}\xi,\dots,s_{1,-1}^{(2k)}\xi$; this $k$ exists because
$V(\mu(u))$ is finite-dimensional. Let us prove that for any $r\geq k+1$
we have
$$
s_{1,-1}^{(2r)}\xi=
a_1^{(2r)}\xi_1+\cdots+a_k^{(2r)}\xi_k
\tag 4.2
$$
for some complex coefficients $a_i^{(2r)}$, where 
$\xi_i:=s_{1,-1}^{(2i)}\xi$,
$i=1,\dots,k$.
This is true for $r=k+1$
by our choice of $k$. Taking in (1.17) $i=k=l=1$, $j=-1$ 
and exchanging $u$ and $v$ we get
$$
\aligned
[s_{11}(u),s_{1,-1}(v)]={}&\frac{2u+1}{u^2-v^2}s_{1,-1}(v)s_{11}(u)\\
-{}&\frac{u+v+1}{u^2-v^2}s_{1,-1}(u)s_{11}(v)
+\frac{1}{u+v}s_{1,-1}(u)s_{-1,-1}(v).
\endaligned
\tag 4.3
$$
Take the coefficient at $u^{-3}v^{-2p}$ in this relation.
Using (4.1) and the fact that the series $s_{11}(u)+s_{-1,-1}(u)$
is even in $u$ (see (1.18)), we obtain
$$
[s_{11}^{(3)},s_{1,-1}^{(2p)}]=2s_{1,-1}^{(2p+2)}
+s_{1,-1}^{(2p)}(2s_{11}^{(2)}+s_{11}^{(1)}) 
-s_{1,-1}^{(2)}(3s_{11}^{(2p)}+s_{-1,-1}^{(2p)}).
\tag 4.4
$$
Write $\mu(u)=1+\mu^{(1)}u^{-1}+\mu^{(2)}u^{-2}+\cdots$. Then
applying both sides of (4.4) to $\xi$ we get
$$
2s_{1,-1}^{(2p+2)}\xi=s_{11}^{(3)}s_{1,-1}^{(2p)}\xi
-(\mu^{(3)}+2\mu^{(2)}+\mu^{(1)}) s_{1,-1}^{(2p)}\xi+
(3\mu^{(2p)}+\mu_-^{(2p)}) s_{1,-1}^{(2)}\xi,
\tag 4.5
$$
where $\mu_-^{(r)}$ is the eigenvalue of $s_{-1,-1}^{(r)}$ at $\xi$.
In particular, for $p=1,\dots,k-1$ we have
$$
s_{11}^{(3)}\xi_p=2\xi_{p+1}+(\mu^{(3)}+2\mu^{(2)}+\mu^{(1)}) \xi_p
-(3\mu^{(2p)}+\mu_-^{(2p)})  \xi_1.
\tag 4.6
$$
Applying (4.2) with $r=k+1$ gives
$$
\aligned
s_{11}^{(3)}\xi_k={}&2(a_1^{(2k+2)}\xi_1+\cdots+a_k^{(2k+2)}\xi_k)\\
+{}&(\mu^{(3)}+2\mu^{(2)}+\mu^{(1)}) \xi_k
-(3\mu^{(2k)}+\mu_-^{(2k)}) \xi_1.
\endaligned
\tag 4.7
$$
Now (4.2) follows from (4.5)--(4.7) by an induction on $r$.

Rewriting (4.2) in terms of generating series and using (4.1) we obtain
$$
s_{1,-1}(u)\xi=(2u+1)(a_1(u)\xi_1+\cdots+a_k(u)\xi_k),
\tag 4.8
$$
where the $a_i(u)$ are even series in $u^{-1}$:
$a_i(u)=u^{-2i}+a_i^{(2k+2)}u^{-2k-2}+\cdots$.

Taking the coefficient at $u^{-3}$ in (4.3) and replacing $v$ with $u$
we get
$$
\align
s_{11}^{(3)}s_{1,-1}(u)={}&s_{1,-1}(u)
(s_{11}^{(3)}+2s_{11}^{(2)}+s_{11}^{(1)}+2u^2)\\
{}-{}&(s_{1,-1}^{(2)}+(u+1)s_{1,-1}^{(1)})s_{11}(u)
+(s_{1,-1}^{(2)}-us_{1,-1}^{(1)})s_{-1,-1}(u).
\endalign
$$
Apply both sides to the vector $\xi$. Then using (4.8) and the symmetry
relation (1.18) we find that
$$
\aligned
(s_{11}^{(3)}-\mu^{(3)}-2\mu^{(2)}-\mu^{(1)}-2u^2)
(a_1(u)\xi_1+\cdots+a_k(u)\xi_k)&\\
={}{}-\left(\frac{2u+1}{2u}\mu(u)+\frac{2u-1}{2u}\mu(-u)\right)&\xi_1.
\endaligned
\tag 4.9
$$
Using (4.6) and (4.7) and taking the coefficient at $\xi_p$,
$p=2,\dots,k$ in this relation we get
$$
a_{p-1}(u)-u^2a_p(u)+\text{const}\cdot a_k(u)=0.
$$
Therefore, for any $p=1,\dots,k$ one has $a_p(u)=A_p(u)a_k(u)$
for a monic polynomial $A_p(u)$ of degree $k-p$ in $u^2$. Hence,
taking the coefficient at $\xi_1$ in (4.9) we obtain that
$$
\left(\frac{2u+1}{2u}\mu(u)+\frac{2u-1}{2u}\mu(-u)\right)=B(u)a_k(u)
\tag 4.10
$$
for a polynomial $B(u)$ of degree $k$ in $u^2$.

Similarly, taking the coefficient at $u^{-2}$ in (4.3) we get
$$
s_{11}^{(2)}s_{1,-1}(u)=s_{1,-1}(u)
(s_{11}^{(2)}+2s_{11}^{(1)}+1)
-s_{1,-1}^{(1)}(s_{11}(u)-s_{-1,-1}(u)).
\tag 4.11
$$
Applying both sides to $\xi$, using (4.8) and the symmetry
relation (1.18) we obtain
$$
(s_{11}^{(2)}-\mu^{(2)}-2\mu^{(1)}-1)
(a_1(u)\xi_1+\cdots+a_k(u)\xi_k)
=-u^{-1}(\mu(u)-\mu(-u))\xi_1.
\tag 4.12
$$
In particular, $s_{11}^{(2)}\xi_p$ is a linear combination of $\xi_p$ and
$\xi_1$. Finally, taking the coefficient at $\xi_1$ in (4.12) we find that
$u^{-1}(\mu(u)-\mu(-u))=C(u)a_k(u)$ for a polynomial $C(u)$ of degree
$k-1$ in $u^2$. Together with (4.10) this proves that the
desired series is $\psi(u)=(u^{2k}a_{k}(u))^{-1}$. $\square$
\medskip

We shall consider (2.16) as a representation
of $\Y^-(2)$
obtained by restriction. Equivalently, using (3.11) one may regard
this representation as the tensor
product of $\Y(2)$-module (2.16) with the trivial $\Y^-(2)$-module.

Let $\mu(u)$ be a polynomial in $u^{-1}$ of degree ${}\leq 2k$.
Consider an expansion
$$
\mu(u)=(1-\gamma_1u^{-1})\cdots (1-\gamma_{2k}u^{-1}).
\tag 4.13
$$
Re-enumerating the numbers $\gamma_i$ if necessary,
we may assume that the following condition is satisfied (cf. (2.19)):
$$
\aligned
&\text{for every}\ \  i=1,\dots, k\ \  \text{we have: } \\
&\text{if the set}\ \ \{\gamma_p+\gamma_q\ |\ 2i-1\leq p<q \leq 2k\}
\cap \ZZ_+ 
\ \ \text{is not empty}\\
&\text{then}\ \  \gamma_{2i-1}+\gamma_{2i}\ \  
\text{is its minimal element.}
\endaligned
\tag 4.14
$$

The following is an analogue of Proposition 2.6.

\bigskip
\proclaim
{\bf Proposition 4.2} Let $V(\mu(u))$ 
be an irreducible highest weight representation
of $\Y^-(2)$ where $\mu(u)$ is a polynomial
of degree ${}\leq 2k$ in $u^{-1}$. Suppose that decomposition (4.13) is
given and condition (4.14) is satisfied.
Then $V(\mu(u))$ is isomorphic to the tensor
product
$$
L(\gamma_1,-\gamma_2)\ot\cdots\ot L(\gamma_{2k-1},-\gamma_{2k}),
\tag 4.15
$$
regarded as $\Y^-(2)$-module.
\endproclaim

\Proof Denote representation (4.15) by $V$.
Let $\xi_i$ be the highest weight vector of 
$L(\gamma_{2i-1},-\gamma_{2i})$. We obtain from (2.23) that
the vector $\xi:=\xi_1\ot\cdots\ot \xi_k$ satisfies
the relations
$$
\aligned
t_{-1,1}(u)\xi&=0,\\
t_{11}(u)\xi&=(1-\gamma_{2} u^{-1})\cdots (1-\gamma_{2k} u^{-1})\xi,\\
t_{-1,-1}(u)\xi_i&=(1+\gamma_{1} u^{-1})
\cdots (1+\gamma_{2k-1} u^{-1})\xi.
\endaligned
\tag 4.16
$$
Hence, by (1.15)
$$
s_{11}(u)\xi=\mu(u)\xi\qquad\text{and}\qquad
s_{-1,1}(u)\xi=0.\tag 4.17
$$ 
Thus, the cyclic  $\Y^-(2)$-span of the vector $\xi$ in $V$ is
a highest weight module with the highest weight 
$\mu(u)$.
Therefore, to prove Proposition 4.2 it suffices to show that
the representation $V$ of $\Y^-(2)$ is irreducible. We use induction on
$k$.

For $k=1$ this is obvious, because the restriction of 
$L(\gamma_{2i-1},-\gamma_{2i})$ to $\spa(2)$ is irreducible
and so its extension to $\Y^-(2)$ is.

Suppose now that $k>1$. If $\wt V$ is a nonzero submodule of $V$
then $\wt V$ contains a nonzero vector $\eta$ such that
$s_{-1,1}(u)\eta=0$. Write $\eta$ in the form
$$
\eta=\sum_{r=0}^p(E_{1,-1})^r \xi_1\ot \eta_r,\tag 4.18
$$
for some $p\geq 0$, 
where 
$\eta_r\in L(\gamma_3,-\gamma_4)\ot\cdots\ot
L(\gamma_{2k-1},-\gamma_{2k})$; 
and
if $\gamma_1+\gamma_2\in\ZZ_+$ then $p\leq  \gamma_1+\gamma_2$. We may
assume that $\eta_p\ne 0$. By (1.31) we have
$$
\aligned
s_{-1,1}(u)\left((E_{1,-1})^r \xi_1\ot \eta_r\right)=
-t_{-1,-1}(u)t_{-1,1}(-u)(E_{1,-1})^r \xi_1&\ot s_{-1,-1}(u)\eta_r\\
+t_{-1,-1}(u)t_{-1,-1}(-u)(E_{1,-1})^r \xi_1&\ot s_{-1,1}(u)\eta_r\\
-t_{-1,1}(u)t_{-1,1}(-u)(E_{1,-1})^r \xi_1&\ot s_{1,-1}(u)\eta_r\\
+t_{-1,1}(u)t_{-1,-1}(-u)(E_{1,-1})^r \xi_1&\ot s_{1,1}(u)\eta_r.
\endaligned
\tag 4.19
$$
Rewriting (2.25) and (2.26) we obtain
$$
t_{-1,-1}(u)(E_{1,-1})^r \xi_1=(1+(\gamma_1-r)u^{-1})(E_{1,-1})^r \xi_1
\tag 4.20
$$
and
$$
t_{-1,1}(u)(E_{1,-1})^r \xi_1=u^{-1}\ts r(\gamma_1+\gamma_2-r+1)
(E_{1,-1})^{r-1} \xi_1.
\tag 4.21
$$
Taking the coefficient at $(E_{1,-1})^p \xi_1$ in the relation
$s_{-1,1}(u)\eta=0$ gives
$$
(1-(\gamma_1-p)u^{-2})s_{-1,1}(u)\eta_p=0.
$$
So, $s_{-1,1}(u)\eta_p=0$. 
By the induction
hypotheses, the representation 
$L(\gamma_3,-\gamma_4)\ot\cdots\ot L(\gamma_{2k-1},-\gamma_{2k})$ 
of $\Y^-(2)$ is irreducible, and hence
$$
\eta_p=\ \text{const}\cdot \xi_2\ot\cdots\ot\xi_k,\tag 4.22
$$
where `const' is a nonzero constant; see Remark 3.4. 

Suppose now that $p\geq 1$ and take the coefficient at 
$(E_{1,-1})^{p-1} \xi_1$ in the relation $s_{-1,1}(u)\eta=0$. 
Using (4.19)--(4.21) we obtain
$$
\aligned
(1-(\gamma_1-p&+1)^2u^{-2})s_{-1,1}(u)\eta_{p-1}\\
{}+{}&u^{-1}\ts p(\gamma_1+\gamma_2-p+1)(1+(\gamma_1-p+1)u^{-1}) 
s_{-1,-1}(u)\eta_p\\
{}+{}&u^{-1}\ts p(\gamma_1+\gamma_2-p+1)(1-(\gamma_1-p)u^{-1}) 
s_{11}(u)\eta_p=0.
\endaligned
\tag 4.23
$$
The symmetry relation (1.18) gives
$$
s_{-1,-1}(u)=\frac{2u+1}{2u} s_{11}(-u)-\frac{1}{2u} s_{11}(u)
\tag 4.24
$$
and
$$
s_{-1,1}^{(2r-1)}=2s_{-1,1}^{(2r)},\qquad r\geq 1.
\tag 4.25
$$
Using (4.24) and multiplying (4.23) by $u^{2k}/(2u+1)$ we bring
it to the form
$$
\aligned
&(u^2-(\gamma_1-p+1)^{2})u^{2k-2}\ts\frac{s_{-1,1}(u)}{2u+1}\eta_{p-1}
+p(\gamma_1+\gamma_2-p+1)\\
&\times\frac{1}{2u}
\left((u+\gamma_1-p+1) u^{2k-2}s_{11}(-u)
+(u-\gamma_1+p-1) u^{2k-2}s_{11}(u)\right)\eta_p=0.
\endaligned
\tag 4.26
$$
Due to (2.18) and (1.15) the elements $s_{ij}^{(r)}$ with $r>2k-2$
act as zero operators
in $L(\gamma_3,-\gamma_4)\ot\cdots\ot L(\gamma_{2k-1},-\gamma_{2k})$.
Hence, the operators $u^{2k-2}s_{ij}(u)$ are polynomials in $u$, 
as well as $u^{2k-2}s_{-1,1}(u)/(2u+1)$ is, because of (4.25).
Further, by (4.22) and (4.16)
$$
u^{2k-2}s_{11}(u)\eta_p=(u-\gamma_3)\cdots (u-\gamma_{2k})\eta_p.
$$
Putting $u=\gamma_1-p+1$ in (4.26) we get
$$
p(\gamma_1+\gamma_2-p+1)(\gamma_1+\gamma_3-p+1)
\cdots (\gamma_1+\gamma_{2k}-p+1)=0.
$$
However, this
is impossible because of condition (4.14).
Therefore, $p$ has to be equal to $0$ and so, $\eta$ is
a multiple of $\xi$.

To complete the proof we have to show that the submodule of $V$
generated by $\xi$ coincides with $V$.

The representation of $\Y(2)$ dual to (4.15) is isomorphic to
$$
L'(-\gamma_1,\gamma_2)\ot\cdots\ot L'(-\gamma_{2k-1},\gamma_{2k});
\tag 4.27
$$
see (2.31). On the other hand, 
using the antiautomorphism (1.26) we can equip the vector space $V^*$
dual to $V$ with an action of $\Y^-(2)$:
$$
(y\cdot f)(v)=f(\sigma(y)v),\qquad y\in\Y^-(2),\quad f\in V^*,\quad v\in V.
$$
Since the antiautomorphism (1.26) is a restriction of (1.7), 
the representation $V^*$ is isomorphic to the restriction
of (4.27) to $\Y^-(2)$.

Similarly to the case of representation $V$, any nonzero submodule of $V^*$
contains a nonzero vector $\eta'$ such that $s_{1,-1}(u)\eta'=0$.
Modifying the previous argument for lowest weight representations
of $\Y^-(2)$
we prove that $\eta'$ coincides (up to a nonzero factor) with
the vector $\xi'=\xi'_1\ot\cdots\ot \xi'_k$,
where $\xi'_i$ is the lowest weight vector of 
$L'(-\gamma_{2i-1},\gamma_{2i})$.
Now, if the submodule in $V$ generated by $\xi$ is proper, then
its annihilator in $V^*$ is a nonzero submodule which does not
contain $\xi'$. Contradiction. $\square$
\medskip

We are now able to give necessary and sufficient conditions
for $V(\mu(u))$ to be finite-dimensional (cf. Proposition 2.7).

\bigskip
\proclaim
{\bf Proposition 4.3} Let $V(\mu(u))$ be an irreducible representation
of $\Y^-(2)$ with the highest weight 
$\mu(u)$. Then $V(\mu(u))$ is
finite-dimensional if and only if there exists a formal series
$\psi(u)\in 1+u^{-2}\C[[u^{-2}]]$ such that 
$\psi(u)\mu(u)$ is a polynomial
in $u^{-1}$ with a decomposition
$$
\psi(u)\mu(u)=(1-\gamma_1u^{-1})\cdots (1-\gamma_{2k}u^{-1}),
\tag 4.28
$$
and $\gamma_{2i-1}+\gamma_{2i}\in\ZZ_+$, $i=1,\dots,k$ 
for a certain re-enumeration
of the $\gamma_i$.
\endproclaim

\Proof Let $\dim V(\mu(u))<\infty$.
By Proposition 4.1,  
taking the composition of our representation with
an automorphism of the form (1.29) we may assume that (4.28)
is satisfied for a certain series $\psi(u)$. Re-enumerate the $\gamma_i$
in such a way that condition (4.14)
be satisfied. Now use Proposition 4.2. Since representation
(4.15) is finite-dimensional this gives the relations
$\gamma_{2i-1}+\gamma_{2i}\in\ZZ_+$ for all $i=1,\dots,k$.

Conversely, let (4.28) hold with 
$\gamma_{2i-1}+\gamma_{2i}\in\ZZ_+$, $i=1,\dots,k$. 
Then representation
(4.15) is finite-dimensional and the composition of $V(\mu(u))$
with the automorphism (1.29) corresponding to the series $\psi(u)$
is isomorphic to the irreducible
quotient of the $\Y^-(2)$-cyclic span in $V$ of the vector $\xi$;
see (4.17).
Therefore, $\dim V(\mu(u))<\infty$. $\square$
\medskip

Proposition 4.3 can be reformulated with the use of monic
polynomials in $u$ as follows (cf. Theorem 2.8).

\bigskip
\proclaim
{\bf Theorem 4.4} The irreducible highest weight representation
$V(\mu(u))$ 
of $\Y^-(2)$ is
finite-dimensional if and only if there exists a monic
polynomial $P(u)\in\C[u]$ such that
$
P(u)=P(-u+1)
$
and
$$
\frac{\mu(-u)}{\mu(u)}=\frac{P(u+1)}{P(u)}.\tag 4.29
$$
In this case $P(u)$ is unique.
\endproclaim

\Proof If $\dim V(\mu(u))<\infty$ then by 
Proposition 4.3 relation (4.28) holds for a certain even
series $\psi(u)$ and
$\gamma_{2i-1}+\gamma_{2i}\in\ZZ_+$ for all $i$.
Set
$$
\aligned
\sigma(u)&=(1+\gamma_{1} u^{-1})(1+\gamma_{3} u^{-1})
\cdots (1+\gamma_{2k-1} u^{-1}),\\
\tau(u)&=(1-\gamma_{2} u^{-1})(1-\gamma_{4} u^{-1})
\cdots (1-\gamma_{2k} u^{-1}).
\endaligned
\tag 4.30
$$
Then there exists a monic polynomial $Q(u)$ in $u$
such that
$$
\frac{\sigma(u)}{\tau(u)}=\frac{Q(u+1)}{Q(u)};\tag 4.31
$$
see the proof of Theorem 2.8. Define $P(u)$
by
$$
P(u)=Q(u)Q(-u+1)(-1)^{\deg Q}.\tag 4.32
$$
Then $P(u)$ is monic and $P(u)=P(-u+1)$. Moreover,
$$
\frac{\tau(-u)}{\sigma(-u)}=\frac{Q(-u)}{Q(-u+1)}.\tag 4.33
$$
Therefore, since $\psi(u)\mu(u)=\sigma(-u)\tau(u)$, 
relation (4.29) holds
due to (4.31) and (4.33).

Conversely, let (4.29) be satisfied for
a monic polynomial $P(u)$ such that $P(u)=P(-u+1)$.
Then it has an even degree and the multiset of 
its roots is invariant with respect to
the transformation $u\mapsto -u+1$ (the symmetry in $\C$ with
the center at 1/2). Therefore, we may write the roots in the form
$\{-\delta_1,\dots,-\delta_s, \delta_1+1,\dots,\delta_s+1\}$ for
a certain $s$. Set
$$
\mu'(u)=(1-(\delta_1+1)u^{-1})\cdots (1-(\delta_s+1) u^{-1})
(1+\delta_1u^{-1})\cdots (1+\delta_s u^{-1}).
$$
By Proposition 4.3 the representation $V(\mu'(u))$
is finite-dimensional. Moreover, the polynomial $\mu'(u)$ satisfies
(4.29) and so, the ratio $\mu(u)/\mu'(u)$ is an even formal series
in $u^{-1}$. This proves that $V(\mu(u))$ is also finite-dimensional.

The uniqueness of $P(u)$ is established 
in the proof of Theorem 2.8.
$\square$
\medskip

The following notation is motivated by Theorem 4.4 (cf. (2.35)). 
For two
formal series $\mu(u)$ and $\nu(u)$ in $u^{-1}$ we shall write
$$
\mu(u)\Ra \nu(u),\tag 4.34
$$
if there exists a monic polynomial $P(u)$ in $u$ such that
$P(u)=P(-u+1)$ and
$$
\frac{\mu(u)}{\nu(u)}=\frac{P(u+1)}{P(u)}.
$$
Obviously, this is only possible if $\mu(u)$ and $\nu(u)$
satisfy the relation $\mu(u)\mu(-u)=\nu(u)\nu(-u)$.

Proposition 4.2 and Theorem 4.4 imply the following parameterization
of representations of the special twisted
Yangian $\SY^-(2)$ (cf. Corollary 2.9). 

\bigskip
\proclaim
{\bf Corollary 4.5} There is a one-to-one correspondence
between finite-dimensional irreducible representations
of the special twisted Yangian $\SY^-(2)$ and the monic polynomials 
$P(u)$ in $u$ such that $P(u)=P(-u+1)$.
Every such representation is isomorphic to a representation of the form 
(4.15). $\square$
\endproclaim

\bigskip
\noindent
{\it Remark 4.6.} Tensoring $V(\mu(u))$ by factors of the form 
$L(\gamma,\gamma)$ is equivalent to multiplying
the highest weight $\mu(u)$ by $1-\gamma^2u^{-2}$.
However, the representations $V(\mu(u))$ and $V((1-\gamma^2u^{-2})\mu(u))$
of $\Y^-(2)$ belong to the same similarity class.
Therefore, we might add to Corollary 4.5 the condition that in (4.15)
all $\gamma_{2i-1}+\gamma_{2i}$ are positive integers. $\square$
\medskip

The following corollary provides a criterion of irreducibility
for the representations of $\Y^-(2)$ (or $\SY^-(2)$)
of the form (4.15) (cf. Corollary 2.11). We formulate it
in terms of the strings $S(\alpha,\beta)$ introduced in Section 2.
Note that two strings $S(\alpha,-\beta)$ and $S(\beta,-\alpha)$
are symmetric to each other with respect to the center $-1/2$
in $\C$. Therefore, $S(\alpha_1,-\beta_1)$ and $S(\alpha_2,-\beta_2)$
are in general position whenever 
$S(\beta_1,-\alpha_1)$ and $S(\beta_2,-\alpha_2)$ are.

\bigskip
\proclaim
{\bf Corollary 4.7}  Consider an irreducible finite-dimensional
representation of $\Y(2)$ of the form
$$
V=L(\gamma_1,-\gamma_2)\ot\cdots\ot L(\gamma_{2k-1},-\gamma_{2k}).
\tag 4.35
$$
Then its restriction to $\Y^-(2)$ is irreducible
if and only if for each pair $i< j$ the strings
$S(\gamma_{2i},-\gamma_{2i-1})$ and $S(\gamma_{2j-1},-\gamma_{2j})$
are in general position.
\endproclaim

\Proof By Corollary 2.11, the strings $S(\gamma_{2i-1},-\gamma_{2i})$,
$i=1,\dots,k$
are pairwise in general position. Permuting the tensor factors in
(4.35) if necessary, we may assume that
$
\gamma_1+\gamma_2\leq \cdots \leq \gamma_{2k-1}+\gamma_{2k}.
$
Suppose that 
the strings $S(\gamma_{2i},-\gamma_{2i-1})$ and
$S(\gamma_{2j-1},-\gamma_{2j})$ are in general position
for all $i <j$.
Then one easily checks that (4.14) is satisfied and so,
by Proposition 4.2 representation (4.35) of $\Y^-(2)$
is irreducible.

Conversely, let $V$ be irreducible as a $\Y^-(2)$-module. Consider
the module $V'$ obtained by replacement of some of the tensor factors
$L(\gamma_{2i-1},-\gamma_{2i})$ respectively by 
$L(\gamma_{2i},-\gamma_{2i-1})$. Then $V'$ is isomorphic to $V$.
Indeed, $\dim V=\dim V'$ and one easily checks using (1.15)
that $V$ is isomorphic
to the irreducible quotient of the $\Y^-(2)$-submodule in $V'$
generated by the tensor product of the highest weight vectors of
the factors. So, $V'$ is irreducible as a $\Y^-(2)$-module and
hence, as a $\Y(2)$-module. By Corollary 2.11 the corresponding strings
are pairwise in general position. $\square$
\medskip

The following theorem together with Theorem 3.3
gives a description
of finite-dimensional irreducible representations of $\Y^-(2n)$
(cf. Theorem 2.12). We use notation (2.35) and (4.34).

\bigskip
\proclaim
{\bf Theorem 4.8} The irreducible highest weight representation
$V(\mu(u))$, $\mu(u)=(\mu_1(u),\dots,\mu_n(u))$ of 
$\Y^-(2n)$ is finite-dimensional if and only if
the following relation holds
$$
\mu_1(-u)\Ra \mu_1(u)\to \mu_2(u)\to\cdots\to\mu_n(u).\tag 4.36
$$
\endproclaim

\Proof Suppose that $\dim V(\mu(u))<\infty$. By Proposition 3.5
we have
$$
\mu_1(u)\to \mu_2(u)\to\cdots\to\mu_n(u). \tag 4.37
$$
The subalgebra in $\Y^-(2n)$ 
generated by the coefficients of the series 
$s_{ij}(u)$ with $i,j=-1,1$
is isomorphic to
$\Y^-(2)$. The cyclic span of the highest weight vector of $V(\mu(u))$
with respect to this subalgebra is a representation with
the highest weight $\mu_1(u)$. Its irreducible quotient
is finite-dimensional and so, by Theorem 4.4 we have
$\mu_1(-u)\Ra \mu_1(u)$.

Conversely, let (4.36) hold. Then
$$
\frac{\mu_1(-u)}{\mu_1(u)}=\frac{P_1(u+1)}{P_1(u)},
\tag 4.38
$$
and
$$
\frac{\mu_{i-1}(u)}{\mu_{i}(u)}=\frac{P_i(u+1)}{P_i(u)} 
\tag 4.39
$$
for certain polynomials $P_i(u)$ of the form
$$
P_i(u)=(u+\delta_1^{(i)})\cdots (u+\delta_{s_i}^{(i)}),
\qquad i=2,\dots,n
\tag 4.40
$$
and
$$
P_1(u)=(u+\delta_1^{(1)})\cdots (u+\delta_{s_1}^{(1)})
(u-\delta_1^{(1)}-1)\cdots (u-\delta_{s_1}^{(1)}-1);
\tag 4.41
$$
see the proof of Theorem 4.4.

Consider the irreducible highest weight representation
$L(\lambda(u))$ with $\lambda(u)=(\lambda_{-n}(u),\dots,\lambda_n(u))$
of the Yangian $\Y(2n)$ where for $i=1,\dots, n$
$$
\aligned
\lambda_i(u)=&\prod_{a=1}^{i}(1+\delta_1^{(a)}u^{-1})
\cdots(1+\delta_{s_a}^{(a)}u^{-1})
\\
\times&\prod_{a=i+1}^{n}(1+(\delta_1^{(a)}+1)u^{-1})
\cdots(1+(\delta_{s_a}^{(a)}+1)u^{-1})
\endaligned
\tag 4.42
$$
and
$$
\lambda_{-i}(u)=
\prod_{a=1}^{n}(1+(\delta_1^{(a)}+1)u^{-1})
\cdots(1+(\delta_{s_a}^{(a)}+1)u^{-1}). \tag 4.43
$$
Then $L(\lambda(u))$ is finite-dimensional by Theorem 2.12.
One easily derives from relations (1.1) and (1.15) that 
the cyclic $\Y^-(2n)$-span of its highest weight vector
is a representation with the highest weight
$\mu'(u)=(\mu'_1(u),\dots,\mu'_n(u))$ where
$\mu'_i(u)=\lambda_i(u)\lambda_{-i}(-u)$. So,
the representation $V(\mu'(u))$ of $\Y^-(2n)$
is finite-dimensional. However, the polynomials $\mu'_i(u)$
satisfy (4.38) and (4.39). Hence, there exists an
automorphism of $\Y^-(2n)$ of the form (1.29)
such that its composition with the representation
$V(\mu'(u))$ is isomorphic to
$V(\mu(u))$ which completes the proof. $\square$
\medskip

Theorem 4.8 implies
the following parameterization of representations of the special twisted
Yangian $\SY^-(2n)$
(cf. Corollary 2.13).

\bigskip
\proclaim
{\bf Corollary 4.9} There is a one-to-one correspondence
between finite-dimensional irreducible representations
of the special twisted Yangian $\SY^-(2n)$
and the families $\{P_1(u),\dots,P_{n}(u)\}$ of
monic polynomials in $u$ with $P_1(u)=P_1(-u+1)$.
Every such representation is
isomorphic to a subquotient of a representation 
of the form $L(\lambda(u))$. $\square$
\endproclaim

\bigskip
\bigskip
\noindent
{\bf 5. Representations of $\Y^+(2n)$}
\bigskip

The results for the $D$ series are very similar to the
previous case. However, some differences
take place for $n=1$ since the
twisted Yangian $\Y^+(2)$ corresponds to the abelian Lie algebra
$\oa(2)$. In particular, contrary to the case of $\Y^-(2)$,
it is not true that
every finite-dimensional irreducible representations of
$\Y^+(2)$ is obtained by restriction of a representation of
$\Y(2)$. Also, the step from $\Y^+(2)$ to $\Y^+(2n)$, $n\geq 2$
requires an extra care.

We omit the most of the proofs in this section which are parallel 
to the case of
$\Y^-(2)$ but provide all necessary additional arguments.  

Irreducible highest weight representations $V(\mu(u))$ of $\Y^+(2)$
are parameterized by formal series $\mu(u)\in 1+u^{-1}\C[[u^{-1}]]$.

\bigskip
\proclaim
{\bf Proposition 5.1} If $\ts\dim V(\mu(u)) <\infty$ then there exists
an even formal series $\psi(u)\in 1+u^{-2}\C[[u^{-2}]]$ such that
$(1+1/2u^{-1})\psi(u)\mu(u)$ is a polynomial in $u^{-1}$.
\endproclaim

\Proof We repeat the argument of the proof of Proposition 4.1 
almost word by word and only
point out a few differences.
First, here
the symmetry relation (1.18) gives
$$
s_{1,-1}(u)=s_{1,-1}(-u),
$$
and so, for any $r\geq 1$
$$
s_{1,-1}^{(2r-1)}=0.
$$ 
Relation (4.3) is replaced by
$$
\aligned
[s_{11}(u),s_{1,-1}(v)]={}&\frac{2u+1}{u^2-v^2}s_{1,-1}(v)s_{11}(u)\\
-{}&\frac{u+v+1}{u^2-v^2}s_{1,-1}(u)s_{11}(v)
-\frac{1}{u+v}s_{1,-1}(u)s_{-1,-1}(v).
\endaligned
\tag 5.1
$$
Instead of (4.8) we obtain here
$$
{s}_{1,-1}(u)\xi=a_1(u)\xi_1+\cdots+a_k(u)\xi_k.
$$
The counterpart of relation (4.10) has exactly the same form.
However, for the remaining part of the proof instead of (4.11)
we need to consider
the following formula obtained from (5.1) by taking the coefficient at 
$u^{-4}$:
$$
\aligned
s_{11}^{(4)}s_{1,-1}(u)={}&s_{1,-1}(u)
(s_{11}^{(4)}+2s_{11}^{(3)}+s_{11}^{(2)}+2s_{11}^{(1)}u^2+u^2)\\
{}-{}&s_{1,-1}^{(2)}((u+1)s_{11}(u)-us_{-1,-1}(u)).
\endaligned
$$
Applying both sides to $\xi$ and using the symmetry relation (1.18)
we find that 
$$
\left(u+\frac12\right)\mu(u)-\left(u-\frac12\right)\mu(-u)
=C(u)a_k(u)
$$
for a polynomial $C(u)$ of degree $k$ in $u^2$. Together with
(4.10) this proves that for $\psi(u)=(u^{2k}a_{k}(u))^{-1}$
the series $(1+1/2u^{-1})\psi(u)\mu(u)$ is a polynomial in $u^{-1}$
of degree ${}\leq 2k+1$.
$\square$
\medskip

Note that the irreducible representation of $\Y^+(2)$ with the
highest weight $\mu(u)=(1+\gamma u^{-1})(1+1/2 u^{-1})^{-1}$, where
$\gamma\in\C$ is
one-dimensional. It is obtained by extending the one-dimensional
representation of the abelian Lie algebra $\oa(2)$ given by
$
F_{11}\ts\xi=(\gamma-1/2)\ts\xi, 
$
where $\xi$ is the basis vector; see (3.10). 
We shall denote this representation of $\Y^+(2)$ by $V(\gamma)$.

Let $(1+1/2 u^{-1})\mu(u)$ be polynomial in $u^{-1}$ of 
degree ${}\leq 2k+1$. Consider an expansion
$$
\mu(u)=(1-\gamma_1u^{-1})\cdots (1-\gamma_{2k+1}u^{-1})
\left(1+\frac12 u^{-1}\right)^{-1}.
\tag 5.2
$$
Re-enumerating the numbers $\gamma_i$ if necessary,
we may assume that the following condition is satisfied (cf. (4.14)):
$$
\aligned
&\text{for every}\ \  i=1,\dots, k\ \  \text{we have: } \\
&\text{if the set}\ \ \{\gamma_p+\gamma_q\ |\ 2i-1\leq p<q \leq 2k+1\}
\cap \ZZ_+ 
\ \ \text{is not empty}\\
&\text{then}\ \  \gamma_{2i-1}+\gamma_{2i}\ \  
\text{is its minimal element.}
\endaligned
\tag 5.3
$$
In the following proposition we use tensor products of
representations constructed in accordance with (3.11). 
As before, we denote by
$L(\alpha,\beta)$ a highest weight representation of the Lie algebra
$\gl(2)$ which is also regarded as a $\Y(2)$-module.

\bigskip
\proclaim
{\bf Proposition 5.2} Let $V(\mu(u))$ 
be an irreducible highest weight representation
of $\Y^+(2)$ where $(1+1/2 u^{-1})\mu(u)$ is a polynomial
of degree ${}\leq 2k+1$ in $u^{-1}$. Suppose that decomposition (5.2) is
given and condition (5.3) is satisfied.
Then $V(\mu(u))$ is isomorphic to the tensor
product
$$
L(\gamma_1,-\gamma_2)\ot\cdots\ot L(\gamma_{2k-1},-\gamma_{2k})
\ot V(-\gamma_{2k+1}).
\tag 5.4
$$
\endproclaim

\Proof Denote representation (5.4) by $V$.
Let $\xi_i$ be the highest weight vector of 
$L(\gamma_{2i-1},-\gamma_{2i})$ and $\xi_{k+1}$ be the basis
vector of the one-dimensional representation $V(-\gamma_{2k+1})$. 

Using (2.23) and (1.31) we find that
the vector $\xi:=\xi_1\ot\cdots\ot \xi_{k+1}$ satisfies
the relations
$$
s_{11}(u)\xi=\mu(u)\xi\qquad\text{and}\qquad
s_{-1,1}(u)\xi=0.
$$ 
Thus, the cyclic  $\Y^+(2)$-span of the vector $\xi$ in $V$ is
a highest weight module with the highest weight 
$\mu(u)$.
So, it suffices to show that
the representation $V$ is irreducible. We use induction on $k$.

For $k=0$ this is obvious.
Suppose that $k\geq 1$. Proceeding exactly as in the proof of
Proposition 4.2 we obtain that if the vector $\eta$
is given by (4.18) with
$$
\eta_r\in L(\gamma_3,-\gamma_4)\ot\cdots\ot L(\gamma_{2k-1},-\gamma_{2k})
\ot V(-\gamma_{2k+1})
$$ 
then
$$
\eta_p=\ \text{const}\cdot \xi_2\ot\cdots\ot\xi_{k+1},\tag 5.5
$$
where `const' is a nonzero constant. Further, instead of (4.26) we get here
$$
\aligned
&(u^2-(\gamma_1-p+1)^{2})u^{2k-1}\ts s_{-1,1}(u)\eta_{p-1}
+p(\gamma_1+\gamma_2-p+1)u^{2k-1}\times{}\\
&\left(\frac{2u+1}{2u}(u-\gamma_1+p-1) s_{11}(u)
-\frac{2u-1}{2u}(u+\gamma_1-p+1) s_{11}(-u)\right)\eta_p=0.
\endaligned
\tag 5.6
$$
The elements $s_{-1,1}^{(r)}$ and $s_{1,-1}^{(r)}$ act trivially
in $V(-\gamma_{2k+1})$. Hence, by (1.31) and (3.11) we have
$$
\aligned
s_{-1,1}(u)(\zeta\ot \xi_{k+1})={}&
t_{-1,1}(u)t_{-1,-1}(-u)\zeta\ot s_{11}(u)\xi_{k+1}\\
{}+{}&t_{-1,-1}(u)t_{-1,1}(-u)\zeta\ot s_{-1,-1}(u)\xi_{k+1},
\endaligned
$$
where 
$\zeta\in L(\gamma_3,-\gamma_4)\ot\cdots\ot L(\gamma_{2k-1},-\gamma_{2k})$.
Using the defining relations (1.1) 
and the symmetry relation (1.18) we may rewrite this as follows
$$
\aligned
s_{-1,1}(u)(\zeta\ot \xi_{k+1})={}&
t_{-1,-1}(-u)t_{-1,1}(u)\zeta\ot \frac{2u+1}{2u}s_{11}(u)\xi_{k+1}\\
{}+{}&t_{-1,-1}(u)t_{-1,1}(-u)\zeta\ot 
\frac{2u-1}{2u}s_{11}(-u)\xi_{k+1}.
\endaligned
\tag 5.7
$$
Now, the definition of $V(\gamma)$ and (2.18) imply that
$u^{2k-1}s_{-1,1}(u)\eta_{p-1}$ 
is a polynomial in $u$.
Further, by (5.5) we have
$$
\frac{2u+1}{2u}u^{2k-1}s_{11}(u)\eta_p=
(u-\gamma_3)\cdots (u-\gamma_{2k+1})\eta_p.
$$
So, putting $u=\gamma_1-p+1$ in (5.6) we get
$$
p(\gamma_1+\gamma_2-p+1)(\gamma_1+\gamma_3-p+1)
\cdots (\gamma_1+\gamma_{2k+1}-p+1)=0.
$$
However, this
is impossible because of condition (5.3).
Therefore, $p$ has to be equal to $0$ and so, $\eta$ is
a multiple of $\xi$.

To show that the submodule of $V$
generated by $\xi$ coincides with $V$ we consider
the dual representation $V^*$ defined by
$$
(y\cdot f)(v)=f(\sigma(y)v),\qquad y\in\Y^+(2),\quad f\in V^*,\quad v\in V,
$$
where $\sigma$ is the antiautomorphism of $\Y^+(2)$ given by (1.26).
One easily verifies that $V(\gamma)^*\simeq V(1-\gamma)$ and so,
the representation $V^*$ is isomorphic to
$$
L'(-\gamma_1,\gamma_2)\ot\cdots\ot L'(-\gamma_{2k-1},\gamma_{2k})
\ot V(\gamma_{2k+1}+1),
$$
where $L'(-\gamma_{2i-1},\gamma_{2i})$ is the lowest weight representation
of $\gl(2)$;
see (2.31).

Repeating the previous argument
we prove that any nonzero submodule in $V^*$
contains the tensor product of the lowest weight vectors
of the $L'(-\gamma_{2i-1},\gamma_{2i})$ and the basis
vector of $V(\gamma_{2k+1}+1)$. The proof is completed
exactly as that of Proposition 4.2. $\square$
\medskip

The following proposition is proved in the same way as Proposition 4.3.

\bigskip
\proclaim
{\bf Proposition 5.3} Let $V(\mu(u))$ be an irreducible representation
of $\Y^+(2)$ with the highest weight 
$\mu(u)$. Then $V(\mu(u))$ is
finite-dimensional if and only if there exists a formal series
$\psi(u)\in 1+u^{-2}\C[[u^{-2}]]$ such that 
$(1+1/2u^{-1})\psi(u)\mu(u)$ is a polynomial
in $u^{-1}$ with a decomposition
$$
(1+\frac12 u^{-1})\psi(u)\mu(u)=(1-\gamma_1u^{-1})\cdots
(1-\gamma_{2k+1}u^{-1}), \tag 5.8
$$
and $\gamma_{2i-1}+\gamma_{2i}\in\ZZ_+$, $i=1,\dots,k$ 
for a certain re-enumeration
of the $\gamma_i$. $\square$
\endproclaim

By analogy with Theorem 4.4 we can reformulate these conditions in terms
of monic polynomials in $u$.

\bigskip
\proclaim
{\bf Theorem 5.4} The irreducible highest weight representation
$V(\mu(u))$ 
of $\Y^+(2)$ is
finite-dimensional if and only if there exist a monic
polynomial $P(u)\in\C[u]$ with
$
P(u)=P(-u+1)
$
and $\gamma\in\C$ such that $P(-\gamma)\ne 0$ and
$$
\frac{\mu(-u)}{\mu(u)}=\frac{P(u+1)}{P(u)}
\frac{(u+\gamma)(2u+1)}{(u-\gamma)(2u-1)}.
\tag 5.9
$$
In this case the pair $(P(u),\gamma)$ is unique.
\endproclaim

\Proof If $\dim V(\mu(u))<\infty$ then by 
Proposition 5.3 relation (5.8) holds for a certain even
series $\psi(u)$. Re-enumerate the numbers $\gamma_i$ to satisfy
condition (5.3). By Proposition 5.2 we have 
$\gamma_{2i-1}+\gamma_{2i}\in\ZZ_+$ for $i=1,\dots,k$ since
$V(\mu(u))$ is finite-dimensional. 
Define $P(u)$
as in the proof of Theorem 4.4 by (4.30)--(4.32) and set
$\gamma=\gamma_{k+1}$.
The relation $P(-\gamma)\ne 0$ is implied by (5.3).

Conversely, let (5.9) be satisfied for
a monic polynomial $P(u)$ such that $P(u)=P(-u+1)$ and $\gamma\in\C$
with $P(-\gamma)\ne 0$.
Let
$\{-\delta_1,\dots,-\delta_s, \delta_1+1,\dots,\delta_s+1\}$
be the roots of $P(u)$; see the proof of Theorem 4.4.
Set
$$
\align
\mu'(u)={}&(1-(\delta_1+1)u^{-1})\cdots (1-(\delta_s+1) u^{-1})\\
{}&(1+\delta_1u^{-1})\cdots (1+\delta_s u^{-1})(1-\gamma u^{-1})
\left(1+\frac12 u^{-1}\right)^{-1}.
\endalign
$$
By Proposition 5.3 the representation $V(\mu'(u))$
is finite-dimensional. Since $\mu'(u)$ satisfies
(5.9) the ratio $\mu(u)/\mu'(u)$ is an even formal series
in $u^{-1}$. This proves that $V(\mu(u))$ is also finite-dimensional.

To prove the uniqueness of the pair $(P(u),\gamma)$ suppose that
$(R(u),\delta)$ is another pair, where $R(u)$ is a monic polynomial
in $u$ such that $R(u)=R(-u+1)$, $R(-\delta)\ne 0$ and
$$
\frac{P(u+1)}{P(u)}
\frac{(u+\gamma)}{(u-\gamma)}=\frac{R(u+1)}{R(u)}
\frac{(u+\delta)}{(u-\delta)}.
\tag 5.10
$$
If $\gamma=\delta$ then we get $P(u)=R(u)$; see the proof of
Theorem 2.8. Let $\gamma\ne \delta$. 
We prove
by induction on the sum
of the degrees of the polynomials $P(u)$ and $R(u)$
that (5.10) is impossible. If $P(u)=R(u)=1$
then this is obvious. Further, suppose $P(u)$ is of degree ${}\geq 2$.
Take a root $u_0$ of $P(u)$ such that
$u_0+1$ is not a root. Rewriting (5.10) in the form
$$
P(u)R(u+1)(u+\delta)(u-\gamma)=P(u+1)R(u)(u-\delta)(u+\gamma)
\tag 5.11
$$
we find that either $R(u_0)=0$, or $u_0=\delta$.

If $R(u_0)=0$ then
we may write
$$
P(u)=P'(u)(u-u_0)(u+u_0-1), \qquad R(u)=R'(u)(u-u_0)(u+u_0-1)
$$
and use the induction hypotheses for the polynomials $P'(u)$ and $R'(u)$.

If $u_0=\delta$ but $R(u_0)\ne 0$ then we write
$$
P(u)=P'(u)(u-u_0)(u+u_0-1)
$$
and obtain from (5.11) that
$$
P'(u)R(u+1)(u+u_0-1)(u-\gamma)=P'(u+1)R(u)(u-u_0+1)(u+\gamma).
$$
Denote $\delta'=u_0-1$. Then we have $\delta'\ne\gamma$ because
$P(-\delta')=0$. Moreover, $R(-\delta')=R(u_0)\ne 0$ and we may
apply the induction hypotheses to the polynomials $P'(u)$ and
$R(u)$. $\square$

\bigskip
\proclaim
{\bf Corollary 5.5} There is a one-to-one correspondence
between finite-dimensional irreducible representations
of the special twisted Yangian $\SY^+(2)$ and the 
pairs $(P(u),\gamma)$, where
$P(u)$ is a monic polynomial 
in $u$ such that $P(u)=P(-u+1)$ and $\gamma\in\C$ with $P(-\gamma)\ne 0$.
Every such representation is isomorphic to a representation of the form 
(5.4). $\square$
\endproclaim

\medskip

The following analogue of Corollary 4.7 holds.

\bigskip
\proclaim
{\bf Corollary 5.6}  Consider an irreducible finite-dimensional
representation of $\Y(2)$ of the form
$$
L(\gamma_1,-\gamma_2)\ot\cdots\ot L(\gamma_{2k-1},-\gamma_{2k}).
\tag 5.12
$$
Then the representation
$$
L(\gamma_1,-\gamma_2)\ot\cdots\ot L(\gamma_{2k-1},-\gamma_{2k})
\ot V(-\gamma)
\tag 5.13
$$
of $\Y^+(2)$ is irreducible
if and only if for each pair $1\leq i< j\leq k$ the strings
$S(\gamma_{2i},-\gamma_{2i-1})$ and $S(\gamma_{2j-1},-\gamma_{2j})$
are in general position, and for each $i=1,\dots,k$ one has
$$
\gamma\not\in S(\gamma_{2i-1},-\gamma_{2i})\qquad\text{and}\qquad
\gamma\not\in S(\gamma_{2i},-\gamma_{2i-1}).
\tag 5.14
$$
\endproclaim

\Proof By Corollary 2.11, the strings $S(\gamma_{2i-1},-\gamma_{2i})$,
$i=1,\dots,k$
are pairwise in general position. Permuting the tensor factors in
(5.12) if necessary, we may assume that
$
\gamma_1+\gamma_2\leq \cdots \leq \gamma_{2k-1}+\gamma_{2k}.
$
If the condition on the strings holds then one easily checks that
(5.3) is satisfied for $\gamma_{2k+1}=\gamma$
and so, by Proposition 5.2
representation (5.13) is irreducible.

To prove the `only if' part
consider first the case $k=1$.
Let the representation 
$L(\gamma_1,-\gamma_2)\ot V(-\gamma)$
be irreducible but condition (5.14) is violated.
Note that $L(\gamma_2,-\gamma_1)\ot V(-\gamma)$ is also
irreducible. Indeed, these two representations have the same 
dimension and
the tensor product of the highest weight vector
of $L(\gamma_2,-\gamma_1)$ and the basis vector of $V(-\gamma)$
generates a highest weight module with the same highest weight
as $L(\gamma_1,-\gamma_2)\ot V(-\gamma)$. So we may assume that
$\gamma\in S(\gamma_{1},-\gamma_{2})$ and if $\gamma+\gamma_1\in\ZZ_+$
then $\gamma+\gamma_2\leq \gamma+\gamma_1$. Then
by Proposition 5.2 the
representation $L(\gamma,-\gamma_2)\ot V(-\gamma_1)$
is irreducible.
It is isomorphic
to $L(\gamma_1,-\gamma_2)\ot V(-\gamma)$ because they have
the same
highest weight. Comparing the dimensions we come to a contradiction.

Consider now the general case. Let representation (5.13)
be irreducible. If (5.14) is violated then 
we move the corresponding tensor factor in (5.13) to the right
to get the tensor product of the form
$L(\gamma_{2i-1},-\gamma_{2i})\ot V(-\gamma)$. 
However, this representation
is reducible. Contradiction.

Finally, replacing some of the tensor factors
$L(\gamma_{2i-1},-\gamma_{2i})$ in (5.13) respectively by 
$L(\gamma_{2i},-\gamma_{2i-1})$ we get an isomorphic
representation of $\Y^+(2)$ (see the proof of Corollary 4.7).
Hence, the representation of $\Y(2)$ obtained from (5.12)
by this replacement is also irreducible. Therefore, by Corollary 2.11
the corresponding strings
are pairwise in general position.
$\square$ \medskip

Note that $V(\gamma)$ for $\gamma=1/2$ is the trivial representation
of $\Y^+(2)$. Thus, in this particular case
Corollary 5.6 provides
a criterion of irreducibility of the restriction of (5.12) to $\Y^+(2)$.

\bigskip
\noindent
{\it Example 5.7.} Consider the restriction of the representation
$L(\gamma_{1},-\gamma_{2})$ of $\Y(2)$ with $\gamma_1+\gamma_2\in\ZZ_+$
to $\Y^+(2)$. By Corollary 5.6,
it is reducible if and only if the string $S(\gamma_{1},-\gamma_{2})$
contains $-1/2$. In this case
$\gamma_1,\gamma_2\in 1/2+\ZZ_+$.
Denote the highest weight vector
of $L(\gamma_{1},-\gamma_{2})$ by $\xi$ and set $\xi_p=(E_{1,-1})^p\xi$.
Suppose first that $\gamma_1\geq \gamma_2$.
Then the following are $\Y^+(2)$-submodules in $L(\gamma_{1},-\gamma_{2})$:
$$
\align
L_1=&\ \ \text{span of}\ \ \{\xi_0,\dots,\xi_{\gamma_2-1/2}\},\\
L_2=&\ \ \text{span of}\ \
\{\xi_{\gamma_1+1/2},\dots,\xi_{\gamma_1+\gamma_2}\}.
\endalign
$$
By Proposition 5.2 one has $L_1\simeq L(\gamma_2,-1/2)\ot V(-\gamma_1)$
and $L_2\simeq L(\gamma_2,-1/2)\ot V(\gamma_1+1)$. The quotient
$L(\gamma_{1},-\gamma_{2})/L$, where $L=L_1\oplus L_2$, is isomorphic to
$L(\gamma_1,\gamma_2+1)$.

If $\gamma_1\leq \gamma_2$ then $L(\gamma_{1},-\gamma_{2})$ contains
the submodule
$$
L=\ \ \text{span of}\ \ \{\xi_{\gamma_1+1/2},\dots,\xi_{\gamma_2-1/2}\}
$$
which is isomorphic to $L(\gamma_2,\gamma_1+1)$ and the quotient
$L(\gamma_{1},-\gamma_{2})/L$ is isomorphic to the direct sum
of $L(\gamma_1,-1/2)\ot V(-\gamma_2)$ and 
$L(\gamma_1,-1/2)\ot V(\gamma_2+1)$. $\square$
\medskip

Note that the mapping
$$
s_{ij}(u)\mapsto s_{i'j'}(u), \tag 5.15
$$
where $i'=-i$ for $i=-1,1$ and $i'=i$ otherwise,
defines an automorphism of the algebra $\Y^+(2n)$. 
Consider the case $n=1$ and denote
the composition of this automorphism with
a finite-dimensional irreducible representation $V(\mu(u))$ of $\Y^+(2)$
by $V(\mu(u))^{\sharp}$. By Theorem 3.3,
this is a representation with the highest weight
which we denote
by $\mu(u)^{\sharp}$. We shall say
that $\mu(u)^{\sharp}$ is 
{\it well-defined\/} if
$\mu(u)$ satisfies the conditions of Proposition 5.3. 
Since the automorphisms (1.29) and (5.15) commute,
the ratio $\mu(u)/\mu(u)^{\sharp}$ takes the same value for all
representations from the similarity class containing $V(\mu(u))$.

Equivalently, $V(\mu(u))$ contains a unique (up to scalar multiples)
nonzero vector $\zeta$ such that $s_{1,-1}(u)\zeta=0$ (cf. Remark 3.4).
Then $\zeta$ is an eigenvector for $s_{-1,-1}(u)$ and $\mu(u)^{\sharp}$
is defined by the formula
$$
s_{-1,-1}(u)\ts \zeta=\mu(u)^{\sharp}\ts\zeta.
\tag 5.16
$$

\bigskip
\proclaim
{\bf Proposition 5.8} Let $\mu(u)$ be a series given by
(5.2) such that condition (5.3) is satisfied and
$\gamma_{2i-1}+\gamma_{2i}\in\ZZ_+$ for $i=1,\dots,k$.
Then
$$
\mu(u)^{\sharp}=
(1-\gamma_1u^{-1})\cdots (1-\gamma_{2k}u^{-1})
(1+(\gamma_{2k+1}+1)u^{-1})\left(1+\frac12 u^{-1}\right)^{-1}.
\tag 5.17
$$
\endproclaim

\Proof By Proposition 5.2 the representation $V(\mu(u))$
is isomorphic to the tensor product (5.4). In the notation
of the proof of Proposition 5.2 set
$$
\zeta=(E_{1,-1})^{\gamma_1+\gamma_2}\xi_1\ot\cdots\ot
(E_{1,-1})^{\gamma_{2k-1}+\gamma_{2k}}\xi_k\ot\xi_{k+1}.
$$
Then using (1.31) we verify that $\zeta$ satisfies the relations
$$
s_{-1,-1}(u)\zeta=(1-\gamma_1u^{-1})\cdots (1-\gamma_{2k}u^{-1})
(1+(\gamma_{2k+1}+1)u^{-1})\left(1+\frac12 u^{-1}\right)^{-1}\zeta
$$
and $s_{1,-1}(u)\zeta=0$. Thus, $V(\mu(u))^{\sharp}$ is isomorphic to
$V(\mu(u)^{\sharp})$ with $\mu(u)^{\sharp}$ given by (5.17). $\square$
\medskip

The following theorem together with Theorem 3.3
gives a description
of finite-dimensional irreducible representations of $\Y^+(2n)$
(cf. Theorem 4.8). We use notation (2.35) and (4.34).

\bigskip
\proclaim
{\bf Theorem 5.9} The irreducible highest weight representation
$V(\mu(u))$, $\mu(u)=(\mu_1(u),\dots,\mu_n(u))$ of 
$\Y^+(2n)$, $n\geq 2$ is finite-dimensional if and only 
if $\mu_1(u)^{\sharp}$ is well-defined and either one of the following
four relations holds:
$$
\align
\mu_1(-u)\Ra \mu_1&(u)\to \mu_2(u)\to\cdots\to\mu_n(u),\tag 5.18\\
\frac{2u-1}{2u+1}\ts\mu_1(-u)\Ra \mu_1&(u)\to \mu_2(u)\to\cdots\to\mu_n(u),
\tag 5.19\\
\mu_1(-u)^{\sharp}\Ra \mu_1&(u)^{\sharp}\to \mu_2(u)\to\cdots\to\mu_n(u),
\tag 5.20\\
\frac{2u-1}{2u+1}\ts\mu_1(-u)^{\sharp}\Ra \mu_1&(u)^{\sharp}\to 
\mu_2(u)\to\cdots\to\mu_n(u).\tag 5.21
\endalign
$$
\endproclaim

\Proof Suppose that $\dim V(\mu(u))<\infty$.
Then by Proposition 3.5 we have
$$
\mu_1(u)\to \mu_2(u)\to\cdots\to\mu_n(u).\tag 5.22
$$
Denote the composition of $V(\mu(u))$ with the 
automorphism (5.15)
by $V(\mu(u))^{\sharp}$.
By Theorem 3.3 the representation $V(\mu(u))^{\sharp}$
is highest weight. To find its highest weight vector consider
the cyclic span $V=\Y^+(2)\ts\xi$, where
$\xi$ is the
highest weight vector of $V(\mu(u))$, and $\Y^+(2)$
is identified with the subalgebra in $\Y^+(2n)$ 
generated by the coefficients of the series 
$s_{ij}(u)$ with $i,j=-1,1$.
We claim that $V$ is
an irreducible $\Y^+(2)$-module. Indeed, let
$J$ be the left ideal in $\Y^+(2n)$ generated by the
$s_{kl}^{(r)}$ with $k<l$ and $(k,l)\ne (-1,1)$.
Then $V$ is annihilated by $J$. This follows from the fact that
$J\ts\xi=0$ and the relations
$$
[s_{ij}(u),s_{kl}(v)]\equiv 0 \mod J,
\tag 5.23
$$
where $i,j=\pm 1$ and $k<l$, $(k,l)\ne (-1,1)$. 
To verify (5.23) we may assume, applying (1.18) if necessary,
that $l\geq 2$.
However, in this case (5.23) is
immediate from (1.17).
Now, if $\wt V$ is a nonzero submodule in $V$
then there exists a nonzero vector $\eta\in \wt V$
such that  $s_{-1,1}(u)\eta=0$. Then $s_{kl}(u)\eta=0$ for all $k<l$
since $V$ is annihilated by $J$. So,
$\eta$ is a multiple of $\xi$; see Remark 3.4.
This proves that $\wt V=V$. 

Thus, $V$ is isomorphic to 
the representation $V(\mu_1(u))$ of $\Y^+(2)$. Since 
$\dim V<\infty$ there exists
a nonzero vector $\zeta\in V$ such that
$$
s_{-1,-1}(u)\ts\zeta=\mu_1(u)^{\sharp}\ts\zeta \qquad\text{and}\qquad
s_{1,-1}(u)\zeta=0;
$$
see (5.16).
Moreover, we have $s_{kk}(u)\zeta=\mu_k(u)\zeta$ for $k=2,\dots, n$.
This is implied by the following
relations which are immediate from (1.17):
$$
[s_{ij}(u),s_{kk}(v)]\equiv 0 \mod J,
$$
where $i,j=\pm 1$ and $k\geq 2$. Since $\zeta$ is annihilated by $J$
we may conclude that $\zeta$ is the highest weight vector
of $V(\mu(u))^{\sharp}$ and so, the latter is a representation with
the highest weight $(\mu_1(u)^{\sharp},\mu_2(u),\dots,\mu_n(u))$.
Applying Proposition 3.5 again, we find that
$$
\mu_1(u)^{\sharp}\to \mu_2(u)\to\cdots\to\mu_n(u).\tag 5.24
$$

Further, since the representation $V(\mu_1(u))$ of $\Y^+(2)$
is finite-dimensional, then by Proposition 5.3 we may assume, applying an
automorphism of the form (1.29) if necessary, that
$(1+1/2 u^{-1})\mu_1(u)$ is a polynomial in $u^{-1}$. By (5.22)
we have
$$
\frac{\mu_1(u)}{\mu_2(u)}=\frac{P(u+1)}{P(u)}\tag 5.25
$$
for a monic polynomial $P(u)$ in $u$. Hence,
applying an
automorphism (1.29) again, 
we may assume that both $(1+1/2 u^{-1})\mu_1(u)$
and $(1+1/2 u^{-1})\mu_2(u)$ are polynomials in $u^{-1}$:
$$
\align
\mu_1(u)={}&(1-\gamma_1u^{-1})\cdots (1-\gamma_{2k+1}u^{-1})
\left(1+\frac12 u^{-1}\right)^{-1},\\
\mu_2(u)={}&(1-\delta_1u^{-1})\cdots (1-\delta_{2k+1}u^{-1})
\left(1+\frac12 u^{-1}\right)^{-1}.
\endalign
$$
Re-enumerating the $\gamma_i$
if necessary, we may assume that condition (5.3) holds. Moreover,
$\gamma_{2i-1}+\gamma_{2i}\in\ZZ_+$ for $i=1,\dots,k$ since
$\dim V(\mu_1(u))<\infty$; see Proposition 5.2.
Further, (5.25) implies that there exists a
re-enumeration of the $\delta_i$ such that
$$
\delta_i-\gamma_i\in\ZZ_+ \qquad\text{for}\quad i=1,\dots,2k+1.
\tag 5.26
$$
To see this, choose a string-type subset
$\{a+1,a+2,\dots,b\}$ of roots of $P(u)$ such that $a$ and $b+1$
are not roots. Then we find from (5.25) that $\gamma_i=a$ 
for a certain $i$, and $b$ equals one of the $\delta$'s which we
re-denote by $\delta_i$.
Therefore, 
$\delta_i-\gamma_i\in\ZZ_+$ and we can write
$$
\mu_1(u)=(1-au^{-1})\mu'_1(u),\qquad \mu_2(u)=(1-bu^{-1})\mu'_2(u)
$$
and
$$
\frac{P(u+1)}{P(u)}=\frac{P'(u+1)}{P'(u)}\cdot\frac{u-a}{u-b},
$$
where $P(u)=(u-a-1)\cdots (u-b)P'(u)$. Then we proceed by induction
applying the same argument to the series $\mu'_1(u)$, $\mu'_2(u)$
and the polynomial $P'(u)$.

Now we use Proposition 5.8. Exactly as above, we obtain from (5.24)
that there exists a permutation $\pi\in\Sym_{2k+1}$ such that
$$
\delta_{\pi(i)}-\gamma_i\in\ZZ_+,\ \  i=1,\dots,2k\quad \text{and}\quad
\delta_{\pi(2k+1)}+\gamma_{2k+1}+1\in\ZZ_+.
\tag 5.27
$$
Now, (5.26) and (5.27) imply that $\delta_{\pi(i)}-\delta_i\in\ZZ$
for $i=1,\dots,2k$. But then this has to be true for $i=2k+1$ as well.
Thus, $2\gamma_{2k+1}+1\in\ZZ$, that is, the number 
$\gamma:=\gamma_{2k+1}$ belongs to $\frac12 \ZZ$.

Suppose first that $\gamma\in\frac12+\ZZ$ and $\gamma\geq -\frac12$.
Then exactly as in the proof of Theorem 4.4 we obtain that
$\mu_1(-u)\Ra\mu_1(u)$. Similarly, if 
$\gamma\in\frac12+\ZZ$ and $\gamma\leq -\frac12$ then
$\mu_1(-u)^{\sharp}\Ra\mu_1(u)^{\sharp}$. Further, let $\gamma\in\ZZ$
and $\gamma\geq 0$. Then writing
$$
(1+\frac12 u^{-1})\mu_1(u)=(1-\gamma_1u^{-1})\cdots (1-\gamma_{2k}u^{-1})
(1-\gamma u^{-1})(1-0u^{-1})
$$
we see that $(1-1/2 u^{-1})\mu_1(-u)\Ra (1+1/2 u^{-1})\mu_1(u)$,
or, equivalently,
$$
\frac{2u-1}{2u+1}\ts\mu_1(-u)\Ra \mu_1(u).
$$
Similarly, if $\gamma\in\ZZ$
and $\gamma\leq -1$ then using Proposition 5.8 we obtain
$$
\frac{2u-1}{2u+1}\ts\mu_1(-u)^{\sharp}\Ra \mu_1(u)^{\sharp}.
$$
Combining these relations with (5.22) or (5.24) we complete the 
proof of the `only if'
part.
\medskip

Suppose now that relation (5.18) is satisfied. Then we repeat the
corresponding
argument of the proof of Theorem 4.8 to show that
$V(\mu(u))$ is finite-dimensional. Similarly, if (5.20) holds
then the irreducible
representation with the highest weight
$(\mu_1(u)^{\sharp},\mu_2(u),\dots,\mu_n(u))$ is finite-dimensional,
and so $V(\mu(u))$ is, since the automorphism (5.15) is involutive.

If (5.19) holds then
$$
\frac{\mu_1(-u)}{\mu_1(u)}=\frac{P_1(u+1)}{P_1(u)}\cdot\frac{2u+1}{2u-1}
\tag 5.28
$$
for a monic polynomial $P_1(u)$ of the form (4.41), and 
relations (4.39) hold
for polynomials $P_i(u)$ of the form (4.40).
Define the $\lambda_j(u)$ by (4.42) and (4.43). Then the
representation $L(\lambda(u))$ of $\Y(2n)$ is
finite-dimensional by Theorem 2.12. 

Let $V(\mu_0)$ denote the irreducible representation of the Lie
algebra $\oa(2n)$ with the highest weight $\mu_0=(-1/2,\dots,-1/2)$. It is
finite-dimensional; see (3.6) and (3.7). Extend $V(\mu_0)$ to
a representation of $\Y^+(2n)$ by using (1.23)
and consider the tensor product
$
L(\lambda(u))\ot V(\mu_0).
$ 
The tensor product of the
highest weight vectors of $L(\lambda(u))$ and $V(\mu_0)$ generates
a $\Y^+(2n)$-submodule with the highest weight
$\mu'(u)=(\mu'_1(u),\dots,\mu'_n(u))$ where
$\mu'_i(u)=\lambda_i(u)\lambda_{-i}(-u)(1+1/2 u^{-1})^{-1}$;
see (3.10) and (3.11). So the
representation $V(\mu'(u))$ of $\Y^+(2n)$ is finite-dimensional.
However, the polynomials $\mu'_i(u)$
satisfy (5.28) and (4.39). Hence, there exists an
automorphism of $\Y^+(2n)$ of the form (1.29)
such that its composition with the representation
$V(\mu'(u))$ is isomorphic to
$V(\mu(u))$ which completes the proof in this case.

Finally, if relation (5.21) holds then the same argument shows that 
the representation
$V(\mu_1(u)^{\sharp},\mu_2(u),\dots,\mu_n(u))$ is finite-dimensional,
and so $V(\mu(u))$ is. $\square$
\medskip

We conclude this section with the
following description of representations of the
special twisted Yangian $\SY^+(2n)$ implied by Theorem 5.9
(cf. Corollary 4.9).

\bigskip
\proclaim
{\bf Corollary 5.10} There is
a one-to-one correspondence between
finite-dimensional irreducible representations
of the special twisted Yangian $\SY^+(2n)$, $n\geq 2$ and the families
$\{P_1(u),\dots,P_{n}(u),\varepsilon\}$, where the $P_i(u)$ are monic
polynomials in $u$ with $P_1(u)=P_1(-u+1)$, and the parameter $\varepsilon$
takes values in $\{1,2,3,4\}$; the values $\varepsilon=1$ and
$\varepsilon=3$
are identified if $P_1(1/2)\ne 0$.
Every such representation is
isomorphic to a subquotient of a representation 
of the form $L(\lambda(u))$ or $L(\lambda(u))\ot V(\mu_0)$.
\endproclaim

\Proof Each of the four conditions (5.18)--(5.21) allows one
to naturally associate a family $\{P_1(u),\dots,P_{n}(u)\}$
of monic polynomials to the corresponding similarity class
of finite-dimensional irreducible representations of $\Y^+(2n)$;
see the proof of Theorem 5.9.
We associate the values of $\varepsilon$ respectively to
conditions (5.18)--(5.21). It follows from the proof
of Theorem 5.9 that two families of monic polynomials
with values of $\varepsilon=1$ and $\varepsilon=3$
correspond to the 
same similarity class if
and only if $P_1(1/2)\ne 0$. One easily checks that identifying
these values we obtain the desired one-to-one correspondence.
The second claim also follows from the proof of Theorem 5.9.
$\square$

\bigskip
\bigskip
\noindent
{\bf 6. Representations of $\Y^+(2n+1)$}
\bigskip

The description of finite-dimensional irreducible representations of
$\Y^+(2n+1)$ is based on the simplest nontrivial case $n=1$.
Our main instrument in this case is an investigation of
the restriction of an irreducible highest weight representation
of $\Y^+(3)$ to the subalgebra $\Y^+(2)$ and the use of the results 
of the previous section.

\medskip

Irreducible highest weight representations $V(\mu(u))$ of $\Y^+(3)$
are parameterized by pairs of formal series 
$\mu(u)=(\mu_0(u),\mu_1(u))$, where
$\mu_0(u)\in 1+u^{-2}\C[[u^{-2}]]$ and
$\mu_1(u)\in 1+u^{-1}\C[[u^{-1}]]$.

\bigskip
\proclaim
{\bf Proposition 6.1} If $\ts\dim V(\mu(u)) <\infty$ then
there exists a formal series 
$\varphi(u)\in 1+u^{-1}\C[[u^{-1}]]$ such that
$\varphi(u)\mu_0(u)$ and 
$\varphi(u)\mu_1(u)$ are polynomials in $u^{-1}$.
\endproclaim

\Proof Denote by $\xi$ the highest weight vector of $V(\mu(u))$.
Let $k$ denote the minimum nonnegative integer such that
$s_{10}^{(k+1)}\xi$ is a linear combination of the vectors
$\xi_1:=s_{10}^{(1)}\xi$, $\dots$, $\xi_k:=s_{10}^{(k)}\xi$.
Relations (1.17) give
$$
\aligned
[s_{11}&(u),s_{10}(v)]=\frac{1}{u-v}
\left(s_{10}(v)s_{11}(u)-s_{10}(u)s_{11}(v)\right)\\
{}+{}&\frac{u-v+1}{u^2-v^2}s_{1,-1}(v)s_{01}(u)
-\frac{1}{u+v}s_{1,-1}(u)s_{-1,0}(v)
-\frac{1}{u^2-v^2}s_{1,-1}(u)s_{01}(v).
\endaligned
\tag 6.1
$$
Since the vector $\xi$ is annihilated by $s_{01}(u)$
(and hence by $s_{-1,0}(u)$) this makes the consideration here
very similar to the case of $\Y(2)$; see Proposition 2.5.
In particular, repeating the corresponding argument of the
proof of Proposition 2.5 we can show that
$$
s_{10}(u)\xi=a_1(u)\xi_1+\cdots+a_k(u)\xi_k
$$
for certain series $a_i(u)$; cf. (2.13).
By (1.17) we have
$$
\aligned
[s_{00}(u),s_{10}(v)]={}&\frac{u+v+1}{u^2-v^2}s_{10}(u)s_{00}(v)\\
{}-{}&\frac{1}{u+v}s_{0,-1}(u)s_{00}(v)
-\frac{2v+1}{u^2-v^2}s_{10}(v)s_{00}(u).
\endaligned
\tag 6.2
$$
Taking the coefficients at $u^{-2}$ in (6.1) and (6.2) 
and using again the argument of the proof of Proposition 2.5 we
find that the components of $\mu(u)$ can be written in the form
$$
\mu_0(u)=B(u)a_k(u)\qquad\text{and}\qquad \mu_1(u)=C(u)a_k(u),
$$
where $B(u)$ and $C(u)$ are monic polynomials in $u$ 
of degree $k$. Thus, the series $\varphi(u)$ can be defined by
$\varphi(u)=(a_k(u)u^k)^{-1}$.
$\square$
\medskip

Write the polynomials $\varphi(u)\mu_0(u)$ and 
$\varphi(u)\mu_1(u)$ in the form
$$
\aligned
\varphi(u)\mu_0(u)={}&(1+\alpha_1 u^{-1})\cdots (1+\alpha_k u^{-1}),\\
\varphi(u)\mu_1(u)={}&(1+\beta_1 u^{-1})\cdots (1+\beta_k u^{-1}),
\endaligned
\tag 6.3
$$
with $\alpha_i,\beta_i\in\C$. Since $\mu_0(u)$ is an even series in
$u^{-1}$ we can find an appropriate
automorphism of the form (1.29) such that its
composition with $V(\mu(u))$ is an irreducible representation with
the following highest weight which we shall again denote by $\mu(u)$: 
$$
\aligned
\mu_0(u)={}&(1-\alpha_1^2 u^{-2})\cdots (1-\alpha_k^2 u^{-2}),\\
\mu_1(u)={}&(1-\alpha_1 u^{-1})\cdots (1-\alpha_k u^{-1})
(1+\beta_1 u^{-1})\cdots (1+\beta_k u^{-1}).
\endaligned
\tag 6.4
$$
Our next task is to find necessary conditions for the $\alpha_i$ and
$\beta_i$ so that the representation $V(\mu(u))$ with the highest weight
given by (6.4) be finite-dimensional.

Consider the following representation of $\Y(3)$:
$$
L(\alpha_1,\alpha_1,\beta_1)\ot\cdots \ot L(\alpha_k,\alpha_k,\beta_k).
\tag 6.5
$$
It follows from
(1.15) that the $\Y^+(3)$-cyclic span of
the tensor product of the highest weight vectors of
the $L(\alpha_i,\alpha_i,\beta_i)$
is a representation of
$\Y^+(3)$ with the highest weight $\mu(u)=(\mu_0(u),\mu_1(u))$
given by (6.4).
The generators $t_{ij}^{(r)}$ of $\Y(3)$
with $r\geq k+1$ act trivially in representation (6.5)
(cf. (2.18)).
Therefore, by (1.15)
the generators $s_{ij}^{(r)}$ of $\Y^+(3)$ with $r\geq 2k+1$ act
trivially in $V(\mu(u))$. The series
$$
\sss_{ij}(u):=u^{2k}\ts s_{ij}(u)
\tag 6.6
$$
is then a polynomial operator in $u$ and we may evaluate $u$
to get well-defined operators in $V(\mu(u))$. Set
$$
\aligned
a_0(u)=u^{2k}\ts\mu_0(u)={}&(u^2-\alpha_1^2)\cdots (u^2-\alpha_k^2),\\ 
a_1(u)=u^{2k}\ts\mu_1(u)={}&(u-\alpha_1)\cdots (u-\alpha_k)
(u+\beta_1)\cdots (u+\beta_k).
\endaligned
\tag 6.7
$$

Using Proposition 1.1 we shall now equip 
the vector space $V:=V(\mu(u))$ with another
structure of $\Y^+(3)$-module. The matrix elements of the Sklyanin
comatrix $\wh S(u)$ can be found from [M3, Section 6]. 
In particular,
$$
\aligned
\wh s_{00}(u)={}&s_{11}(-u)s_{11}(u-1)-s_{1,-1}(-u)s_{-1,1}(u-1),\\
\wh s_{-1,0}(u)={}&s_{0,-1}(-u)s_{-1,1}(u-1)-s_{01}(-u)s_{11}(u-1),\\ 
\wh s_{11}(u)={}&s_{11}(-u)s_{00}(u-1)-
s_{10}(-u)s_{-1,0}(u-1).
\endaligned
\tag 6.8
$$
Let us take the composition of the representation $V(\mu(u))$
with the automorphism $S(u)\to \wh S(-u+1/2)$ of $\Y^+(3)$;
see (1.27). One can easily see from (6.8) that
the resulting representation
is highest weight and generated by
the vector $\xi$. The highest
weight of this representation is
$$
(\mu_1(u-1/2)\mu_1(-u-1/2),\ \mu_1(u-1/2)\mu_0(-u-1/2)).
\tag 6.9
$$
We shall use the notation
$$
\alpha^*_i=\frac12 -\beta_i,\quad \beta^*_i=\frac12 -\alpha_i,
\qquad\text{for}\quad i=1,\dots,k.
\tag 6.10
$$
Set
$$
\aligned
\mu^*_0(u)={}&(1-(\alpha^*_1)^2 u^{-2})\cdots (1-(\alpha^*_k)^2 u^{-2}),\\
\mu^*_1(u)={}&(1-\alpha^*_1 u^{-1})\cdots (1-\alpha^*_k u^{-1})
(1+\beta^*_1 u^{-1})\cdots (1+\beta^*_k u^{-1}).
\endaligned
\tag 6.11
$$
Take the composition of the irreducible representation of $\Y^+(3)$
with the highest weight (6.9) and the automorphism (1.29)
corresponding to the series
$$
\psi_0(u)=\frac{\mu^*_0(u)}{\mu_1(u-1/2)\mu_1(-u-1/2)}.
$$
The resulting representation is isomorphic to $V(\mu^*(u))$,
$\mu^*(u)=(\mu^*_0(u),\mu^*_1(u))$.

Thus, the vector space $V$ carries two representations of $\Y^+(3)$
isomorphic to $V(\mu(u))$ and $V(\mu^*(u))$, respectively.
To distinguish them we shall denote by $s^*_{ij}(u)$ the action
of the generators $s_{ij}(u)$ in $V(\mu^*(u))$. By the above definition
we have
$$
s^*_{ij}(u)=\psi_0(u)\ts \wh s_{ij}(-u+\frac12). \tag 6.12
$$
Denote by $\sst_{ij}(u)$
the corresponding polynomial operator (6.6) in $V(\mu^*(u))$. 
By analogy with (6.7)
set
$$
a^*_0(u)=u^{2k}\ts\mu^*_0(u),\qquad
a^*_1(u)=u^{2k}\ts\mu^*_1(u).
$$

For the proof of the next proposition
fix an index $i\in\{1,\dots,k\}$ and set $\alpha:=\alpha_i$, 
$\beta^*:=\frac12-\alpha_i$.
Given nonnegative integer $p$ introduce 
the following vector in $V$:
$$
\eta_{p}=
\sss_{10}(-\alpha+p-1)\cdots \sss_{10}(-\alpha) \ts \xi.
$$

\bigskip
\proclaim
{\bf Proposition 6.2} We have the relations
$$
\align
\sss_{00}(u) \eta_{p}={}&a_0(u)\frac{u^2-(\alpha-p)^2}
{u^2-\alpha^2}\ts \eta_p,
\tag 6.13\\
\sss_{-1,0}(u) \eta_{p}={}&p\ts a_0(u) a_1(-\alpha+p-1)
\frac{u+\alpha-p}
{u^2-\alpha^2}\ts \eta_{p-1}, 
\tag 6.14\\
\sst_{-1,1}(u) \eta_{p}={}&0, 
\tag 6.15\\
\sst_{11}(u) \eta_{p}={}&a^*_1(u)\frac{u+\beta^*+p}
{u+\beta^*}\ts \eta_p. 
\tag 6.16
\endalign
$$
\endproclaim

\Proof We use induction on $p$ to prove (6.13)--(6.16)
simultaneously. 
For $p=0$ all the relations are obvious.
Suppose now that $p\geq 1$.
Using the induction hypotheses and (6.2) we obtain
$$
\align
\sss_{00}(u)\ts\eta_p={}&\sss_{00}(u)\sss_{10}(-\alpha+p-1)\ts\eta_{p-1}\\
{}={}&\frac{u^2-(\alpha-p)^2}
{u^2-(\alpha-p+1)^2}\sss_{10}(-\alpha+p-1)
\sss_{00}(u)\ts\eta_{p-1},
\endalign
$$
which proves (6.13). Further, (1.28) gives
$$
[\wh s_{11}(u),s_{10}(v)]=-\frac {1}{u-v-2}
\left(\wh s_{1,-1}(u)s_{-1,0}(v)+ 
\wh s_{10}(u)s_{00}(v)+\wh s_{11}(u)s_{10}(v)\right).
\tag 6.17
$$
Hence, using the induction
hypotheses and (6.12) we can write
$$
\align
\sst_{11}(u) \ts\eta_p={}&\sst_{11}(u)
\sss_{10}(-\alpha+p-1)\ts\eta_{p-1}\\
{}={}&\frac{u+\beta^*+p}
{u+\beta^*+p-1}\sss_{10}(-\alpha+p-1)
\sst_{11}(u)\ts\eta_{p-1},
\endalign
$$
which proves (6.16). The proof of (6.15) is the same; it suffices
to use the expansion for $[\wh s_{-1,1}(u),s_{10}(v)]$
implied by (1.28). Finally, to prove (6.14) consider the following
relation implied by (1.17):
$$
\aligned
[s_{-1,0}&(u),s_{10}(v)]=\frac{1}{u-v}
\left(s_{10}(u)s_{-1,0}(v)-s_{10}(v)s_{-1,0}(u)\right)\\
+&\frac{u-v-1}{u^2-v^2}s_{11}(v)s_{00}(u)
-\frac{1}{u+v}s_{-1,-1}(u)s_{00}(v)
+\frac{1}{u^2-v^2}s_{11}(u)s_{00}(v).
\endaligned
\tag 6.18
$$
Using the induction
hypotheses we obtain
$$
\aligned
\sss_{-1,0}(u) \ts\eta_p={}&\sss_{-1,0}(u) 
\sss_{10}(-\alpha+p-1)\ts\eta_{p-1}\\
{}={}&\frac{u+\alpha-p}
{u+\alpha-p+1}\sss_{10}(-\alpha+p-1)
\sss_{-1,0}(u)\ts\eta_{p-1}\\
{}+{}&\frac{u+\alpha-p}{u^2-(\alpha-p+1)^2}
\sss_{11}(-\alpha+p-1)
\sss_{00}(u)\ts\eta_{p-1}. 
\endaligned
\tag 6.19
$$
Using (6.13) and applying further (6.19) to 
$\sss_{-1,0}(u) \ts\eta_{p-1}$ etc. we find that
$$
\sss_{-1,0}(u)\eta_p=a_0(u)\frac{u+\alpha-p}
{u^2-\alpha^2}\ts \eta^{(1)}_p,
\tag 6.20
$$
where
$$
\eta^{(1)}_p=\sum_{q=1}^p \sss_{10}(-\alpha+p-1)\cdots 
\sss_{11}(-\alpha+q-1)\cdots \sss_{10}(-\alpha)\ts\xi
$$
and $\sss_{11}(-\alpha+q-1)$ takes the $q$th position from the right.

On the other hand, by (6.8)
$$
\wh s_{11}(u+1)=s_{11}(-u-1)s_{00}(u)-
s_{10}(-u-1)s_{-1,0}(u).
$$
Let us multiply both sides of this relation by $(u(u+1))^{2k}$ and
apply them to the vector $\eta_{p-1}$.
By (6.12) we have 
$$
\wh s_{11}(u+1)=s^*_{11}(-u-\frac 12)\ts \psi_0(u+\frac 12)^{-1}.
$$
So, we can use (6.16) to calculate 
$(u(u+1))^{2k}\ts\wh s_{11}(u+1)\eta_{p-1}$. 
Performing this calculation and
using (6.13) and (6.20)
we come to the relation
$$
a_1(-u-1)\ts \eta_{p-1}=\frac {u-\alpha+p-1}{u-\alpha}\ts
\sss_{11}(-u-1)\ts\eta_{p-1}-
\frac {1}{u-\alpha}\ts
\sss_{10}(-u-1)\ts\eta^{(1)}_{p-1}.
\tag 6.21
$$
Note that 
$$
\eta^{(1)}_{p}=\sss_{11}(-\alpha+p-1)\eta_{p-1}+
\sss_{10}(-\alpha+p-1)\eta^{(1)}_{p-1}.
$$
Therefore,
putting $u=\alpha-p$ into (6.21) we get
$$
\eta^{(1)}_{p}=p\ts a_1(-\alpha+p-1)\ts \eta_{p-1}.
\tag 6.22
$$
Together with (6.20) this completes the proof of (6.14). $\square$
\medskip

Given nonnegative integers $p_1,\dots,p_k$ introduce 
vectors in $V$ of a more general form:
$$
\eta_{p_1,\dots,p_k}=\left(\prod_{i=1,\dots,k}^{\rightarrow}
\sss_{10}(-\alpha_i+p_i-1)\cdots \sss_{10}(-\alpha_i) \right) \xi.
\tag 6.23
$$
Set 
$$
S^{\natural}_{-1,0}(u)=\frac{\sss_{-1,0}(u)}{(u+\alpha_1-p_1)\cdots
(u+\alpha_k-p_k)}.
$$

\bigskip
\proclaim
{\bf Proposition 6.3} {\rm (i)} The series $S^{\natural}_{-1,0}(u)
\eta_{p_1,\dots,p_k}$ is a polynomial in $u$ with values in $V$ and for
every
$i=1,\dots,k$
$$
S^{\natural}_{-1,0}(\alpha_i-p_i) \eta_{p_1,\dots,p_k}
=-a_1(-\alpha_i+p_i-1)\prod_{j=1}^k (\alpha_i-\alpha_j-p_i)
\ts\eta_{p_1,\dots,p_i-1,\dots,p_k}.
\tag 6.24
$$
{\rm (ii)} One has the relation
$$
\sss_{00}(u) \ts\eta_{p_1,\dots,p_k}=(u^2-(\alpha_1-p_1)^2)\cdots
(u^2-(\alpha_k-p_k)^2) \eta_{p_1,\dots,p_k}.
\tag 6.25
$$
\endproclaim

\Proof First we verify (ii) by using (6.2) and repeating the corresponding
argument of the proof of Proposition 6.2. Next, by an obvious induction
with the use of (6.18) and (6.25) (cf. the proof of (6.20))
the show that
$$
\aligned
{}&\sss_{-1,0}(u)\ts \eta_{p_1,\dots,p_k}=(u+\alpha_1-p_1)\cdots
(u+\alpha_k-p_k)\\
{}\times{}&\sum_{j=1}^k (u-\alpha_1)\cdots (u-\alpha_{j-1})
(u-\alpha_{j+1}+p_{j+1})\cdots (u-\alpha_{k}+p_{k})\ts
\eta^{(j)}_{p_1,\dots,p_k}, 
\endaligned
\tag 6.26
$$
where $\eta^{(j)}_{p_1,\dots,p_k}$ is the vector in $V$ obtained from
$\eta_{p_1,\dots,p_k}$ by replacing the product
$$
\sss_{10}(-\alpha_j+p_j-1)\cdots \sss_{10}(-\alpha_j)
\tag 6.27
$$ 
in (6.23)
with the sum
$$
\sum_{q=1}^{p_j} \sss_{10}(-\alpha_j+p_j-1)\cdots 
\sss_{11}(-\alpha_j+q-1)\cdots \sss_{10}(-\alpha_j);
$$
here $\sss_{11}(-\alpha_j+q-1)$ takes the $q$th position from the right.
This proves that the series $S^{\natural}_{-1,0}(u)
\eta_{p_1,\dots,p_k}$ is a polynomial in $u$. Moreover, (6.26) and (6.22)
imply (6.24) for $i=k$. To complete the proof we show that
$\eta_{p_1,\dots,p_k}$ is not changed if we arbitrary permute
the products (6.27) in the definition (6.23). We only need to show this
for permutations of adjacent products. By (6.25) it suffices to verify
the following auxiliary claim: if a vector $\eta\in V$ satisfies
$$
\sss_{00}(u)\ts\eta=(u^2-\gamma_1^2)\cdots (u^2-\gamma_k^2)\ts\eta
$$
then
$$
\sss_{10}(\gamma_1)\sss_{10}(\gamma_2)\ts\eta=
\sss_{10}(\gamma_2)\sss_{10}(\gamma_1)\ts\eta. \tag 6.28
$$
Indeed, relations (1.17) give
$$
[s_{10}(u),s_{10}(v)]=
-\frac{1}{u+v}\left(s_{1,-1}(u)s_{00}(v)-s_{1,-1}(v)s_{00}(u)\right).
\tag 6.29
$$
Therefore, if $\gamma_1+\gamma_2\ne 0$ then (6.28) holds because
$\sss_{00}(\gamma_i)\eta=0$. Further, we obtain from (6.29) that
$$
[\sss_{10}(-u),\sss_{10}(u)]=
\frac{d\ts\sss_{1,-1}(u)}{du}\sss_{00}(u)-\sss_{1,-1}(u)
\frac{d\ts\sss_{00}(u)}{du}.
$$
So, if $\gamma_1+\gamma_2= 0$ then (6.28) holds since
in this case $\eta$ is annihilated by both the values of
$\sss_{00}(u)$ and its derivative at $u=\gamma_2$.
$\square$
\medskip

We suppose, as before, that the $p_i$ are nonnegative integers.

\bigskip
\proclaim
{\bf Proposition 6.4} The vector $\eta_{p_1,\dots,p_k}\in V$ 
satisfies the relations
$$
\align
\sst_{-1,1}(u)\eta_{p_1,\dots,p_k}={}&0, \tag 6.30\\
\sst_{11}(u)\eta_{p_1,\dots,p_k}={}&(u-\alpha^*_1)\cdots (u-\alpha^*_k)
(u+\beta^*_1+p_1)\cdots (u+\beta^*_k+p_k)\ts\eta_{p_1,\dots,p_k}.
\tag 6.31
\endalign
$$
\endproclaim

\Proof Relation (6.31) follows from (6.17) and Proposition 6.3
by an easy induction; cf. the proof of (6.16). The proof of (6.30) 
is the same with the use of (6.12) and the expansion for 
$[\wh s_{-1,1}(u),s_{10}(v)]$ given by (1.28). $\square$
\medskip

Recall that the representation $L(\gamma,\gamma,\delta)$ of $\Y(3)$
if finite-dimensional if and only if $\gamma-\delta\in\ZZ_+$; see
Section 2.

\bigskip
\proclaim
{\bf Proposition 6.5} Let the irreducible representation $V(\mu(u))$
of $\Y^+(3)$ with the components of $\mu(u)$ given by (6.4)
be finite-dimensional and let $\wt{\mu}(u)$ be
defined by (6.4) with the corresponding parameters $\wt{\alpha}_i$ and
$\wt{\beta}_i$, where
$$
\wt{\alpha}_i=\alpha_i+l_i,\quad \wt{\beta}_i=\beta_i-m_i,
\qquad i=1,\dots,k
\tag 6.32
$$
and $l_i,m_i\in\ZZ_+$. Then the representation $V(\wt{\mu}(u))$
is also finite-dimensional.
\endproclaim

\Proof For each $i$
the representation 
$L(\alpha_i+l_i,\alpha_i+l_i,\alpha_i)\ot V(\mu(u))$
of $\Y^+(3)$ is finite-dimensional. The tensor product of the
highest weight vectors of 
$L(\alpha_i+l_i,\alpha_i+l_i,\alpha_i)$ and $V(\mu(u))$ generates
a $\Y^+(3)$-module with the highest weight $\mu'(u)$, where
$$
\aligned
\mu'_0(u)={}&\mu_0(u)(1-(\alpha_i+l_i)^2 u^{-2}),\\
\mu'_1(u)={}&\mu_1(u)(1-(\alpha_i+l_i)u^{-1})(1+\alpha_i u^{-1});
\endaligned
$$
see (1.31) and (3.11).
The composition of this module and the automorphism (1.29) with 
the series $\psi(u)=(1-\alpha_i^2 u^{-2})^{-1}$ yields a 
finite-dimensional representation
with the highest weight obtained from $\mu(u)$ by replacing
$\alpha_i$ with $\alpha_i+l_i$.

Similarly, considering tensor products of the form
$L(\beta_i,\beta_i,\beta_i-m_i)\ot V(\mu(u))$ we show that
the representation with the highest weight obtained from $\mu(u)$ by
replacing $\beta_i$ with $\beta_i-m_i$ if also finite-dimensional.
$\square$
\medskip

We are now in a position to establish necessary and sufficient conditions 
for $V(\mu(u))$ to be finite-dimensional.

\bigskip
\proclaim
{\bf Proposition 6.6} The irreducible highest weight representation
$V(\mu(u))$ of $\Y^+(3)$ with the components of $\mu(u)$ given by (6.4) is
finite-dimensional if and only if there exist re-enumerations of the
$\alpha_i$ and $\beta_i$ such that either one of the following two
conditions is satisfied
$$
\align
&\alpha_i-\beta_i\in\ZZ_+ \quad\text{for}\quad i=1,\dots,k; \tag 6.33\\
&\alpha_i-\beta_i\in\ZZ_+ \quad\text{for}\quad i=1,\dots,k-1 
\quad\text{and}\quad \alpha_k\in\frac12+\ZZ_+,\ \ \beta_k\in -\ZZ_+.
\tag 6.34
\endalign
$$
\endproclaim

\Proof Suppose that $\dim V(\mu(u))<\infty$. 
As before, we equip the vector space
$V=V(\mu(u))$ with one more structure of $\Y^+(3)$-module.
This
module is isomorphic to $V(\mu^*(u))$ with $\mu^*(u)$ given by (6.11).
Our basic idea is to construct nonzero vectors
which generate highest weight $\Y^+(2)$-submodules in $V(\mu^*(u))$
and then use the conditions on the highest weight 
imposed by Proposition 5.3.

Consider the mutliset
$
M=\{\alpha^*_1,\dots,\alpha^*_k,\beta^*_1,\dots,\beta^*_k\}.
$
Let us represent it as a disjoint union
$$
M=\underset{z\in\C/\ZZ}\to{\bigsqcup} M_z,\qquad M_z=(z+\ZZ)\cap M.
$$
Let
$$
M_z=\{\alpha^*_{i_1},\dots,\alpha^*_{i_r},
\beta^*_{j_1},\dots,\beta^*_{j_s}\}
\tag 6.35
$$
with both sequences 
written in the increasing order of the real parts, 
where $I_z=\{i_1,\dots,i_r\}$ and $J_z=\{j_1,\dots,j_s\}$ are subsets of 
$\{1,\dots,k\}$ depending on $z$.

Assume first that for each $M_z$ the following conditions are satisfied
$$
\aligned
&\alpha^*_{i_1}=\cdots=\alpha^*_{i_r}, \qquad
\beta^*_{j_1}=\cdots=\beta^*_{j_s},\qquad
\text{Re}\ts \beta^*_{j_1}<0,\\
&\text{and if}\quad \beta^*_{j'}\in M_{-z}\quad\text{then}\quad
\alpha^*_{i_1}+\beta^*_{j'}>0.
\endaligned
\tag 6.36
$$
Consider the vector $\eta_{p_1,\dots,p_k}$ (see (6.23)) with
the numbers $p_i$ defined as follows. For $M_z$ given by (6.35)
set
$p_{j_b}=\alpha_{j_b}+\alpha_{j'}$, $b=1,\dots,s$ if 
there exists $\beta^*_{j'}\in M_{-z}$; otherwise set
$p_{j_b}=0$. 
Starting with $\eta_{p_1,\dots,p_k}$
we shall apply repeatedly 
appropriate operators $S^{\natural}_{-1,0}(\alpha_i-q_i)$
to the vectors of the form $\eta_{q_1,\dots,q_k}$ to get
the highest weight vector 
$\xi=\eta_{0,\dots,0}$ with a certain coefficient.
By conditions (6.36) and the choice
of the numbers $p_i$ this coefficient is
nonzero, as follows from (6.24). This proves that 
$\eta_{p_1,\dots,p_k}\ne 0$.

On the other hand, by Proposition 6.4 the vector $\eta_{p_1,\dots,p_k}$
generates a $\Y^+(2)$-submodule in $V(\mu^*(u))$ with the highest weight
$$
(1-\alpha^*_1 u^{-1})\cdots (1-\alpha^*_k u^{-1})
(1+(\beta^*_1+p_1) u^{-1})\cdots (1+(\beta^*_k+p_k) u^{-1}).
\tag 6.37
$$
This submodule is finite-dimensional. However, for all $i$ we have
$-\beta^*_i-p_i-\frac12\not\in\ZZ_+$ and for all $i<j$ we have
$-\beta^*_i-p_i-\beta^*_j-p_j\not\in\ZZ_+$.
So, Proposition 5.3 implies that there exist
re-enumerations of the $\alpha^*_i$ and $\beta^*_i$ such that either
$$
\alpha^*_i-\beta^*_i\in\ZZ_+ \qquad\text{for}\quad i=1,\dots,k; \tag 6.38
$$
or
$$
\alpha^*_i-\beta^*_i\in\ZZ_+ \quad\text{for}\quad i=1,\dots,k-1; 
\qquad\text{and}\quad \alpha^*_k\in\frac12+\ZZ_+.\tag 6.39
$$

Let now the
parameters $\alpha^*_i$ and $\beta^*_i$ be arbitrary.
We can choose nonnegative integers
$l_i$ and $m_i$ in such a way that 
conditions (6.36) are satisfied for any $z$
by the shifted parameters
$\alpha^*_i+l_i$ and $\beta^*_i-m_i$. Therefore, 
by Proposition~6.5 either
(6.38), or (6.39) with $\ZZ_+$ replaced by $\ZZ$
will still hold. 

Let us consider the multisets $M_{z_0}$ for different $z_0$.

Suppose first that $z_0\not\equiv 0,\frac12\mod\ZZ$.
Then in (6.35)
$s$ may only be equal to $r$ or $r+1$.
Let us prove that
$$
\alpha^*_{i_t}-\beta^*_{j_t}\in\ZZ_+ \qquad\text{for}\quad t=1,\dots,r.
\tag 6.40
$$
Indeed, if this condition is violated then there exists an index
$m\in\{1,\dots,r\}$ such that the number of the elements $\beta^*_{j_t}$
satisfying $\alpha^*_{i_m}-\beta^*_{j_t}\in\ZZ_+$ is less than $m$.
Now we use Proposition 6.5 again. Applying shifts of the form (6.32),
if necessary, we may assume that conditions (6.36) are satisfied
for all $z\not\equiv z_0\mod\ZZ$ while for 
the elements of $M_{z_0}$ we have:
$$
\aligned
\alpha^*_{i_1}=\cdots=\alpha^*_{i_m}, &\qquad
\alpha^*_{i_{m+1}}=\cdots=\alpha^*_{i_r},\\
\beta^*_{j_1}=\cdots=\beta^*_{j_l}, &\qquad
\beta^*_{j_{l+1}}=\cdots=\beta^*_{j_s},
\endaligned
\tag 6.41
$$
with $\text{Re}\ts \beta^*_{j_1}<0$ and if $\beta^*_{j'}\in M_{-z_0}$
then $\alpha^*_{i_r}+\beta^*_{j'}>0$.
Here $l<m$ and the differences $\alpha^*_{i_1}-\beta^*_{j_1}$,
$\beta^*_{j_s}-\alpha^*_{i_1}$ and $\alpha^*_{i_r}-\beta^*_{j_s}$
are positive integers. 

Consider the vector $\eta_{p_1,\dots,p_k}$ with
the $p_i$ defined as in the previous argument except for
$i\in J_{z_0}\cup J_{-z_0}$. In the latter case
set $p_{j_{b}}=\alpha_{j_{b}}+\alpha_{j'}$ for
$b=l+1,\dots,s$ if there exists $\beta^*_{j'}\in M_{-z_0}$; otherwise
set $p_{j_b}=0$ (shifting further
$\beta^*_{j'}\to \beta^*_{j'}-m_{j'}$ we may suppose that
the $p_{j_{b}}$ are nonnegative integers).
Let each of the remaining $p_i$ with $i\in J_{z_0}\cup J_{-z_0}$
be equal to zero.
Using again (6.24) we find that
$\eta_{p_1,\dots,p_k}\ne 0$ and so by Proposition~6.4 it generates
a finite-dimensional $\Y^+(2)$-submodule in $V(\mu^*(u))$
with the highest weight given by (6.37). However,
one easily checks that the condition
of Proposition~5.3 is not satisfied. 
Contradiction.

Further, if $z_0\equiv 0\mod\ZZ$ the same argument shows that
(6.40) holds. The only difference is that here $M_{z_0}=M_{-z_0}$
and one should take $j'$ to be equal to $j_1$ or $j_s$ 
depending on whether $\beta^*_{j_s}\in 1+\ZZ_+$ or 
$\beta^*_{j_s}\in -\ZZ_+ $.

Finally, let $z_0\equiv \frac12\mod\ZZ$. In (6.35)
$s$ must be equal to $r$ or $r-1$. We shall prove that either
condition (6.40) holds (then $s=r$), or
$$
\alpha^*_{i'_t}-\beta^*_{j_t}\in\ZZ_+ \quad\text{for}\quad t=1,\dots,r-1
\qquad\text{and}\quad \alpha^*_{i'_r}\in\frac12+\ZZ_+
\tag 6.42
$$
for a permutation $(i'_1,\dots,i'_r)$ of the elements
of $I_{z_0}$.

Suppose that (6.40) is violated. Then we can again represent $M_{z_0}$
in the form (6.41). 
If either
$\beta^*_{j_s}\in \frac32 +\ZZ_+$ with 
$\alpha^*_{i_1}\in -\frac 12-\ZZ_+$,
or
$
\beta^*_{j_s}\in \frac12 -\ZZ_+ 
$
then we come to a contradiction exactly as above.
So, we must have
$
\beta^*_{j_s}\in \frac32 +\ZZ_+ ,\quad \alpha^*_{i_1}\in \frac 12+\ZZ_+.
$
Taking $j'=j_1$ and repeating the previous argument
we obtain the inequality $s-l\leq r-m$. 
It can only be satisfied if $s=r-1$ and $l=m-1$.
Note that this conclusion is valid for all possible values of 
the parameter $m$.
Let us now assume that $m$
takes the minimum value. Then (6.42) will be satisfied
for the indices $i'_1,\dots,i'_r$ defined as follows:
$$
\alignat2
i'_t={}&i_t\qquad&&\text{for}\quad t=1,\dots,m-1,\\
i'_t={}&i_{t+1}\qquad&&\text{for}\quad t=m,\dots,r-1,\\
i'_r={}&i_m. &&
\endalignat
$$

Thus, we have proved that if $V(\mu(u))$ is finite-dimensional then
there exist re-enumerations of the $\alpha_i$ and $\beta_i$ such that
either (6.38), or (6.39) holds.
Conditions (6.38) are equivalent to (6.33), while (6.39) means that
$$
\alpha_i-\beta_i\in\ZZ_+ \quad\text{for}\quad i=1,\dots,k-1; 
\qquad\text{and}\quad \beta_k\in -\ZZ_+;
\tag 6.43
$$
see (6.10).
Applying this result to the representation $V(\mu^*(u))$ 
instead of $V(\mu(u))$ we obtain that either
relation (6.33) holds, or there exist permutations
$\sigma,\tau\in\Sym_k$ such that
$$
\alpha_{\sigma(i)}-\beta_{\tau(i)}\in\ZZ_+ \quad\text{for}\quad
i=1,\dots,k-1;  \qquad\text{and}\quad \alpha_{\sigma(k)}\in \frac12+\ZZ_+.
\tag 6.44
$$
If $\sigma(k)=k$ then (6.34) is satisfied. Otherwise, set $l=\sigma(k)$
and rewrite (6.44) as follows:
$$
\alpha_i-\beta_{\rho(i)}\in\ZZ_+ \quad\text{for}\quad
i\ne l;  \qquad\text{and}\quad \alpha_l\in \frac12+\ZZ_+,
$$
where $\rho\in\Sym_k$. Consider the cycle of $\rho$ of the form
$k\to i_1\to\cdots\to i_s\to k$. Suppose that it does not contain $l$. Then
setting $\beta'_i=\beta_{\rho(i)}$ if $i$ belongs to the cycle,
and $\beta'_i=\beta_i$ otherwise, 
we obtain by (6.43) that
$\alpha_i-\beta'_i\in\ZZ_+$ for all $i=1,\dots,k$, that is, (6.33) holds.

Suppose now that the cycle contains $l$: $k\to i_1\to\cdots\to i_t\to l$.
Set
$\alpha'_k=\alpha_l$, $\alpha'_{\rho(i)}=\alpha_i$ for $i=k,i_1,\dots,i_t$
and $\alpha'_i=\alpha_i$ for the remaining indices $i$. 
Then by (6.43) relation (6.34) is satisfied for 
the $\alpha'_i$ and $\beta_i$
which completes the proof of the `only if' part of the proposition.
\medskip

Now let (6.33) hold. Then representation (6.5)
of $\Y^+(3)$ is finite-dimensional. However, $V(\mu(u))$ is isomorphic to
its subquotient and therefore $\dim V(\mu(u))<\infty$. 

Finally, let (6.34) hold. Consider the highest weight representation
$V(\mu)$ of the Lie algebra $\oa(3)$ and extend it to $\Y^+(3)$ by
using the homomorphism (1.23). We obtain a representation
with the highest weight given by (3.10).
By (3.9) the representation $V(\mu)$ is finite-dimensional
for $\mu=\beta_k-\frac12$. Therefore the representation of $\Y^+(3)$
of the form
$$
L(\alpha_1,\alpha_1,\beta_1)\ot\cdots \ot 
L(\alpha_{k-1},\alpha_{k-1},\beta_{k-1}) \ot L(\alpha_k,\alpha_k,\frac12)
\ot V(\beta_k-\frac12)
\tag 6.45
$$
is also finite-dimensional. However, we find from (1.31) and (3.11) that
$V(\mu(u))$ is isomorphic to a subquotient of (6.45) and so,
$\dime V(\mu(u))<\infty$. $\square$
\medskip

The following theorem together with Theorem 3.3
gives a description
of finite-dimensional irreducible representations of $\Y^+(2n+1)$
(cf. Theorems 4.8 and 5.9). We use notation (2.35).

\bigskip
\proclaim
{\bf Theorem 6.7} The irreducible highest weight representation
$V(\mu(u))$, $\mu(u)=(\mu_0(u),\dots,\mu_n(u))$ of 
$\Y^+(2n+1)$ is finite-dimensional if and only 
if either one of the following two relations
holds:
$$
\align
\mu_0&(u)\to \mu_1(u)\to\cdots\to\mu_n(u),\tag 6.46\\
\frac{2u}{2u+1}\ts\mu_0&(u)\to \mu_1(u)\to\cdots\to\mu_n(u).
\tag 6.47
\endalign
$$
\endproclaim

\Proof Suppose that $\dim V(\mu(u))<\infty$.
Then by Proposition 3.5 we have
$$
\mu_1(u)\to \mu_2(u)\to\cdots\to\mu_n(u).
$$
The subalgebra in $\Y^+(2n+1)$ 
generated by the coefficients of the series 
$s_{ij}(u)$ with $i,j=-1,0,1$
is isomorphic to
$\Y^+(3)$. The cyclic span of the highest weight vector of $V(\mu(u))$
with respect to this subalgebra is a representation with
the highest weight $(\mu_0(u),\mu_1(u))$. Its irreducible quotient
is finite-dimensional and so, by
Proposition 6.6 we have $\mu_0(u)\to \mu_1(u)$ or 
$\dfrac{2u}{2u+1}\ts\mu_0(u)\to \mu_1(u)$ depending on whether condition
(6.33) or (6.34) is satisfied; see the proof of Theorem 2.8.

Conversely, let relation (6.46) hold. Then 
for $i=1,\dots,n$ we have
$$
\frac{\mu_{i-1}(u)}{\mu_{i}(u)}=\frac{P_i(u+1)}{P_i(u)}
\tag 6.48
$$
for monic polynomials $P_i(u)$.
Let
$$
P_i(u)=(u+\delta_1^{(i)})\cdots (u+\delta_{s_i}^{(i)}),
\qquad i=1,\dots,n.
\tag 6.49
$$

Consider the irreducible highest weight representation
$L(\lambda(u))$ with $\lambda(u)=(\lambda_{-n}(u),\dots,\lambda_n(u))$
of the Yangian $\Y(2n+1)$ where
the $\lambda_i(u)$ are defined by (4.42)
for $i=1,\dots, n$ and by (4.43)
for $i=0,\dots,n$.
Then $L(\lambda(u))$ is finite-dimensional by Theorem 2.12
and we conclude by repeating 
the corresponding argument of the proof of Theorem 4.8.

Finally, let (6.47) hold. Then there exist monic polynomials
(6.49) such that (6.48) holds for $i=2,\dots,n$
and
$$
\frac{\mu_{0}(u)}{\mu_{1}(u)}=\frac{P_1(u+1)}{P_1(u)}
\frac{2u+1}{2u}. 
\tag 6.50
$$

Let $V(\mu_0)$ denote the irreducible representation of the Lie
algebra $\oa(2n+1)$ with the highest weight $\mu_0=(-1/2,\dots,-1/2)$. It
is
finite-dimensional; see (3.6) and (3.9). Extend $V(\mu_0)$ to
a representation of $\Y^+(2n+1)$ by using (1.23)
and consider the tensor product
$
L(\lambda(u))\ot V(\mu_0),
$
where $L(\lambda(u))$ is the defined above representation of $\Y(2n+1)$.
The tensor product of the
highest weight vectors of $L(\lambda(u))$ and $V(\mu_0)$ generates
a $\Y^+(2n+1)$-submodule with the highest weight
$\mu'(u)=(\mu'_0(u),\dots,\mu'_n(u))$ where
$\mu'_i(u)=\lambda_i(u)\lambda_{-i}(-u)(1+1/2 u^{-1})^{-1}$
for $i=1,\dots,n$ and $\mu'_0(u)=\lambda_0(u)\lambda_{0}(-u)$;
see (1.31) and (3.11). So the
representation $V(\mu'(u))$ of $\Y^+(2n+1)$ is finite-dimensional
and the $\mu'_i(u)$ satisfy (6.48) 
for $i=2,\dots, n$ and (6.50).
It remains to repeat
the corresponding argument of the proof of Theorem 5.9. $\square$
\medskip

Theorem 6.7 implies the following parameterization
of representations of the
special twisted Yangian $\SY^+(2n+1)$
(cf. Corollaries 4.9 and 5.10).

\bigskip
\proclaim
{\bf Corollary 6.8} There is a one-to-one correspondence between
finite-dimensional irreducible representations
of the special twisted Yangian $\SY^+(2n+1)$ and the families
$\{P_1(u),\dots,P_{n}(u),\varepsilon\}$, where the $P_i(u)$ are monic
polynomials in $u$ and the parameter $\varepsilon$ takes values in
$\{1,2\}$.
Every such representation is
isomorphic to a subquotient of a representation 
of the form $L(\lambda(u))$ or $L(\lambda(u))\ot V(\mu_0)$. $\square$
\endproclaim

\bigskip
\noindent
{\it Remark 6.9.} One can use Propositions 6.3 and 6.4 to get
a `branching rule' for the restriction of a `generic' finite-dimensional
irreducible representation $V$ of $\Y^+(3)$ to the subalgebra $\Y^+(2)$
and thus to construct an analogue of Gelfand-Tsetlin basis in $V$;
cf. [M2], [NT2]. $\square$

\bigskip
\bigskip
\noindent
{\bf References}
\bigskip

\itemitem{[BL]} {D. Bernard and A. LeClair},
{\sl Quantum group symmetries and non-local currents
in 2D QFT}. {Commun. Math. Phys.} {\bf 142} (1991), 99--138.

\itemitem{[CP1]}
{V. Chari and A. Pressley},
{\sl Yangians and $R$-matrices},
{L'Enseign. Math.}
{\bf 36}
(1990),
267--302.

\itemitem{[CP2]}
{V. Chari and A. Pressley},
{\sl Fundamental representations of Yangians and
singularities of $R$-matrices},
{J. Reine Angew. Math.}
{\bf 417}
(1991),
87--128.

\itemitem{[CP3]}
{V. Chari and A. Pressley},
{\sl Yangians: their representations and characters,}
{Acta Appl. Math.} {\bf 44} (1996), 39--58.

\itemitem{[C1]}
{I. V. Cherednik},
{\sl A new interpretation of Gelfand--Tzetlin bases}, {Duke Math. J.}
{\bf 54}
(1987),
563--577.

\itemitem{[C2]}
{I. V. Cherednik},
{\sl Quantum groups as hidden symmetries of classic
representation theory},
in \lq Differential Geometric Methods in Physics
(A. I. Solomon, Ed.)',
World Scientific,
Singapore,
1989,
pp. 47--54.

\itemitem{[Di]}
{J. Dixmier},
{\sl Alg\`ebres Enveloppantes}, 
{Gauthier-Villars, Paris},
1974.

\itemitem{[D1]}
{V. G. Drinfeld},
{\sl Hopf algebras and the quantum Yang--Baxter equation}, {Soviet Math.
Dokl.}
{\bf 32}
(1985),
254--258.

\itemitem{[D2]}
{V. G. Drinfeld},
{\sl A new realization of Yangians and quantized affine
algebras},
{Soviet Math. Dokl.}
{\bf 36}
(1988),
212--216.

\itemitem{[KR]}
{A. N. Kirillov and N. Yu. Reshetikhin},
{\sl Yangians, Bethe ansatz and combinatorics},
{Lett. Math. Phys.}
{\bf 12}
(1986),
199--208.

\itemitem{[KK]} {T. H. Koornwinder and V. B. Kuznetsov},
{\sl Gauss hypergeometric function and quadratic $R$-matrix algebras},
{St. Petersburg Math.~J.} {\bf 6} (1994), 161--184.

\itemitem{[K]} {V. E. Korepin}, {\sl The analysis of 
the bilinear relation of
the six-vertex model}, {Soviet Phys. Dokl.} {\bf 27} (1982), 612--613.

\itemitem{[KS1]}
{P. P. Kulish and E. K. Sklyanin},
{\sl Quantum spectral transform method: recent developments},
in \lq Integrable Quantum Field Theories', {Lecture Notes in Phys.}
{\bf 151}
Springer,
Berlin-Heidelberg,
1982,
pp. 61--119.

\itemitem{[KS2]}
{P. P. Kulish and E. K. Sklyanin},
{\sl Algebraic structures related to reflection equations}, {J. Phys.}
{\bf A25}
(1992),
5963--5975.

\itemitem{[KJC]}
{V. B. Kuznetsov, M. F. J\o{rgensen}, P. L. Christiansen},
{\sl New boundary conditions for integrable lattices},
{J. Phys. A} {\bf 28} (1995), 4639.

\itemitem{[M1]}
{A. Molev},
{\sl Representations of twisted Yangians}, {Lett. Math. Phys.}
{\bf 26}
(1992),
211--218.

\itemitem{[M2]}
{A. Molev},
{\sl Gelfand--Tsetlin basis for representations of Yangians},
{Lett. Math. Phys.}
{\bf 30}
(1994),
53--60.

\itemitem{[M3]}
{A. Molev},
{\sl Sklyanin determinant, Laplace operators, and characteristic identities
for classical Lie algebras}, 
{J. Math. Phys.}
{\bf 36}
(1995),
923--943.

\itemitem{[M4]} {A. Molev}, {\sl Noncommutative symmetric 
functions and Laplace operators
for classical Lie algebras}, {Lett. Math. Phys.}
{\bf 35} (1995), 135-143.

\itemitem{[MN]} {A. Molev and M. Nazarov},
{\sl Capelli identities for classical Lie algebras},
Preprint CMA 003-97, Australian National University, Canberra.

\itemitem{[MNO]}
{A. Molev, M. Nazarov and G. Olshanski},
{\sl Yangians and classical Lie algebras}, 
Russian Math. Surveys
{\bf 51}:2
(1996),
205--282.

\itemitem{[N1]}
{M. Nazarov},
{\sl Quantum Berezinian and the classical Capelli identity},
{Lett. Math. Phys.}
{\bf 21}
(1991),
123--131.

\itemitem{[N2]}
{M. Nazarov},
{\sl Yangians and Capelli identities},
in
{``A.~A.~Kirillov Seminar on Representation Theory"}, 
%(S.~Gindikin, Ed.) 
%{Amer. Math. Soc. Transl.}
%{\bf ???},
AMS,
Providence,
1997, pp 139--164;
q-alg/9601027.

\itemitem{[NO]}
{M. Nazarov and G. Olshanski},
{\sl Bethe subalgebras in twisted Yangians}, 
Commun. Math. Phys.
{\bf 178}
(1996),
483--506.

\itemitem{[NT1]}
{M. Nazarov and V. Tarasov},
{\sl Yangians and Gelfand--Zetlin bases}, Proc. RIMS, Kyoto Univ.,
{\bf 30} (1994), 459--478.

\itemitem{[NT2]}
{M. Nazarov and V. Tarasov},
{\sl Representations of Yangians with Gelfand--Zetlin bases}, 
to appear in {J. Reine Angew. Math.}; q-alg/9502008.

\itemitem{[Ok]}
{A. Okounkov},
{\sl Quantum immanants and higher Capelli identities},
{Transformation Groups}
{\bf 1}
(1996),
99--126.

\itemitem{[O1]}
{G. Olshanski},
{\sl Representations of infinite-dimensional classical
groups, limits of enveloping algebras, and Yangians},
in \lq Topics in Representation Theory (A. A. Kirillov, Ed.), {Advances in
Soviet Math.} {\bf 2},
AMS,
Providence RI,
1991,
pp. 1--66.

\itemitem{[O2]}
{G. Olshanski},
{\sl Twisted Yangians and infinite-dimensional classical Lie algebras},
in \lq Quantum Groups (P. P. Kulish, Ed.)', {Lecture Notes in Math.}
{\bf 1510},
Springer,
Berlin-Heidelberg,
1992,
pp. 103--120.

\itemitem{[S]}
{E. K. Sklyanin},
{\sl Boundary conditions for integrable quantum systems}, {J. Phys.}
{\bf A21}
(1988),
2375--2389.

\itemitem{[TU]}
{K. Takemura and D. Uglov}, {\sl The orthogonal eigenbasis and norms
of eigenvectors in the Spin Calogero-Sutherland Model}, Preprint RIMS 1114,
Kyoto, 1996; solv-int/9611006.

\itemitem{[TF]}
{L. A. Takhtajan and L. D. Faddeev},
{\sl Quantum inverse scattering method and the Heisenberg
$XYZ$-model},
{Russian Math. Surv.}
{\bf 34}
(1979),
no. 5,
11--68.

\itemitem{[T1]}
{V. O. Tarasov},
{\sl Structure of quantum $L$-operators for the
$R$-matrix of of the $XXZ$-model},
{Theor. Math. Phys.}
{\bf 61}
(1984),
1065--1071.

\itemitem{[T2]}
{V. O. Tarasov},
{\sl Irreducible monodromy matrices for the $R$-matrix of the
$XXZ$-model and lattice local quantum Hamiltonians}, {Theor. Math. Phys.}
{\bf 63}
(1985),
440--454.

\itemitem{[UK]}
{D. B. Uglov and V. E. Korepin}, {\sl The Yangian symmetry of
the Hubbard model}, {Phys. Lett. A} {\bf 190} (1994), 238--242.

\itemitem{[Z1]}
{R. B. Zhang}, {\sl Representations of super Yangian},
{J. Math. Phys.} {\bf 36} (1995), no. 7, 3854--3865.

\itemitem{[Z2]}
{R. B. Zhang}, {\sl The $\gl(M|N)$ super Yangian and
its finite-dimensional representations}, {Lett. Math. Phys.} 
{\bf 37} (1996), 419--434.

\enddocument